\documentclass[a4paper,12pt]{article}
\usepackage{amssymb}
\usepackage{latexsym}
\usepackage{amsmath}
\usepackage{amsthm}
\usepackage{amsfonts}
\usepackage{amssymb}
\usepackage{graphicx}

\newtheorem{Th}{Theorem}
\newtheorem{Prop}{Proposition}

\newtheorem{Lm}{Lemma}
\newtheorem{Lma}{Lemma}[section]
\newtheorem{Dfi}{Definition}

\newcommand{\be}{\begin{equation}}
\newcommand{\ee}{\end{equation}}
\newcommand{\R}{\mathbb{R}}
\newcommand{\N}{\mathbb{N}}
\newcommand{\C}{\mathbb{C}}
\newcommand{\Z}{\mathbb{Z}}

\newcommand\res{\mathop{\hbox{\vrule height 7pt width .5pt depth 0pt
\vrule height .5pt width 6pt depth 0pt}}\nolimits}

\newcommand{\reset}{\setcounter{equation}{0}\setcounter{Th}{0}\setcounter{Prop}{0}\setcounter{Co}{0}
\setcounter{Lm}{0}\setcounter{Rm}{0}}

\def\lf{\left}
\def\rg{\right}

\def\al{\alpha}
\def\la{\lambda}

\def\ep{\varepsilon}
\def\ds{\displaystyle}
\def\ov{\overline}
\def\Om{\Omega}
\def\om{\omega}
\def\p{\partial}
\def\bn{\vec{n}}

\def\bbe{\vec{e}}

\def\bw{\vec{w}}

\def\bP{\vec{\Phi}}

\def\bP{\vec{\Phi}}

\def\scII{\mathbb I}

\def\res{\mathop{\hbox{\vrule height 7pt width .5pt 
depth 0pt\vrule height .5pt width 6pt depth 0pt}}\nolimits}

\DeclareMathOperator{\supp}{supp}

\begin{document}
\title{Tori in $S^3$ minimizing locally the conformal volume.}
\author{ Tristan Rivi\`ere\footnote{Forschungsinstitut f\"ur Mathematik, ETH Zentrum,
CH-8093 Z\"urich, Switzerland.}}
\date{ }
\maketitle

{\bf Abstract :}{\it We prove that the conformal immersions of  complex two tori into $S^3$ which locally minimize their conformal volume in their conformal class
all satisfy some elliptic PDE. We prove that they are either minimal tori, CMC flat tori,  elliptic conformally constrained minimal tori or critical point of the area under
some fixed conformally congruent area. On the way to establish this result we prove that tori which are critical points of the area for perturbations within a given conformal class
and which are degenerate points of the conformal class mapping - i.e. isothermic - are either minimal surfaces or flat CMC tori. These results are all proved in the general framework of weak immersions.}
\medskip

\noindent{\bf Math. Class. 49Q10, 53A05, 53A30, 35J20}
\section{Introduction}

The notion of conformal volume introduced by P.Li and S.T.Yau in \cite{LiYa} has stimulated a broad interest in mathematics beyond the differential geometry of submanifolds.

Let $G(S^3)$ be the M\"obius group of conformal transformations of $S^3$ we are considering immersions $\vec{\Phi}$ of the two torus $T^2$ which are critical
for the conformal volume
\be
\label{0-0}
V_c(\vec{\Phi}):=\sup_{\Psi\in G(S^3)} A(\Psi\circ\vec{\Phi})
\ee
within  their conformal class. Precisely for any smooth path $\vec{\Phi}_t$ such that $\vec{\Phi}_0=\vec{\Phi}$ and such that the conformal class defined by
the metric $\vec{\Phi}_t^\ast g_{S^3}$, where $g_{S^3}$ denotes the standard metric on the $3-$sphere $S^3$, is equal to the one defined by $\vec{\Phi}^\ast g$
we have
\be
\label{0-1}
\lf.\frac{d}{dt}V_c(\vec{\Phi}_t)\rg|_{t=0}=0
\ee
Absolute minimizers of $V_c$ in a given family of conformally equivalent metrics, when they exist, do satisfy (\ref{0-1}) for instance. Beside the case of rectangular tori,
which is partly solved in \cite{MR}, it is not known which conformal class possess a minimizer for the conformal volume. It is expected that not every class posses
a minimizer of $V_c$. Indeed, while taking a minimizing sequence $\vec{\Phi}_k$ and the optimal $\Psi_k$, if it exist,
satisfying $V_c(\vec{\Phi}_k):= A(\Psi_k\circ\vec{\Phi}_k)$ one cannot a-priori exclude that the sequence of immersions $\Psi_k\circ\vec{\Phi}_k$
 weakly converges to a geodesic sphere. This degeneracy in the minimization process is believed by the author to be the only possible one and, when the minimal
  conformal volume within a conformal class is strictly larger than $4\pi$, it should be achieved by an immersion of $T^2$ satisfying (\ref{0-1}). 
  
The purpose of the present work is to give a characterization of 2-dimensional tori in $S^3$ which are local minimizers of the conformal volume in their conformal class. In other words we aim
to give a characterization of complex 2-tori (tori equipped with a given complex structure) minimizing locally the conformal volume for immersions defining this given complex structure.

In \cite{Ri4} the author derived the infinitesimal characterization of the variations of riemann surfaces :   a PDE satisfied by critical points of functional such as the area for variations  within a given conformal class was derived
in this work. This was a delicate issue
at least for isothermic surfaces which are the degenerate points of the conformal class mapping and for which a direct application of Lagrange multiplier theory was excluded a-priori. It is proved that, for any $\vec{\Phi}$ satisfying (\ref{0-1}) for any variation $\vec{\Phi}_t$ within the conformal class of $\vec{\Phi}$, there exists an holomorphic quadratic form
$Q$ of the underlying riemann surface such that
\be
\label{0-2}
H=\Re<Q,h^0>_{g_{\vec{\Phi}}}\quad,
\ee
where $g_{\vec{\Phi}}:=\vec{\Phi}^\ast g_{S^3}$ is the induced metric on the surface, $H$ is the mean curvature of the immersion and $h^0$ is the Weingarten form given in complex coordinates by
\be
\label{0-3}
h^0:=2\,\pi_{\vec{n}}(\p^2_{z^2}\vec{\Phi})\ dz^2\quad,
\ee
where $\pi_{\vec{n}}$ is the projection onto the normal direction to the surface in $TS^3$. Solutions to (\ref{0-2}) are called {\it conformally constrained minimal surfaces}.

Hopf tori - i.e. the preimages by the Hopf fibration from $S^3$ into $S^2$ of closed curves on $S^2$ - are examples of {\it conformally constrained minimal surfaces} (see \cite{BPP}). These surfaces are not necessarily smooth and it suffices to have a weak notion of second fundamental form as well as a notion of conformal coordinates  in order
to give a meaning to equation (\ref{0-2}). In \cite{Ri3} (see also \cite{Ri1}) the author introduced a notion of weak immersions compatible with variational purposes, for which conformal coordinates exist and for which the mapping assigning
the conformal class is smooth for some topology.  The motivation originally was to produce a suitable framework for studying the variations
of Willmore Lagrangian. It appears moreover that this is a ''minimal'' requirement for ensuring the smoothness of the mapping assigning the conformal class. This is the framework we shall consider here. 

 A {\it weak immersion} 
of $T^2$ into $S^3$ is a map
$\vec{\Phi}$ from $T^2$ into $S^3$ such that
\begin{itemize}
\item[i)]
\[
\vec{\Phi}\in W^{1,\infty}(\Sigma,{\R}^4)\quad,
\]
\item[ii)] there exists a constant $C_{\vec{\Phi}}>1$ such that
\[
\forall x\in \Sigma\quad\forall X\in T_x\Sigma\quad\quad C_{\vec{\Phi}}^{-1}\ g_0(X)\le |d\vec{\Phi}(X)|^2\le  C_{\vec{\Phi}}\ g_0(X)
\]
where $g_0$ is some given smooth reference metric on $T^2$ i.e. in other words the metric on $T^2$ equal to the pull back by $\vec{\Phi}$ of the canonical metric of ${\R}^4$ is equivalent
to any reference metric on $T^2$,
\item[iii)]
\[
\vec{n}_{\vec{\Phi}}\in W^{1,2}(T^2,S^3)
\]
where $\vec{n}_{\vec{\Phi}}$ is the Gauss map associated to $\vec{\Phi}$ i.e. the unit vector perpendicular to the surface in $TS^3$ and positively oriented. 
\end{itemize}

This space is usually denoted ${\mathcal E}_{T^2}$. This is a Banach manifold modeled on $W^{1,\infty}\cap W^{2,2}(T^2,{\R}^3)$. Each element in ${\mathcal E}_{T^2}$ defines uniquely 
a conformal class and, for a fixed choice of generators of the $\pi_1(T^2)$ the mapping which to $\vec{\Phi}$ assigns the corresponding Teichm\"uller class is $C^1$ (see \cite{Ri3}).

%A further motivation for considering the space of {\it weak immersion} comes from the fact that the area is twice differentiable in this space and this is the lowest requirement
%for this fact to be true since the second derivative involves the $L^2$ norm of the second fundamental form.

The analysis of the equation (\ref{0-2}) is made particularly complex due to the fact that it includes both {\it hyperbolic} and {\it elliptic} regimes. {\it Strictly  elliptic }regime is observed
in the domain given by
\[
{\mathfrak E}_s(\vec{\Phi}):=\lf\{x\in T^2\ ;\ 2\,|Q|_{g_{\vec{\Phi}}}(x) <1 \rg\}
\]
The splitting between these two regimes can be seen formally in the following way. On a two torus an holomorphic quadratic form is either zero or never vanish. Then, locally there exist
complex coordinates in which $Q=4^{-1}\,dz^2$. In these coordinates the equation (\ref{0-2}) reads
\[
\Delta\vec{\Phi}+2\,\vec{\Phi}\, e^{2\la}=\p_{x_1}(e^{-2\la}\,\p_{x_1}\vec{\Phi})-\p_{x_2}(e^{-2\la}\,\p_{x_2}\vec{\Phi})
\]
where $g_{\vec{\Phi}}=e^{2\la}\ [dx_1^2+dx_2^2]$. Observe that  $|Q|_{g_{\vec{\Phi}}}= 2^{-1}\, e^{-2\la}$ hence the domain of strict ellipticity ${\mathfrak E}_s(\vec{\Phi})$ 
corresponds to the set of points $x$ where $\la(x)>0$ in these special coordinates and hence the principal symbol of the {\it conformally constrained minimal surfaces} equation  which reads
\be
\label{0-4}
\p_{x_1}\lf((1-e^{-2\la})\, \p_{x_1}\vec{\Phi}\rg)+\p_{x_2}\lf((1+e^{-2\la})\, \p_{x_2}\vec{\Phi}\rg)+2\,\vec{\Phi}\, e^{2\la}=0
\ee
is exactly invertible on this domain. Following \cite{MS} and \cite{He}, we proved in \cite{Ri1} that for any weak immersion $\vec{\Phi}$ in any conformal coordinates, the induced metric is continuous and hence ${\mathfrak E}_s(\vec{\Phi})$ is an open subset
of $T^2$. The elliptic nature of the {\it conformally constrained minimal surfaces} equation on the domain of strict ellipticity is reinforced by the following result
\begin{Th}
\label{th-0.1}
Let $\vec{\Phi}$ be a weak immersion satisfying the conformally constrained minimal surface equation
\[
H=\Re<Q,h^0>_{g_{\vec{\Phi}}}
\]
for some holomorphic quadratic form on $T^2$ equipped with the conformal class defined by $\vec{\Phi}$. Then $\vec{\Phi}$ is analytic in conformal coordinates
within the domain of strict ellipticity ${\mathfrak E}_s(\vec{\Phi})$ given by
\[
{\mathfrak E}_s(\vec{\Phi}):=\lf\{x\in T^2\ ;\ 2\,|Q|_{g_{\vec{\Phi}}}(x) <1 \rg\}\quad.
\]
\hfill $\Box$
\end{Th}
The theorem is clearly optimal since Hopf tori, which are solutions to the {\it conformally constrained minimal surface} equation, satisfy $2\,|Q|_{g_{\vec{\Phi}}}(x) \equiv 1$ on $T^2$ and can be nowhere
$C^{2}$ by taking the lifting of a nowhere $C^{2}$ curve $\Gamma$ on $S^2$ but still satisfying $\int_\Gamma|\kappa|^2\ dl<+\infty$ where $\kappa$ is the geodesic curvature of the curve in $S^2$
which ensures that this is a weak immersion.
\begin{Dfi}
\label{df-0-1}
A weak immersion is said to satisfy the elliptic conformally constrained minimal surface equation (resp. strictly elliptic) if it satisfies
\[
H=\Re<Q,h^0>_{g_{\vec{\Phi}}}
\]
for some holomorphic quadratic form on $T^2$ equipped with the conformal class defined by $\vec{\Phi}$ and 
\[
{\mathfrak E}(\vec{\Phi}):=\lf\{x\in T^2\ ;\ 2\,|Q|_{g_{\vec{\Phi}}}(x) \le 1 \rg\}
\]
concide with the whole torus (resp. ${\mathfrak E}_s(\vec{\Phi})=T^2$).\hfill $\Box$
\end{Dfi}
Strictly elliptic {conformally constrained minimal surface } could be seen as generalization of minimal surfaces sharing many features with minimal surfaces. In particular one deduces from the analysis for proving theorem~\ref{th-0.1} that the space
of strictly elliptic {conformally constrained minimal surface } satisfying $2\,|Q|_{g_{\vec{\Phi}}}(x) < 1-\ep_0$ for any $0< \ep_0\le1$ with uniformly bounded area and controlled conformal class
is compact in $C^l$ topology for any $l\in {\N}$.

Regarding now the similar problem in the euclidian space ${\R}^3$, in \cite{Ri4} it is proved that critical points of the area among weak immersions of a compact surface $\Sigma$ into ${\R}^3$ with prescribed conformal class satisfies also the {\it conformally constrained minimal surface } equation 
\[
\vec{H}=\Re<Q,\vec{h^0}>_{g_{\vec{\Phi}}}
\]
By multiplying the equation by $\vec{\Phi}$ and integrating by parts one gets the following proposition which is reminiscent to the corresponding result for minimal surfaces.
\begin{Prop}
\label{pr-I-1}
There exists no compact strictly elliptic conformally constrained minimal surface  in ${\R}^3$.\hfill $\Box$
\end{Prop}
 Observe that an arbitrary cylinder over a plane curve in ${\R}^3$ is an example of an elliptic (non strictly elliptic) non compact conformally constrained minimal surface satisfying
 $H=\Re(H^0)$ in cylindrical coordinates in which the Weingarten Operator takes the form 
 \[
{\mathbb I}=\ \lf(
 \begin{array}{cc}
 \kappa & 0\\[3mm]
 0& 0
 \end{array}
 \rg)\quad.
 \]
 It is also {\it isothermic} since $\Im(H^0)\equiv 0$. In \cite{BPP} the authors show that surfaces with rotational symmetry are {\it conformally constrained minimal} everywhere away from the axis points. 
 Such surfaces are not necessarily smooth and therefore should not be in general {\it elliptic conformally constrained minimal}.
 
 Going back to immersions in $S^3$, in section IV of the present paper, we establish the following result which identifies the space of {\it isothermic conformally constrained minimal} immersions of tori.
 \begin{Th}
 \label{th-III.1}
 A weak immersion of the torus $T^2$ in $S^3$ is both {\it isothermic} and a critical point of the area under constrained conformal class  if and only
 if it is either minimal or realizes a flat CMC torus. Such tori  are  all solving the strictly elliptic conformally constrained minimal equation. \hfill $\Box$
 \end{Th}
 For instance, the {\it Weingarten Operator} of {\it Hopf tori} in coordinates given by the Hopf fibers and parallel lifts of the curve $\Gamma$ in $S^2$ - whose lift by the {\it Hopf fibration} is equal to the torus - is given by - see \cite{Pi} identity (21) -
  \[
{\mathbb I}=\  \lf(
 \begin{array}{cc}
 2\,\kappa & 1\\[3mm]
 1& 0
 \end{array}
 \rg)\quad.
 \]
It is clearly solving the conformally constrained minimal surface equation $H=\Re(H^0)$. It is moreover isothermic if and only if $\kappa=\kappa_0\ge 0$ is constant - i.e. $\Gamma$ is a circle in $S^2$ and the torus is a CMC Clifford torus -. It solves the {\it isothermic} equation $$\Re((1+i\,\kappa_0)H^0)=0\quad,$$  together with the {\it strictly elliptic conformally constrained minimal} equation 
$$
H=\Re\lf(\frac{\kappa_0}{\kappa_0+i}\,H^0\rg)\quad,
$$
where $\la\equiv 0$ and $Q=4^{-1}\frac{\kappa_0}{\kappa_0-i} dz^2$ satisfies $2|Q|_{g}=\sqrt{\frac{\kappa_0^2}{1+\kappa_0^2}}<1$. Finally theorem~\ref{th-III.1} can be put in perspective with a result by 
J. Richter, see \cite{Ric} and \cite{BuPP}, asserting that {\it conformally constrained Willmore } tori in $S^3$ which are in addition {\it isothermic} are CMC surfaces for some constant sectional curvature.

 \medskip
 
 The main result of the present paper is the following theorem which identifies tori minimizing locally the conformal volume within their conformal class.
\begin{Th}
\label{th-0-2}
Let $\vec{\Phi}$ be weak immersion locally minimizing the conformal volume $V_c$ in it's conformal class. Assume $V_c(\vec{\Phi})\notin 4\pi{\N}$ and that $V_c$ is differentiable at $\vec{\Phi}$, then the following alternative holds :
either 
\begin{itemize}
\item{i)} $\vec{\Phi}$ is a minimal immersion,
\item{ii)} $\vec{\Phi}(T^2)$ is a flat CMC clifford torus congruent to $a\, S^1\times\sqrt{1-a^2}\, S^1$ for some $a\in (0,1)$ ,
\item{iii)} $\vec{\Phi}$ is a critical point of the area under the constraint that $A(\Psi\circ\vec{\Phi})$ is constant for some conformal transformation $\Psi$ of $S^3$ which is not
an isometry
\item{iv)} $\vec{\Phi}$ is an elliptic conformally constrained minimal surface.
\end{itemize}
\hfill$\Box$
\end{Th}

The paper is organized as follows. In section II we compute the first and second variations of the area functional within a fixed conformal class. In section III we prove theorem~\ref{th-0.1}. In section IV we identify the isothermic conformally constrained minimal surfaces. In section V we construct families of deformations which reduces the conformal volume for non-elliptic conformally constrained minimal surfaces. In the last section we collect the informations obtained from the previous ones to prove theorem~\ref{th-0-2}.

\section{The first and second variations of a weak immersion locally minimizing the conformal volume in $S^3$.}

\subsection{Preliminaries : notations and some computations.} 
Consider $\vec{\Phi}_t:=\vec{\Phi}+t\,\vec{w}$ .

\medskip

We have in local coordinates (that we can choose to be conformal for $\vec{\Phi}$)
\be
\label{a-III.4}
\vec{n}_t=\star_{{\R}^4}\lf(\vec{\Phi}_t\wedge\frac{\p_{x_1}\vec{\Phi}_t\wedge\p_{x_2}\vec{\Phi}_t}{|\p_{x_1}\vec{\Phi}_t\wedge\p_{x_2}\vec{\Phi}_t|}\rg)
\ee
We have
\be
\label{a-III.4a}
\vec{n}_t=\vec{n}+t\ (a_1\,\vec{e}_1+a_2\,\vec{e}_2+b\,\vec{\Phi}) +o(t)
\ee
They can be identified as follows:
\begin{align*}
	e^\lambda\, a_1&=\left\langle e^\lambda\,\bbe_1,\left.\frac{d}{dt}\vec{n}_t\right|_{t=0}\right\rangle =-\left\langle \left.\frac{d}{dt}\left(\p_{x_1}\bP_t\right)\right|_{t=0},\bn\right\rangle=-\langle\p_{x_1} \bw,\bn\rangle,
\end{align*}
and similarly one obtains $a_2=-e^{-\lambda}\,\langle\p_{x_2} \bw,\bn\rangle$. Let $v:=\frac{d\vec{\Phi}}{dt}\cdot\vec{n}=\vec{w}\cdot\vec{n}$. We have
\[
b:=\vec{\Phi}\cdot\frac{d\vec{n}}{dt}=-v
\]
Hence we have
\be
\label{vec-1}
\frac{d\vec{n}}{dt}=-e^{-\la}\lf[\langle\p_{x_1} \bw,\bn\rangle\ \vec{e}_1+\langle\p_{x_2} \bw,\bn\rangle\ \vec{e}_2\rg]-\, v\ \vec{\Phi}
\ee
We now compute the time derivative of the second fundamental form. We have
\be
\label{vec-2}
\begin{array}{l}
\ds\frac{d{\mathbb I}}{dt}=\sum_{i,j=1}^2\lf[\frac{d\vec{n}}{dt}\cdot\partial^2_{x_i\,x_j}\vec{\Phi}+\vec{n}\cdot\p^2_{x_i\,x_j}\vec{w}\rg]\ dx_i\otimes dx_j\\[5mm]
\ds\quad=-\sum_{i,j,k=1}^2\lf[e^{-2\la}\,\p_{x_k}\vec{w}\cdot\vec{n}\ \p_{x_k}\vec{\Phi}\cdot\p^2_{x_i x_j}\vec{\Phi}\rg]\ dx_i\otimes dx_j\\[5mm]
\ds\quad+v\ g+\sum_{i,j=1}^2\vec{n}\cdot\p^2_{x_i\,x_j}\vec{w}\ dx_i\otimes dx_j
\end{array}
\ee
Recall that
\be
\label{vec-3}
e^{-2\la}\ \p_{x_k}\vec{\Phi}\cdot\p^2_{x_i x_j}\vec{\Phi}=\delta_{jk}\ \p_{x_i}\la+\delta_{ik}\ \p_{x_j}\la-\delta_{ij}\ \p_{x_k}\la
\ee
Combining (\ref{vec-2}) and (\ref{vec-3}) we obtain
\be
\label{vec-4}
\begin{array}{l}
\ds\frac{d{\mathbb I}}{dt}=-\,(\p_{x_1}\la\,\p_{x_1}\vec{w}\cdot\vec{n}-\p_{x_2}\la\,\p_{x_2}\vec{w}\cdot\vec{n})\ (dx_1^2-dx_2^2)\\[5mm]
\ds\quad-\,(\p_{x_1}\la\,\p_{x_2}\vec{w}\cdot\vec{n}+\p_{x_2}\la\,\p_{x_1}\vec{w}\cdot\vec{n})\ (dx_1\otimes dx_2+dx_2\otimes dx_1)\\[5mm]
\ds\quad+v\ g+\sum_{i,j=1}^2\vec{n}\cdot\p^2_{x_i\,x_j}\vec{w}\ dx_i\otimes dx_j
\end{array}
\ee

Recall  that
$
 H_t=\frac12 \sum_{i,j} (g_t)^{ij}({\scII}_t)_{ij}
$
Hence 
\[
\frac{dH}{dt}=\frac12 \sum_{i,j} (g_t)^{ij}\frac{d({\scII}_t)_{ij}}{dt}+\frac12 \sum_{i,j} \frac{d(g_t)^{ij}}{dt}({\scII}_t)_{ij}
\]
Using (\ref{vec-4}) we have
\be
\label{vec4-b-1}
\frac12 \sum_{i,j} (g_t)^{ij}\frac{d({\scII}_t)_{ij}}{dt}=v+\frac{\vec{n}}{2}\cdot\Delta_g\vec{w}
\ee
%where we denote $(\scII_t)_{ij}=\langle\II_t(\p_{x_i},\p_{x_j}),\bn_t\rangle=-\langle\p_{x_i}\bn_t,\p_{x_j}\vec{\Phi}\rangle$. 
%\begin{multline}
%\label{0III.10i}
%\left.\frac{d}{dt}H_t\right|_{t=0}=-\frac12 \sum_{i,j}\left.\frac{d}{dt}(g_t)^{ij}\right|_{t=0}\,\langle\p_{x_i}\bn,\p_{x_j}\vec{\Phi}\rangle\\
%+g^{ij}\,\left[ \left<\p_{x_i}\left.\frac{d}{dt}\vec{n}_t\right|_{t=0},\p_{x_j}\vec{\Phi}\right>+\langle\p_{x_i}\vec{n},\p_{x_j}\vec{w}\rangle \right].
%\end{multline}
We have $(g_t)_{ij}=\langle\p_{x_i}\vec{\Phi}_t,\p_{x_j}\vec{\Phi}_t\rangle$, thus
\begin{equation}
\label{0III.10f}
\begin{split}
\left.\frac{d}{dt}(g_t)_{ij}\right|_{t=0}&=\langle \p_{x_i}\vec{w},\p_{x_j}\vec{\Phi}\rangle+\langle\p_{x_i}\vec{\Phi},\p_{x_j}\vec{w}\rangle.
\end{split}
\end{equation}
Since $\sum_{i}(g_t)^{ki}(g_t)_{ij}=\delta_{kj}$ and $g_{ij}=e^{2\la}\, I_2$, where $I_2$ is the $(2\times 2)$-identity matrix, we have
\begin{equation}
\label{0III.10g}
  \left.\frac{d}{dt}(g_t)^{kj}\right|_{t=0}\, e^{2\la}+\,e^{-2\la}\,\left.\frac{d}{dt}(g_t)_{kj}\right|_{t=0}=0,
\end{equation}

from which we deduce
\be
\label{0III.10h}
\left.\frac{d}{dt}(g_t)^{kj}\right|_{t=0}=-e^{-4\la}\,\left.\frac{d}{dt}(g_t)_{kj}\right|_{t=0}=-e^{-4\la}\,\left(\langle \p_{x_k}\vec{w},\p_{x_j}\vec{\Phi}\rangle+\langle\p_{x_k}\vec{\Phi},\p_{x_j}\vec{w}\rangle\right).
\ee

%Thus,
%\begin{align}\label{longcalc}
%&\sum_{i,j} g^{ij}\,\left<\p_{x_i}\left.\frac{d}{dt}\vec{n}_t\right|_{t=0},\p_{x_j}\vec{\Phi}\right>\notag\\[5mm]
%&=-e^{-2\lambda}\, \p_{x_1}\left(e^{-2\lambda}\,\langle \langle\p_{x_1} \bw,\bn\rangle \ \p_{x_1}\bP+\langle\p_{x_2} \bw,\bn\rangle\ \p_{x_2}\bP,\p_{x_1}\bP\rangle \right)\notag\\[5mm]
%&\phantom{=}-e^{-2\lambda}\, \p_{x_2}\left(e^{-2\lambda}\,\langle \langle\p_{x_1} \bw,\bn\rangle \ \p_{x_1}\bP+\langle\p_{x_2} \bw,\bn\rangle\ \p_{x_2}\bP,\p_{x_2}\bP\rangle \right)\notag-2\, v\\[5mm]
%&=-e^{-2\lambda}\,\left(\p_{x_1}\langle\p_{x_1}\bw,\bn\rangle+\p_{x_2}\langle\p_{x_2}\bw,\bn\rangle\right)-2\, v.
%\end{align}

The previous computation gives then
\be
\label{vec-5}
\begin{array}{l}
\ds\frac{dH}{dt}=\frac{\vec{n}}{2}\cdot\lf[\Delta_g\vec{w}+2\,{\vec{w}}\rg]\\[5mm]
\ds\quad\quad+\frac{1}{2}\sum_{i,j=1}^2e^{-4\la}\lf[\p_{x_i}\vec{w}\cdot\p_{x_j}\vec{\Phi}+\p_{x_j}\vec{w}\cdot\p_{x_i}\vec{\Phi}\rg]\ \p_{x_j}\vec{n}\cdot\p_{x_i}\vec{\Phi}
\end{array}
\ee
We write
\be
\label{vec-11}
\vec{w}=\sigma_1\,\p_{x_1}\vec{\Phi}+\sigma_2\,\p_{x_2}\vec{\Phi}+v\, \vec{n}\quad.
\ee
and with these notations we compute in one hand 
\be
\label{vec-5-1}
\begin{array}{l}
\ds\frac{\vec{n}}{2}\cdot\Delta_g\vec{w}=\frac{\Delta_gv}{2}-\frac{|{\mathbb I}|^2_g}{2}\, v+e^{-2\la}\sum_{ij=1}^2\p_{x_i}\sigma_j\ {\mathbb I}_{ij}\\[5mm]
\ds\quad\quad +e^{-2\la}\sum_{i=1}^2\sigma_i\,\p_{x_i}(e^{2\la}H)
\end{array}
\ee
and using in particular (\ref{vec-3}) we have
\be
\label{vec-5-2}
\begin{array}{l}
\ds\frac{1}{2}\sum_{i,j=1}^2e^{-4\la}\lf[\p_{x_i}\vec{w}\cdot\p_{x_j}\vec{\Phi}+\p_{x_j}\vec{w}\cdot\p_{x_i}\vec{\Phi}\rg]\ \p_{x_j}\vec{n}\cdot\p_{x_i}\vec{\Phi}\\[5mm]
\ds\quad=|{\mathbb I}|^2_g\, v-\, e^{-2\la}\sum_{i,j=1}^2{\mathbb I}_{ij}\ \p_{x_i}\sigma_j-2\, H\, \sum_{k=1}^2\sigma_k\ \p_{x_k}\la
\end{array}
\ee
Combining (\ref{vec-5}), (\ref{vec-5-1}) and (\ref{vec-5-2}) we obtain
\be
\label{vec-5-3}
\begin{array}{l}
\ds\frac{dH}{dt}=\frac{1}{2}\lf[\Delta_g v+(|{\mathbb I}|^2_g+2)\, v\rg]+\sum_{i=1}^2\sigma_i\,\p_{x_i}H
\end{array}
\ee

Let $\ov{q}:=\ov{q}_{ij}\,dx_i\otimes dx_j=a\,[dx_1^2-dx_2^2]+b\,[dx_1\otimes dx_2+dx_2\otimes dx_1]$ for some functions $a$ and $b$ 
independent of $t$. We have
\be
\label{vec-6}
\lf<\frac{d {\mathbb I}^0}{dt},\ov{q}\rg>_g=\lf<\frac{d {\mathbb I}}{dt},\ov{q}\rg>_g-\, H\,\lf<\frac{d g}{dt},\ov{q}\rg>_g
\ee
We have in one hand
\be
\label{vec-7}
\begin{array}{l}
\ds-\, H\,\lf<\frac{d g}{dt},\ov{q}\rg>_g=-2\, H\ e^{-4\la}\, \lf[  a\,( \p_{x_1}\vec{w}\cdot\p_{x_1}\vec{\Phi}-\p_{x_2}\vec{w}\cdot\p_{x_2}\vec{\Phi})\rg.\\[5mm]
\ds\quad\quad\quad\lf.+b\, ( \p_{x_1}\vec{w}\cdot\p_{x_2}\vec{\Phi}+\p_{x_2}\vec{w}\cdot\p_{x_1}\vec{\Phi})\rg]
\end{array}
\ee
In the other hand we have
\be
\label{vec-8}
\begin{array}{l}
\ds\lf<\frac{d{\mathbb I}}{dt},\ov{q}\rg>=a\, e^{-4\la}\,\lf(\vec{n}\cdot\p^2_{x_1^2}\vec{w}-\vec{n}\cdot\p^2_{x_2^2}\vec{w}\rg)+2\, b\ e^{-4\la}\, \vec{n}\cdot\p^2_{x_1x_2}\vec{w}
\\[5mm]
\ds\quad-\, 2\, a\ e^{-4\la}\,(\p_{x_1}\la\,\p_{x_1}\vec{w}\cdot\vec{n}-\p_{x_2}\la\,\p_{x_2}\vec{w}\cdot\vec{n})\ \\[5mm]
\ds\quad-2\, b\ e^{-4\la}\,(\p_{x_1}\la\,\p_{x_2}\vec{w}\cdot\vec{n}+\p_{x_2}\la\,\p_{x_1}\vec{w}\cdot\vec{n})\ 
\end{array}
\ee
Thus, combining (\ref{vec-6}), (\ref{vec-7}) and (\ref{vec-8}) gives
\be
\label{vec-9}
\begin{array}{l}
\ds\lf<\frac{d {\mathbb I}^0}{dt},\ov{q}\rg>_g=a\, e^{-4\la}\,\lf(\vec{n}\cdot\p^2_{x_1^2}\vec{w}-\vec{n}\cdot\p^2_{x_2^2}\vec{w}\rg)+2\, b\ e^{-4\la}\, \vec{n}\cdot\p^2_{x_1x_2}\vec{w}
\\[5mm]
\ds\quad-\, 2\, a\ e^{-4\la}\,(\p_{x_1}\la\,\p_{x_1}\vec{w}\cdot\vec{n}-\p_{x_2}\la\,\p_{x_2}\vec{w}\cdot\vec{n})\ \\[5mm]
\ds\quad-2\, b\ e^{-4\la}\,(\p_{x_1}\la\,\p_{x_2}\vec{w}\cdot\vec{n}+\p_{x_2}\la\,\p_{x_1}\vec{w}\cdot\vec{n})\ \\[5mm]
\ds\quad-2\, H\ e^{-4\la}\, \lf[  a\,( \p_{x_1}\vec{w}\cdot\p_{x_1}\vec{\Phi}-\p_{x_2}\vec{w}\cdot\p_{x_2}\vec{\Phi})+b\, ( \p_{x_1}\vec{w}\cdot\p_{x_2}\vec{\Phi}+\p_{x_2}\vec{w}\cdot\p_{x_1}\vec{\Phi})\rg]
\end{array}
\ee
%We write
%\be
%\label{vec-11}
%\vec{w}=\sigma_1\,\p_{x_1}\vec{\Phi}+\sigma_2\,\p_{x_2}\vec{\Phi}+v\, \vec{n}\quad.
%\ee
We have then after some computations
\be
\label{vec-9-1}
\begin{array}{l}
\ds \vec{n}\cdot\p^2_{x_1^2}\vec{w}-\vec{n}\cdot\p^2_{x_2^2}\vec{w}= \p_{x_1}(\vec{n}\cdot\p_{x_1}\vec{w})-\p_{x_2}(\vec{n}\cdot\p_{x_2}\vec{w})- \p_{x_1}\vec{n}\cdot\p_{x_1}\vec{w}+\p_{x_2}\vec{n}\cdot\p_{x_2}\vec{w}\\[5mm]
\quad=\p^2_{x_1^2}v-\p^2_{x_2^2}v-\,v\ (|\p_{x_1}\vec{n}|^2-|\p_{x_2}\vec{n}|^2)+\, H\ e^{2\la}(\p_{x_1}\sigma_1-\p_{x_2}\sigma_2)\\[5mm]
\quad +e^{-2\la}\ {\mathbb I}^0_{11}\ (\p_{x_1}(e^{2\la}\sigma_1)+\p_{x_2}(e^{2\la}\sigma_2))+e^{-2\la}\ {\mathbb I}^0_{12}\ (\p_{x_1}(e^{2\la}\sigma_2)-\p_{x_2}(e^{2\la}\sigma_1))\\[5mm]
\quad+\p_{x_1}({\mathbb I}_{11}\sigma_1+{\mathbb I}_{12}\sigma_2)-\p_{x_2}({\mathbb I}_{12}\sigma_1+{\mathbb I}_{22}\sigma_2)
\end{array}
\ee
Similarly we have
\be
\label{vec-9-2}
\begin{array}{l}
\ds 2\,\vec{n}\cdot\p^2_{x_1x_2}\vec{w}=\p_{x_1}(\vec{n}\cdot\p_{x_2}\vec{w})+\p_{x_2}(\vec{n}\cdot\p_{x_1}\vec{w})-\p_{x_1}\vec{n}\cdot\p_{x_2}\vec{w}-\p_{x_2}\vec{n}\cdot\p_{x_1}\vec{w}\\[5mm]
\quad=2\,\p^{2}_{x_1x_2}v-2\,\p_{x_1}\vec{n}\cdot\p_{x_2}\vec{n}\ v+H\ e^{2\la}(\p_{x_2}\sigma_1+\p_{x_1}\sigma_2)\\[5mm]
\quad +e^{-2\la}\ {\mathbb I}^0_{11}\ (\p_{x_2}(e^{2\la}\sigma_1)-\p_{x_1}(e^{2\la}\sigma_2))+e^{-2\la}\ {\mathbb I}^0_{12}\ (\p_{x_1}(e^{2\la}\sigma_1)+\p_{x_2}(e^{2\la}\sigma_2))\\[5mm]
\quad+\p_{x_1}({\mathbb I}_{12}\sigma_1+{\mathbb I}_{22}\sigma_2)+\p_{x_2}({\mathbb I}_{11}\sigma_1+{\mathbb I}_{12}\sigma_2)
\end{array}
\ee
We have also
\be
\label{vec-9-3}
\begin{array}{l}
\p_{x_1}\la\,\p_{x_1}\vec{w}\cdot\vec{n}-\p_{x_2}\la\,\p_{x_2}\vec{w}\cdot\vec{n}=\p_{x_1}\la\,\p_{x_1}v-\p_{x_2}\la\,\p_{x_2}v\\[5mm]
\quad+e^{2\la}\,H\ (\sigma_1\p_{x_1}\la-\sigma_2\p_{x_2}\la)+{\mathbb I}^0_{11}(\sigma_1\p_{x_1}\la+\sigma_2\p_{x_2}\la)+{\mathbb I}^0_{12}(\sigma_2\p_{x_1}\la-\sigma_1\p_{x_2}\la)
\end{array}
\ee
and
\be
\label{vec-9-4}
\begin{array}{l}
\p_{x_1}\la\,\p_{x_2}\vec{w}\cdot\vec{n}+\p_{x_2}\la\,\p_{x_1}\vec{w}\cdot\vec{n}=\p_{x_1}\la\,\p_{x_2}v+\p_{x_2}\la\,\p_{x_1}v\\[5mm]
\quad+e^{2\la}\,H\ (\sigma_1\p_{x_2}\la+\sigma_2\p_{x_1}\la)+{\mathbb I}^0_{11}(\sigma_1\p_{x_2}\la-\sigma_2\p_{x_1}\la)+{\mathbb I}^0_{12}(\sigma_1\p_{x_1}\la+\sigma_2\p_{x_2}\la)
\end{array}
\ee
Finally we have in one hand
\be
\label{vec-9-5}
\p_{x_1}\vec{w}\cdot\p_{x_1}\vec{\Phi}-\p_{x_2}\vec{w}\cdot\p_{x_2}\vec{\Phi}=-\,2\, v\ {\mathbb I}^0_{11}+e^{2\la}\ (\p_{x_1}\sigma_1-\p_{x_2}\sigma_2)
\ee
and in the other hand
\be
\label{vec-9-6}
\p_{x_1}\vec{w}\cdot\p_{x_2}\vec{\Phi}+\p_{x_2}\vec{w}\cdot\p_{x_1}\vec{\Phi}=-\,2\, v\ {\mathbb I}^0_{12}+e^{2\la}\ (\p_{x_1}\sigma_2+\p_{x_2}\sigma_1)
\ee
Combining (\ref{vec-9})...(\ref{vec-9-6}) we obtain
\be
\label{vec-9-7}
\begin{array}{l}
\ds\lf<\frac{d {\mathbb I}^0}{dt},\ov{q}\rg>_g=a\,e^{-2\la}\lf(\p_{x_1}(e^{-2\la}\p_{x_1}v)-\p_{x_2}(e^{-2\la}\p_{x_2}v)\rg)+b\,e^{-2\la}\lf(\p_{x_1}(e^{-2\la}\p_{x_2}v)+\p_{x_2}(e^{-2\la}\p_{x_1}v)\rg)\\[5mm]
\ds\quad-\, e^{-4\la}\ H\, a\ (\p_{x_1}(e^{2\la}\sigma_1)-\p_{x_2}(e^{2\la}\sigma_2))-\, e^{-4\la}\ H\, b\ (\p_{x_1}(e^{2\la}\sigma_2)+\p_{x_2}(e^{2\la}\sigma_1))\\[5mm]
\ds\quad+ e^{-4\la}\ (a\,{\mathbb I}^0_{11}+b\,{\mathbb I}^0_{12})\ (\p_{x_1}\sigma_1+\p_{x_2}\sigma_2)+e^{-4\la}\ (a\,{\mathbb I}^0_{12}-b\,{\mathbb I}^0_{11})\ (\p_{x_1}\sigma_2-\p_{x_2}\sigma_1)\\[5mm]
\ds\quad+e^{-4\la}\ a\ \lf[ \p_{x_1}({\mathbb I}_{11}\sigma_1+{\mathbb I}_{12}\sigma_2)-\p_{x_2}({\mathbb I}_{12}\sigma_1+{\mathbb I}_{22}\sigma_2)\rg]\\[5mm]
\ds\quad +e^{-4\la}\ b\ \lf[\p_{x_1}({\mathbb I}_{12}\sigma_1+{\mathbb I}_{22}\sigma_2)+\p_{x_2}({\mathbb I}_{11}\sigma_1+{\mathbb I}_{12}\sigma_2)\rg]
\end{array}
\ee
Assume now that $a-ib$ is holomorphic - i.e. $\p_{x_1}a+\p_{x_2}b=0$ and $\p_{x_1}b-\p_{x_2}a=0$  then we have for instance
\[
\begin{array}{l}
\ds\quad a\,\lf(\p_{x_1}(e^{-2\la}\p_{x_1}v)-\p_{x_2}(e^{-2\la}\p_{x_2}v)\rg)+b\,\lf(\p_{x_1}(e^{-2\la}\p_{x_2}v)+\p_{x_2}(e^{-2\la}\p_{x_1}v)\rg)\\[5mm]
\ds\quad =\p_{x_1}(a\,e^{-2\la}\p_{x_1}v)-\p_{x_2}(a\,e^{-2\la}\p_{x_2}v)+\p_{x_1}(b\,e^{-2\la}\p_{x_2}v)+\p_{x_2}(b\,e^{-2\la}\p_{x_1}v)
\end{array}
\]
or we have also
\[
\begin{array}{l}
a\ \lf[ \p_{x_1}({\mathbb I}_{11}\sigma_1+{\mathbb I}_{12}\sigma_2)-\p_{x_2}({\mathbb I}_{12}\sigma_1+{\mathbb I}_{22}\sigma_2)\rg]+b\ \lf[\p_{x_1}({\mathbb I}_{12}\sigma_1+{\mathbb I}_{22}\sigma_2)+\p_{x_2}({\mathbb I}_{11}\sigma_1+{\mathbb I}_{12}\sigma_2)\rg]\\[5mm]
=\p_{x_1}\lf(a\lf[{\mathbb I}_{11}\sigma_1+{\mathbb I}_{12}\sigma_2\rg]+b\lf[({\mathbb I}_{12}\sigma_1+{\mathbb I}_{22}\sigma_2\rg]\rg)-\p_{x_2}\lf(a\lf[{\mathbb I}_{12}\sigma_1+{\mathbb I}_{22}\sigma_2\rg]-b \lf[ {\mathbb I}_{11}\sigma_1+{\mathbb I}_{12}\sigma_2\rg]\rg)
\end{array}
\]
and hence
\be
\label{vec-9-8}
\begin{array}{l}
\ds\lf<\frac{d {\mathbb I}^0}{dt},\ov{q}\rg>_g= e^{-2\la}\lf[\p_{x_1}(a\,e^{-2\la}\p_{x_1}v)-\p_{x_2}(a\,e^{-2\la}\p_{x_2}v)+\p_{x_1}(b\,e^{-2\la}\p_{x_2}v)+\p_{x_2}(b\,e^{-2\la}\p_{x_1}v) \rg]   \\[5mm]
\ds\quad-\, e^{-4\la}\ H\, a\ (\p_{x_1}(e^{2\la}\sigma_1)-\p_{x_2}(e^{2\la}\sigma_2))-\, e^{-4\la}\ H\, b\ (\p_{x_1}(e^{2\la}\sigma_2)+\p_{x_2}(e^{2\la}\sigma_1))\\[5mm]
\ds\quad+ e^{-4\la}\ (a\,{\mathbb I}^0_{11}+b\,{\mathbb I}^0_{12})\ (\p_{x_1}\sigma_1+\p_{x_2}\sigma_2)+e^{-4\la}\ (a\,{\mathbb I}^0_{12}-b\,{\mathbb I}^0_{11})\ (\p_{x_1}\sigma_2-\p_{x_2}\sigma_1)\\[5mm]
\ds\quad+e^{-4\la}\ \p_{x_1}\lf(a\lf[{\mathbb I}^0_{11}\sigma_1+{\mathbb I}^0_{12}\sigma_2\rg]+b\lf[({\mathbb I}^0_{12}\sigma_1-{\mathbb I}^0_{11}\sigma_2\rg]\rg)\\[5mm]
\ds\quad-e^{-4\la}\ \p_{x_2}\lf(a\lf[{\mathbb I}^0_{12}\sigma_1-{\mathbb I}_{11}^0\sigma_2\rg]-b \lf[ {\mathbb I}^0_{11}\sigma_1+{\mathbb I}^0_{12}\sigma_2\rg]\rg)\\[5mm]
\ds\quad+e^{-4\la}\ \lf[\p_{x_1}\lf(e^{2\la}\ H\ (a\,\sigma_1+b\,\sigma_2)\rg)+\p_{x_2}\lf(e^{2\la}\ H\ (b\,\sigma_1-a\,\sigma_2)\rg)\rg]
\end{array}
\ee

We now compute
\be
\label{vec-10}
\begin{array}{l}
\ds\frac{d}{dt}<\ov{q},{\mathbb I}^0>_g=\frac{d}{dt}\lf(\sum_{ijkl=1}^2 g^{ki}\,g^{jl} \,\ov{q}_{kl}\,{\mathbb I}^0_{ij}\rg)= 2\ e^{-2\la}\ \sum_{ijk=1}^2\frac{dg^{ij}}{dt}{\mathbb I}^0_{ik}\,\ov{q}_{jk}+\lf<\frac{d {\mathbb I}^0}{dt},\ov{q}\rg>_g
\end{array}
\ee
Using (\ref{0III.10h}) we obtain, using also (\ref{vec-3}),
\be
\label{vec-12}
\begin{array}{l}
\ds\frac{dg^{ij}}{dt}=2\ v\ {\mathbb I}_{ij}\ e^{-4\la}-\ e^{-4\la}\sum_{l=1}^2\lf[\p_{x_i}(\sigma_l\,\p_{x_l}\vec{\Phi})\cdot\p_{x_j}\vec{\Phi}+
 \p_{x_j}(\sigma_l\,\p_{x_l}\vec{\Phi})\cdot\p_{x_i}\vec{\Phi}  \rg]\\[5mm]
 \ds\quad\quad=2\ v\ {\mathbb I}_{ij}\ e^{-4\la}-e^{-2\la}\,(\p_{x_i}\sigma_j+\p_{x_j}\sigma_i)-2\, e^{-2\la}\, \delta_{ij}\ \sum_{l=1}^2\sigma_l\,\p_{x_l}\la
\end{array}
\ee
Combining (\ref{vec-10}) and (\ref{vec-12}) we obtain
\be
\label{vec-13}
\begin{array}{l}
\ds\frac{d}{dt}<\ov{q},{\mathbb I}^0>_g=4\, e^{-6\la}\, v\ \sum_{ijk=1}^2{\mathbb I}_{ij}\,{\mathbb I}^0_{ik}\,\ov{q}_{jk}+\lf<\frac{d {\mathbb I}^0}{dt},\ov{q}\rg>_g\\[5mm]
\ds\quad\quad-2\ e^{-4\la}\ \sum_{i,j,k}(\p_{x_i}\sigma_j+\p_{x_j}\sigma_i)\ {\mathbb I}^0_{ik}\ \ov{q}_{jk}-\, 4\ <\ov{q},{\mathbb I}^0>_g\ \sum_{l=1}^2\sigma_l\ \p_{x_l}\la
\end{array}
\ee
Observe that 
\[
\begin{array}{l}
\ds\sum_{i,j,k=1}^2{\mathbb I}^0_{ki}\,\ov{q}_{kj}\,{\mathbb I}^0_{ij}=\sum_{i,j=1}^2{\mathbb I}^0_{ji}\,\ov{q}_{jj}\,{\mathbb I}^0_{ij}+\sum_{i,j=1}^2{\mathbb I}^0_{j+1i}\,\ov{q}_{j+1j}\,{\mathbb I}^0_{ij}\\[5mm]
\ds\quad=a\,\sum_{i=1}^2{\mathbb I}^0_{1i}\,{\mathbb I}^0_{i1}-a\,\sum_{i=1}^2{\mathbb I}^0_{2i}\,{\mathbb I}^0_{i2}-b\,\sum_{i,j=1}^2{\mathbb I}^0_{j+1i}\,{\mathbb I}^0_{ij}
\end{array}
\]
Since ${\mathbb I}^0_{11}=-{\mathbb I}^0_{22}$ and ${\mathbb I}^0_{12}={\mathbb I}^0_{21}$ we deduce
\[
\begin{array}{l}
\ds\sum_{i,j,k=1}^2{\mathbb I}^0_{ki}\,\ov{q}_{kj}\,{\mathbb I}^0_{ij}=-b\,\sum_{i,j=1}^2{\mathbb I}^0_{j+1i}\,{\mathbb I}^0_{ij}\\[5mm]
\ds\quad=-b\,\sum_{i=1}^2{\mathbb I}^0_{1i}\,{\mathbb I}^0_{i2}-b\,\sum_{i=1}^2{\mathbb I}^0_{2i}\,{\mathbb I}^0_{i1}=-b\,{\mathbb I}^0_{11}\,{\mathbb I}^0_{12}-b\,{\mathbb I}^0_{12}\,{\mathbb I}^0_{22}-b\,{\mathbb I}^0_{21}\,{\mathbb I}^0_{11}-b\,{\mathbb I}^0_{22}\,{\mathbb I}^0_{21}=0
\end{array}
\]
Combining this identity and (\ref{vec-13}) we obtain
\be
\label{vec-14}
\begin{array}{l}
\ds\frac{d}{dt}<\ov{q},{\mathbb I}^0>_g=4\, H\, v\ <\ov{q},{\mathbb I}^0>_g+\lf<\frac{d {\mathbb I}^0}{dt},\ov{q}\rg>_g\\[5mm]
\ds\quad\quad-2\ e^{-4\la}\ \sum_{i,j,k}(\p_{x_i}\sigma_j+\p_{x_j}\sigma_i)\ {\mathbb I}^0_{ik}\ \ov{q}_{jk}-\, 4\ <\ov{q},{\mathbb I}^0>_g\ \sum_{l=1}^2\sigma_l\ \p_{x_l}\la
\end{array}
\ee
A short computation gives
\be
\label{vec-15}
-2\ e^{-4\la}\ \sum_{i,j,k}(\p_{x_i}\sigma_j+\p_{x_j}\sigma_i)\ {\mathbb I}^0_{ik}\ \ov{q}_{jk}=-\, 2\ (\p_{x_1}\sigma_1+\p_{x_2}\sigma_2)\ <\ov{q},{\mathbb I}^0>_g
\ee
Combining (\ref{vec-14}) and (\ref{vec-15}) one obtains
\be
\label{vec-16}
\begin{array}{l}
\ds\frac{d}{dt}<\ov{q},{\mathbb I}^0>_g=4\, H\, v\ <\ov{q},{\mathbb I}^0>_g+\lf<\frac{d {\mathbb I}^0}{dt},\ov{q}\rg>_g\\[5mm]
\ds\quad\quad-2\, e^{-2\la}\ (\p_{x_1}(e^{2\la}\,\sigma_1)+\p_{x_2}(e^{2\la}\, \sigma_2))\ <\ov{q},{\mathbb I}^0>_g
\end{array}
\ee
Together with (\ref{vec-9-8}) this gives
\be
\label{vec-17}
\begin{array}{l}
\ds    \frac{d}{dt}<\ov{q},{\mathbb I}^0>_g=4\, H\, v\ <\ov{q},{\mathbb I}^0>_g -2\, e^{-2\la}\ (\p_{x_1}(e^{2\la}\,\sigma_1)+\p_{x_2}(e^{2\la}\, \sigma_2))\ <\ov{q},{\mathbb I}^0>_g \\[5mm] 
\ds\quad +\,e^{-2\la}\lf[\p_{x_1}(a\,e^{-2\la}\p_{x_1}v)-\p_{x_2}(a\,e^{-2\la}\p_{x_2}v)+\p_{x_1}(b\,e^{-2\la}\p_{x_2}v)+\p_{x_2}(b\,e^{-2\la}\p_{x_1}v)\rg]    \\[5mm]
\ds\quad-\, e^{-4\la}\ H\, a\ (\p_{x_1}(e^{2\la}\sigma_1)-\p_{x_2}(e^{2\la}\sigma_2))-\, e^{-4\la}\ H\, b\ (\p_{x_1}(e^{2\la}\sigma_2)+\p_{x_2}(e^{2\la}\sigma_1))\\[5mm]
\ds\quad+ e^{-4\la}\ (a\,{\mathbb I}^0_{11}+b\,{\mathbb I}^0_{12})\ (\p_{x_1}\sigma_1+\p_{x_2}\sigma_2)+e^{-4\la}\ (a\,{\mathbb I}^0_{12}-b\,{\mathbb I}^0_{11})\ (\p_{x_1}\sigma_2-\p_{x_2}\sigma_1)\\[5mm]
\ds\quad+e^{-4\la}\ \p_{x_1}\lf(a\lf[{\mathbb I}^0_{11}\sigma_1+{\mathbb I}^0_{12}\sigma_2\rg]+b\lf[({\mathbb I}^0_{12}\sigma_1-{\mathbb I}^0_{11}\sigma_2\rg]\rg)\\[5mm]
\ds\quad-e^{-4\la}\ \p_{x_2}\lf(a\lf[{\mathbb I}^0_{12}\sigma_1-{\mathbb I}_{11}^0\sigma_2\rg]-b \lf[ {\mathbb I}^0_{11}\sigma_1+{\mathbb I}^0_{12}\sigma_2\rg]\rg)\\[5mm]
\ds\quad+e^{-4\la}\ \lf[\p_{x_1}\lf(e^{2\la}\ H\ (a\,\sigma_1+b\,\sigma_2)\rg)+\p_{x_2}\lf(e^{2\la}\ H\ (b\,\sigma_1-a\,\sigma_2)\rg)\rg]
\end{array}
\ee

Observe
\be
\label{III.3a}
\begin{array}{l}
\ds\left.\frac{d}{dt}\lf(dvol_{g_t}\rg)\right|_{t=0}=\left.\frac{d}{dt}\lf(det(g_t)_{ij}\rg)^{1/2}\right|_{t=0} dx_1\wedge dx_2\\[5mm]
\ds=\frac12 e^{-2\lambda}\ \left.\frac{d}{dt}\left((g_t)_{11}\ (g_t)_{22}-(g_t)_{12}^2\right)\right|_{t=0} dx_1\wedge dx_2\\[5mm]
\ds=\frac12 \left(\left.\frac{d}{dt}(g_t)_{11}\right|_{t=0}+ \left.\frac{d}{dt}(g_t)_{22}\right|_{t=0}\right) dx_1\wedge dx_2\\[5mm]
\ds=\left(\p_{x_1}\vec{\Phi}\cdot\p_{x_1}\bw+\p_{x_2}\vec{\Phi}\cdot\p_{x_2}\vec{w}\right)\,dx_1\wedge dx_2\\[5mm]
\ds=-2\vec{H}\cdot\vec{w}\ dvol_g+\lf(\p_{x_1}(e^{2\la}\sigma_1)+\p_{x_2}(e^{2\la}\sigma_2)\rg)\, dx_1\wedge dx_2.
\end{array}
\ee
\subsubsection{The first and second derivatives of the area under constrained conformal class}
We consider $\vec{\Phi}$ a weak immersion in ${\mathcal E}_\Sigma$. We fix generators of the $\pi_1$ and we denote by ${\mathcal C}(\vec{\Phi})$ 
the corresponding Teichm\"uller class of immersions which are conformal isotopic to $\vec{\Phi}$ - immersions
$\vec{\Xi}$ such that there exists a  Lipschitz diffeomorphism  $\Psi$ isotopic to the identity and conformal from $(\Sigma,\vec{\Phi}^\ast g_{S^3})$
into $(\Sigma,\vec{\Xi}^\ast g_{S^3})$ -. Let $D_0((-1,1),{\mathcal E}_{\Sigma})$ be the space of
mapping from $(-1,1)$ into ${\mathcal E}_\Sigma$ which are continuous and differentiable at 0. We say that $\vec{\Phi}$ is a critical point of the area
under constrained conformal class if
\be
\label{xsIII.22}
\lf\{
\begin{array}{l}
\ds\forall\ \vec{\Phi}_t\in D_0((-1,1),{\mathcal E}_{\Sigma})\quad\mbox{s.t. }\quad\forall\ t\in(-1,1) \quad{\mathcal C}(\vec{\Phi}_t)={\mathcal C}(\vec{\Phi})\\[5mm]
\ds\quad\mbox{ then}\quad \frac{d}{dt}A(\vec{\Phi}_t)(0)=0
\end{array}
\rg.
\ee
As shown in \cite{Ri4}, a weak immersion of $\Sigma$ is a critical point of the area under constrained conformal class if and only if there exists 
an holomorphic quadratic form $Q$ such that
\be
\label{xsIII.23}
H=\Re<{Q},h^0>_{wp}\quad.
\ee
writing locally $Q=(Q_1+i\,Q_2)\ dz^2$, since $$h^0=e^{2\la}\,[H^0_{\Re}+i\,H^0_{\Im}] \, dz^2=({\mathbb I}^0_{11}-i\,{\mathbb I}^0_{12})\ dz^2=-({\mathbb I}^0_{22}+i\,{\mathbb I}^0_{12})\, dz^2$$, we have that
\[
\Re<{Q},h^0>_{wp}=4\,e^{-4\la}\Re\lf((Q_1-i\,Q_2)\,({\mathbb I}^0_{11}-i\,{\mathbb I}^0_{12})\rg)=4\,e^{-4\la}\,[Q_1{\mathbb I}^0_{11}-Q_2{\mathbb I}^0_{12}]
\]
where we are using that $|dz^2|_g^2=4\,e^{-4\la}$. Let
 $$
 \ov{q}:=2\,\Re(Q)=2\, Q_1\,[dx_1^2-dx_2^2]-2\, Q_2\,[dx_1\,dx_2+dx_2\,dx_1]
 $$ 
 then we have
\[
<\ov{q},{\mathbb I}^0>_g=4\, e^{-4\la}\,[Q_1\,{\mathbb I}^0_{11}-Q_2{\mathbb I}^0_{12}]=\Re<{Q},h^0>_{wp}
\]
Hence $\vec{\Phi}$ is a critical point of the area
under constrained conformal class if and only if there exists $\ov{q}:=2\, Q_1\,[dx_1^2-dx_2^2]-2\, Q_2\,[dx_1\,dx_2+dx_2\,dx_1]$ such that
$Q_1+i\,Q_2$ is holomorphic and
\be
\label{xsIII.24}
H=<\ov{q},{\mathbb I}^0>_g\quad.
\ee

\medskip

Let $\vec{\Phi}_t$ be a mapping in $C^2((-1,1),{\mathcal E}_{\Sigma})$ such that ${\mathcal C}(\vec{\Phi}_t)\equiv{\mathcal C}(\vec{\Phi})$ for $t\in(-1,1)$
and such that $\vec{\Phi}_0=\vec{\Phi}$ is  a critical point of the area under constrained conformal class. We are now computing the second derivative
of $A(\vec{\Phi}_t)$ at $0$. Using (\ref{III.3a}) we obtain
\be
\label{xsIII.25}
\ds\frac{dA(\vec{\Phi}_{t})}{dt}=\frac{d}{d t}\lf[\int_{\Sigma^2}\,dvol_{g_{t}}\rg]=-2\int_{\Sigma}\vec{H}_{t}\cdot\frac{d\vec{\Phi}_{t}}{dt}\ dvol_{g_{t}}\quad,
\ee
where $g_{t}$ denotes the metric on $\Sigma$ given by $g_{t}:=\vec{\Phi}_{t}^\ast g_{S^3}$ and $\vec{H}_{t}$ is the mean curvature vector
of the immersion $\vec{\Phi}_{t}$.

\medskip

\be
\label{xsIII.26-c}
\begin{array}{l}
\ds\frac{d^2A(\vec{\Phi}_{t})}{dt^2}(0)=-2\int_{\Sigma}\frac{d{H}}{dt}(0)\ v\ dvol_{g}-2\int_{\Sigma}H\frac{d\vec{n}}{dt}\cdot\vec{w}\    dvol_{g}  \\[5mm]
\ds\quad\quad-2\int_{\Sigma}{H}\,\vec{n}\cdot\frac{d}{dt}\lf(\frac{d\vec{\Phi}_{t}}{dt}\ dvol_{g_{t}}\rg)
\end{array}
\ee
%A short computation gives
%\be
%\label{vec-18}
%\vec{w}\cdot\frac{d\vec{n}}{dt}=-\sum_{k=1}^2\sigma_k\,\p_{x_k}v-\sum_{ij=1}^2{\mathbb I}_{ij}\sigma_i\,\sigma_j\quad.
%\ee
Combining (\ref{xsIII.26-c}) and (\ref{vec-5-3}) gives
\be
\label{xsIII.26}
\begin{array}{l}
\ds\frac{d^2A(\vec{\Phi}_{t})}{dt^2}(0)=\int_{\Sigma}|dv|^2_g\ dvol_g-\int_{\Sigma}[|{\mathbb I}|^2_g+2]\ v^2\ dvol_g\\[5mm]
\ds\quad\quad-2\int_{\Sigma}v\,\sum_{k=1}^2 \sigma_k\ \p_{x_k} H\ dvol_{g}-2\int_{\Sigma}H\frac{d\vec{n}}{dt}\cdot\vec{w}\    dvol_{g} \\[5mm]
\ds\quad-2\int_{\Sigma}{H}\,\vec{n}\cdot\frac{d}{dt}\lf(\frac{d\vec{\Phi}_{t}}{dt}\ dvol_{g_{t}}\rg)
\end{array}
\ee

%\\[5mm]
%\ds\quad=\int_{\Sigma}\lf[|dw|_g^2-(|\vec{\mathbb I}|^2_g+2)\,w^2\rg]\ dvol_g-2\int_{\Sigma}{H}\,\vec{n}\cdot\frac{d}{dt}\lf(\frac{d\vec{\Phi}_{t}}{dt}\ dvol_{g_{t}}\rg)\
%\quad.
%\end{array}
%\ee
Since ${\mathcal C}(\vec{\Phi}_t)\equiv{\mathcal C}(\vec{\Phi})$ we have that
\[
d{\mathcal C}_{\vec{\Phi}_t}\frac{d \vec{\Phi}}{dt}\equiv0
\]
we have that for any $Q$ holomorphic quadratic form for the conformal class defined by ${\vec{\Phi}}^\ast g_{S^3}$
\be
\label{xsIII.27}
0\equiv\int_{\Sigma}\Re<{Q},\vec{h}_t^0>_{wp}\cdot\frac{d\vec{\Phi}_{t}}{dt}\ dvol_{g_{t}}=\int_{\Sigma}<\ov{q},\vec{\mathbb I}_t^0>_{g_t}\cdot\frac{d\vec{\Phi}_{t}}{dt}\ dvol_{g_{t}}\quad.
\ee
Taking the derivative of this identity gives
\be
\label{xsIII.28}
\int_{\Sigma}\frac{d}{dt}\lf(<\ov{q},\vec{\mathbb I}_t^0>_{g_t}\rg)\cdot\frac{d\vec{\Phi}_{t}}{dt}\ dvol_{g_{t}}=-\int_{\Sigma}<\ov{q},\vec{\mathbb I}_t^0>_{g_t}\cdot\frac{d}{dt}\lf(\frac{d\vec{\Phi}_{t}}{dt}\ dvol_{g_{t}}\rg)
\ee
%At $t=0$, since (\ref{xsIII.24}) holds and then the tangential variations do not affect the area at the second order and denoting $\vec{n}\cdot\frac{d\vec{\Phi}_{t}}{dt}(0)=w$, we have
We then have
\be
\label{xsIII.29}
\begin{array}{l}
\ds-2\,\int_{\Sigma} H\,\vec{n}\cdot\frac{d}{dt}\lf(\frac{d\vec{\Phi}_{t}}{dt}\ dvol_{g_{t}}\rg)=2\,\int_{\Sigma}\frac{d}{dt}\lf(<\ov{q},{\mathbb I}_t^0>_{g_t}\rg)\ v\ dvol_{g}\\[5mm]
\ds\quad\quad\quad+2\int_{\Sigma}\lf(<\ov{q},{\mathbb I}_t^0>_{g_t}\rg)\ \frac{d\vec{n}}{dt}\cdot \vec{w}\ dvol_{g}
\end{array}
\ee
Combining (\ref{xsIII.26}) and (\ref{xsIII.29}) we obtain
\be
\label{vec-18}
\begin{array}{l}
\ds\frac{d^2A(\vec{\Phi}_{t})}{dt^2}(0)=\int_{\Sigma}|dv|^2_g\ dvol_g-\int_{\Sigma}[|{\mathbb I}|^2_g+2]\ v^2\ dvol_g\\[5mm]
\ds\quad\quad-2\int_{\Sigma}v\,\sum_{k=1}^2 \sigma_k\ \p_{x_k} H\ dvol_{g} +2\,\int_{\Sigma}\frac{d}{dt}\lf(<\ov{q},{\mathbb I}_t^0>_{g_t}\rg)\ v\ dvol_{g}       
\end{array}
\ee
The identity (\ref{vec-17}) gives
\be
\label{vec-19}
\begin{array}{l}
\ds -2\int_{\Sigma}v\,\sum_{k=1}^2 \sigma_k\ \p_{x_k} H\ dvol_{g} +2\,\int_{\Sigma}\frac{d}{dt}\lf(<\ov{q},{\mathbb I}_t^0>_{g_t}\rg)\ v\ dvol_{g} \\[5mm]
\ds\quad =\int_{\Sigma}8 \, H^2\,v^2\ \ dvol_g-2\,\int_{\Sigma}
<\ov{q},dv\otimes dv>_g\  dvol_g\\[5mm]
\ds\quad-2\int_{\Sigma}v\,\sum_{k=1}^2 \sigma_k\ \p_{x_k} H\ dvol_{g} -\,4\,\int_{\Sigma} H\, v\, (\p_{x_1}(e^{2\la}\,\sigma_1)+\p_{x_2}(e^{2\la}\, \sigma_2))\ dx_1\wedge dx_2\\[5mm]
\ds\quad -\,2\,\int_{\Sigma} H\, v\ e^{-2\la}\ \lf[   a\ (\p_{x_1}(e^{2\la}\sigma_1)-\p_{x_2}(e^{2\la}\sigma_2))+\, b\ (\p_{x_1}(e^{2\la}\sigma_2)+\p_{x_2}(e^{2\la}\sigma_1)) \rg]\ dx_1\wedge dx_2\\[5mm]
\ds\quad+\,2\,\int_{\Sigma} \, v\ e^{-2\la}\ \lf[  (a\,{\mathbb I}^0_{11}+b\,{\mathbb I}^0_{12})\ (\p_{x_1}\sigma_1+\p_{x_2}\sigma_2)+ (a\,{\mathbb I}^0_{12}-b\,{\mathbb I}^0_{11})\ (\p_{x_1}\sigma_2-\p_{x_2}\sigma_1) \rg]\ dx_1\wedge dx_2\\[5mm]
\ds\quad+\,2\,\int_{\Sigma} \, v\ e^{-2\la}\ \p_{x_1}\lf(a\lf[{\mathbb I}^0_{11}\sigma_1+{\mathbb I}^0_{12}\sigma_2\rg]+b\lf[({\mathbb I}^0_{12}\sigma_1-{\mathbb I}^0_{11}\sigma_2\rg]\rg)\ dx_1\wedge dx_2\\[5mm]
\ds\quad-\,2\,\int_{\Sigma} \, v\ e^{-2\la}\ \p_{x_2}\lf(a\lf[{\mathbb I}^0_{12}\sigma_1-{\mathbb I}_{11}^0\sigma_2\rg]-b \lf[ {\mathbb I}^0_{11}\sigma_1+{\mathbb I}^0_{12}\sigma_2\rg]\rg)\ dx_1\wedge dx_2\\[5mm]
\ds\quad+\,2\,\int_{\Sigma} \, v\ e^{-2\la}\  \lf[\p_{x_1}\lf(e^{2\la}\ H\ (a\,\sigma_1+b\,\sigma_2)\rg)+\p_{x_2}\lf(e^{2\la}\ H\ (b\,\sigma_1-a\,\sigma_2)\rg)\rg]\ dx_1\wedge dx_2
\end{array}
\ee
We have in one hand

\be
\label{vec-20}
\begin{array}{l}
\ds\quad+\,2\,\int_{\Sigma} \, v\ e^{-2\la}\ \p_{x_1}\lf(a\lf[{\mathbb I}^0_{11}\sigma_1+{\mathbb I}^0_{12}\sigma_2\rg]+b\lf[({\mathbb I}^0_{12}\sigma_1-{\mathbb I}^0_{11}\sigma_2\rg]\rg)\ dx_1\wedge dx_2\\[5mm]
\ds\quad-\,2\,\int_{\Sigma} \, v\ e^{-2\la}\ \p_{x_2}\lf(a\lf[{\mathbb I}^0_{12}\sigma_1-{\mathbb I}_{11}^0\sigma_2\rg]-b \lf[ {\mathbb I}^0_{11}\sigma_1+{\mathbb I}^0_{12}\sigma_2\rg]\rg)\ dx_1\wedge dx_2\\[5mm]
\ds\quad=-\,2\,\int_{\Sigma} (a\,{\mathbb I}^0_{11}+b\,{\mathbb I}^0_{12})\ \lf[  \p_{x_1}(v\,e^{-2\la})\ \sigma_1+  \p_{x_2}(v\,e^{-2\la})\ \sigma_2\rg]\ dx_1\wedge dx_2\\[5mm]
\ds\quad\quad-\,2\,\int_{\Sigma} (a\,{\mathbb I}^0_{12}-b\,{\mathbb I}^0_{11})\ \lf[  \p_{x_1}(v\,e^{-2\la})\ \sigma_2- \p_{x_2}(v\,e^{-2\la})\ \sigma_1\rg]\ dx_1\wedge dx_2
\end{array}
\ee
Using the fact that $H=2\, e^{-4\la}\,(a\,{\mathbb I}^0_{11}+b\,{\mathbb I}^0_{12})$, a short computation gives
\be
\label{vec-21}
\begin{array}{l}
\ds-\,2\,\int_{\Sigma} H\, v\, (\p_{x_1}(e^{2\la}\,\sigma_1)+\p_{x_2}(e^{2\la}\, \sigma_2))\ dx_1\wedge dx_2\\[5mm]
\ds+\,2\,\int_{\Sigma} \, v\ e^{-2\la}\ \lf[  (a\,{\mathbb I}^0_{11}+b\,{\mathbb I}^0_{12})\rg] (\p_{x_1}\sigma_1+\p_{x_2}\sigma_2)\ dx_1\wedge dx_2\\[5mm]
\ds-\,2\,\int_{\Sigma} (a\,{\mathbb I}^0_{11}+b\,{\mathbb I}^0_{12})\ \lf[  \p_{x_1}(v\,e^{-2\la})\ \sigma_1+  \p_{x_2}(v\,e^{-2\la})\ \sigma_2\rg]\ dx_1\wedge dx_2\\[5mm]
\ds\quad=\int_{\Sigma}\ v\,e^{2\la}\lf[\sigma_1\,\p_{x_1}H+\sigma_2\,\p_{x_2}H\rg]\ dx_1\wedge dx_2
\end{array}
\ee
%\be
%\label{vec-22}
%\begin{array}{l}
%\ds+\,2\,\int_{\Sigma} \, v\ e^{-2\la}\ \lf[  (a\,{\mathbb I}^0_{12}-b\,{\mathbb I}^0_{11})\rg] (\p_{x_1}\sigma_2-\p_{x_2}\sigma_1)\ dx_1\wedge dx_2\\[5mm]
%\ds\quad\quad-\,2\,\int_{\Sigma} (a\,{\mathbb I}^0_{12}-b\,{\mathbb I}^0_{11})\ \lf[  \p_{x_1}(v\,e^{-2\la})\ \sigma_2- \p_{x_2}(v\,e^{-2\la})\ \sigma_1\rg]\ dx_1\wedge dx_2\\[5mm]
%\ds\quad=\int_{\Sigma}\ v\,e^{2\la}\lf[-\sigma_1\,\p_{x_2}<\ov{q}^\perp,{\mathbb I}^0>_g+\sigma_2\,\p_{x_1}<\ov{q}^\perp,{\mathbb I}^0>_g\rg]\ dx_1\wedge dx_2\\[5mm]
%\ds\quad+\,2\,\int_{\Sigma} \, v\, <\ov{q}^\perp,{\mathbb I}^0>_g\ (\p_{x_1}(e^{2\la}\,\sigma_2)+\p_{x_1}(e^{2\la}\, \sigma_2))\ dx_1\wedge dx_2
%\end{array}
%\ee
Another short computation gives also
\be
\label{vec-23}
\begin{array}{l}
\ds\quad -\,2\,\int_{\Sigma} H\, v\ e^{-2\la}\ \lf[   a\ (\p_{x_1}(e^{2\la}\sigma_1)-\p_{x_2}(e^{2\la}\sigma_2))+\, b\ (\p_{x_1}(e^{2\la}\sigma_2)+\p_{x_2}(e^{2\la}\sigma_1)) \rg]\ dx_1\wedge dx_2\\[5mm]
\ds\quad+\,2\,\int_{\Sigma} \, v\ e^{-2\la}\  \lf[\p_{x_1}\lf(e^{2\la}\ H\ (a\,\sigma_1+b\,\sigma_2)\rg)+\p_{x_2}\lf(e^{2\la}\ H\ (b\,\sigma_1-a\,\sigma_2)\rg)\rg]\ dx_1\wedge dx_2\\[5mm]
\ds\quad=2\,\int_{\Sigma} v\  \lf[(a\,\sigma_1+b\,\sigma_2)\ \p_{x_1}H+\,(b\sigma_1-a\,\sigma_2)\ \p_{x_2}H\rg]\ dx_1\wedge dx_2
\end{array}
\ee
We have also
\be
\label{vec-24}
\begin{array}{l}
\ds\quad-2\int_{\Sigma}v\,\sum_{k=1}^2 \sigma_k\ \p_{x_k} H\ dvol_{g} -\,2\,\int_{\Sigma} H\, v\, (\p_{x_1}(e^{2\la}\,\sigma_1)+\p_{x_2}(e^{2\la}\, \sigma_2))\ dx_1\wedge dx_2\\[5mm]
\ds\quad=2\ \int_{\Sigma}H\ \sum_{k=1}^2\sigma_k\ \p_{x_k}v\ dvol_g
\end{array}
\ee 
Combining (\ref{vec-18})....(\ref{vec-24}) gives
\be
\label{vec-25}
\begin{array}{l}
\ds\frac{d^2A(\vec{\Phi}_{t})}{dt^2}(0)=\int_{\Sigma}\lf[|dv|^2_g-2<\ov{q},dv\otimes dv>_g\,-[|{\mathbb I}|^2_g+2-8\, H^2]\ v^2\rg]\ dvol_g\\[5mm]
\ds\quad+2\ \int_{\Sigma}H\ \sum_{k=1}^2\sigma_k\ \p_{x_k}v\ dvol_g+\int_{\Sigma}\ v\, \sum_{k=1}^2\sigma_k\ \p_{x_k}H\ dvol_g\\[5mm]
\ds\quad+2\,\int_{\Sigma} v\ \lf[(a\,\sigma_1+b\,\sigma_2)\ \p_{x_1}H+\,(b\sigma_1-a\,\sigma_2)\ \p_{x_2}H\rg]\ dx_1\wedge dx_2\\[5mm]
\ds\quad-\,\int_{\Sigma} \, v\ \  \Im\lf<Q,h^0\rg>_{wp} \ (\p_{x_1}\sigma_2-\p_{x_2}\sigma_1)\ dvol_g\\[5mm]
\ds\quad+\,\int_{\Sigma}  e^{2\la}\,\Im\lf<Q,h^0\rg>_{wp}\ \lf[  \p_{x_1}(v\,e^{-2\la})\ \sigma_2- \p_{x_2}(v\,e^{-2\la})\ \sigma_1\rg]\ dvol_g\\[5mm]
\end{array}
\ee
We denote $\vec{w}_T$ to be the following vector-field tangent to the surface $\vec{w}_T:=\vec{w}-\vec{w}\cdot\vec{n}\,\vec{n}$. We have for instance
 $$\sum_{k=1}^2\sigma_k\ \p_{x_k}v=(dv,\vec{w}_T)_{T^\ast\Sigma,T\Sigma}\quad\mbox{ and }\quad \sum_{k=1}^2\sigma_k\ \p_{x_k}H=(dH,\vec{w}_T)_{T^\ast\Sigma,T\Sigma}$$
Observe that we have
\[
 \begin{array}{l}
 \ds -\,v\ (\p_{x_1}\sigma_2-\p_{x_2}\sigma_1)+  e^{2\la}\ \lf[  \p_{x_1}(v\,e^{-2\la})\ \sigma_2- \p_{x_2}(v\,e^{-2\la})\ \sigma_1\rg]\\[5mm]
\ds =-\,v\ e^{-2\la}\ \lf[\p_{x_1}(e^{2\la}\,\sigma_2)-\p_{x_2}(e^{2\la}\,\sigma_1)\rg]+  \lf( \p_{x_1}v\ \sigma_2- \p_{x_2}v\, \sigma_1\rg)
 \end{array}
\]
Denote $\vec{w}_T^{\,\ast}$ the 1-form dual to the vector-field $\vec{w}_T$ for the induced metric on the surface. We have in particular
\[
\ast\, d\vec{w}_T^{\,\ast}=e^{-2\la}\ \lf[\p_{x_1}(e^{2\la}\,\sigma_2)-\p_{x_2}(e^{2\la}\,\sigma_1)\rg]
\]
Finally we denote by $\ov{q}\res dH=2\, \Re(Q)\res dH$ the following contraction between the real part of the holomorphic quadratic form $2\, Q$ and the one form $dH$
\[
\ov{q}\res dH=e^{-2\la}\ a\, (\p_{x_1}H\, dx_1-\p_{x_2}H\, dx_2)+e^{-2\la}\ b\, (\p_{x_1}H\, dx_2+\p_{x_2}H\, dx_1)
\]

We can summarize the results obtained so far in the following lemma.
\begin{Lm}
\label{lm-I}
Let $\vec{\Phi}$ be a weak immersion in ${\mathcal E}_\Sigma(S^3)$. Assume that $\vec{\Phi}$ is a critical point of the area within the conformal class
defined by $\vec{\Phi}$, then there exists an holomorphic quadratic form $Q$ such that
\be
\label{xsIII.32}
H=\Re<{Q},h^0>_{wp}\quad.
\ee
Taking now a path such that $\vec{\Phi}_t$ is a mapping in $C^2((-1,1),{\mathcal E}_{\Sigma})(S^3)$ satisfying that ${\mathcal C}(\vec{\Phi}_t)\equiv{\mathcal C}(\vec{\Phi})$ for $t\in(-1,1)$
and such that $\vec{\Phi}_0=\vec{\Phi}$ and ${d\vec{\Phi}_{t}}/{dt}(0)=\vec{w}$, we have
\be
\label{xsIII.31-zz}
\begin{array}{l}
\ds\frac{d^2A(\vec{\Phi}_{t})}{dt^2}(0)=\int_{\Sigma}\lf[|dv|^2_g-4<\Re(Q),dv\otimes dv>_g\,-[|{\mathbb I}|^2_g+2-8\, H^2]\ v^2\rg]\ dvol_g\\[5mm]
\ds\quad+\int_\Sigma\lf[2H\, dv+v\,dH+4\,v\,\Re(Q)\res dH+\Im\lf<Q,h^0\rg>\, \ast dv   \rg]\cdot \vec{w}_T\ dvol_g\\[5mm]
\ds\quad-\int_\Sigma\Im\lf<Q,h^0\rg>\ v\, \ast d\vec{w}_T^{\,\ast}\ dvol_g\quad.
\end{array}
\ee
where $v:=\vec{w}\cdot\vec{n}$, where $\vec{w}_T$ is the vector-field tangent to the surface obtained by projecting $\vec{w}$ orthogonally (i.e. $\vec{w}_T:=\vec{w}-\vec{w}\cdot\vec{n}\,\vec{n}$) onto $\vec{\Phi}_\ast T\Sigma$. We denote $\vec{w}_T^{\,\ast}$ the 1-form
dual to $\vec{w}$ for the induced metric $g$. We have denoted $\Re(Q)\res dH$  the contraction with respect to the induced metric $g$ between the real part of the holomorphic quadratic form $Q$
and $dH$. The brakets $<\ ,\ >$ denotes the scalar product between quadratic forms induced by the metric $g$. Finally $\cdot$ is the canonical contraction between 1-forms and vector-fields on $\Sigma$.
\hfill $\Box$

\end{Lm}
The computation of the second derivative (\ref{xsIII.31-zz}) was a bit long due to the fact that we were considering general variations of the form $\vec{\Phi}_t$ where ${d\vec{\Phi}_{t}}/{dt}(0)=\vec{w}$ is not necessarily parallel to $\vec{n}$. This
was needed due to regularity issues. Since we are considering variations around a {\it weak immersion} which is not necessarily smooth $\vec{n}$ is only in $W^{1,2}$ and a deformation such as $\vec{\Phi}+t v\ \vec{n}$ would not be within the class
of weak immersions anymore.
\section{Proof of theorem~\ref{th-0.1}.}
\reset

Let $\vec{\Phi}$ be a weak immersion satisfying the {\it strictly elliptic conformally constrained minimal surface equation} that is to say, there exists an holomorphic quadratic form $Q$ such that
\[
H=\Re<Q,h^0>_{g_{\vec{\Phi}}}\quad\mbox{ and }\quad 2\,|Q|_{g_{\vec{\Phi}}}(x) < 1\quad \mbox{on }\quad T^2
\]
We can choose locally complex coordinates such that $Q=4^{-1}\ dz^2$. In these coordinates that we can assume to be defined on the disc $D^2$, the strict ellipticity condition $2\,|Q|_{g_{\vec{\Phi}}}(x) < 1$ becomes $\la>0$ where $g_{\vec{\Phi}}:= e^{2\la}\ [dx_1^2+dx_2^2]$ and the
conformally constrained minimal equation becomes
\be
\label{cc-1}
H=e^{-2\la} H^0_{\Re}\quad.
\ee
From the Codazzi equation (\ref{xsIII.20}) we then have
\be
\label{cc-2}
\Delta(e^{4\la} H)=\p_{x_1}(e^{2\la}\p_{x_1}H)-\p_{x_2}(e^{2\la}\p_{x_2}H)
\ee
Let $u:=e^{4\la} H$, assuming $\vec{\Phi}$ is a weak immersion gives that in complex coordinates $\nabla\la\in W^{1,1}$ and hence $\la\in C^0$ from which we deduce
that $u\in L^2$. The function $u$ satisfies moreover the elliptic PDE
\be
\label{cc-3}
\p_{x_1}((1-e^{-2\la})\,\p_{x_1}u)+\p_{x_2}((1+e^{-2\la})\p_{x_2}u)=-2\ \lf[\p_{x_1}(\p_{x_1}e^{-2\la}\ u)-\p_{x_2}(\p_{x_2}e^{2\la}\ u)\rg]
\ee
This PDE \footnote{Observe that this equation is critical for $u\in L^2$ or even $u\in L^{2,\infty}$. Indeed, $\nabla e^{2\la}$
is in the Lorentz space $L^{2,1}$ (see \cite{Ri1}) and thus $\nabla e^{2\la}\ u\in L^1$. The Laplacian of a function being a divergence of an $L^1$ function implies that this function is in $L^{2,\infty}$ and we are back to the space we are starting from which is the 
definition of being critical.}
 is indeed elliptic since $\la>0$ and $\la\in C^0$ thus we can assume that there exists $c_0>0$ such that $\la \ge c$ on $D^2$ and the symbol in the l.h.s. of (\ref{cc-3}) is invertible. If we can show that $u\in L^p_{loc}(D^2)$ for some $p>2$, this will make the PDE
 subcritical which implies imply by standard bootstrap arguments that $\nabla u\in L^{2,1}$ and then $\nabla H\in L^{2,1}$. Bootstraping this information respectively in 
 \[
 \Delta\vec{\Phi}+\vec{\Phi}\,|\nabla\vec{\Phi}|^2=2\, \vec{H}
 \]
 using Wente integrability by compensation theory (see \cite{Ri1}) we obtain that $\nabla^2\vec{\Phi}\in \cap_{p<+\infty}L^p$....which finally implies, using equation (\ref{0-4}), that $\vec{\Phi}$ is analytic.

We sketch now the proof of the fact that $u\in L^p$ for some $p>2$ which is the only claim that remains to be proven in order to complete the proof of theorem~\ref{th-0.1}. Before establishing the fact that $u\in L^p_{loc}(D^2)$ for some $p>2$ we first prove
the existence of $\al>0$ such that
\be
\label{cc-4}
\sup_{ x^0\in D^2_{1/2}\ ,\ 0<r<1/4} r^{-\al}\ \|u\|_{L^{2,\infty}(B_r(x^0))}<+\infty
\ee
Let $x^0\in D^2_{1/2}$ and $r<r_0$ where $r_0$ is going to be fixed later in the proof. Let $\chi$ be a cut-off function on ${\R}^2$ such that $\chi\equiv 1$ on $B_1^2(0)$ and $\supp(\chi)\subset B_2^2(0)$. We denote by $\chi_r^{x^0}(x)$ the function given by
$\chi_r^{x^0}(x):=\chi(r\,(x-x^0))$. A short computation gives
\be
\label{cc-4-a}
\begin{array}{l}
\ds \p_{x_1}\lf((1-e^{-2\la})\,\p_{x_1}(\chi_r^{x^0}\,u)\rg)+ \p_{x_2}\lf((1+e^{-2\la})\,\p_{x_2}(\chi_r^{x^0}\,u)\rg)\\[5mm]
\ds\quad=2\,\chi_r^{x^0}\, \lf[\p_{x_1}(\p_{x_1}e^{-2\la}\ u)-\p_{x_2}(\p_{x_2}e^{-2\la}\ u)\rg] \\[5mm]
\ds\quad+\,2\,\p_{x_1}\lf((1-e^{-2\la})\,\p_{x_1}(\chi_r^{x^0})\,u\rg)+ \,2\,\p_{x_2}\lf((1+e^{-2\la})\,\p_{x_2}(\chi_r^{x^0})\,u\rg)\\[5mm]
\ds\quad+\,\p_{x_1}(\chi_r^{x^0})\,\p_{x_1}e^{-2\la}\ u-\,\p_{x_2}(\chi_r^{x^0})\,\p_{x_2}e^{-2\la}\ u\\[5mm]
\ds\quad-\,(1-e^{-2\la})\,\p^2_{x^2_1}(\chi_r^{x^0})\,u-\,(1+e^{-2\la})\,\p^2_{x^2_2}(\chi_r^{x^0})\,u
\end{array}
\ee
We shall now use the following 3 lemmas.
\begin{Lm}
\label{lm-cc-2}
Let $\la$ measurable such that $\la\ge c_0>0$. For any $p\in (1,2)$ there exists $C_p(c_0)$ such that for any $X\in L^p(D^2,{\R}^2)$ and $f\in L^1$ there exists a unique solution $v\in W^{1,p}_0(D^2,{\R})$ of  
\[
\lf\{
\begin{array}{l}
\ds \p_{x_1}((1-e^{-2\la})\,\p_{x_1}v)+\p_{x_2}((1+e^{-2\la})\p_{x_2}v)=\mbox{div }X+f\quad\mbox{ in }D^2\\[5mm]
\ds v=0\quad\quad\mbox{ on }\p D^2
\end{array}
\rg.
\]
moreover we have
\be
\label{cc-5-a}
\|\nabla v\|_{L^p(D^2)}\le C(c_0)\ [\|X\|_{L^{p}(D^2)}+\|f\|_1]\quad.
\ee
\hfill $\Box$
\end{Lm}
The proof of this lemma is standard and can be done using the Stampacchia duality method. We then have the following result 
\begin{Lm}
\label{lm-cc-3}
Let $\la$ measurable such that $\la\ge c_0>0$. Let $v\in W^{1,1}(D^2,{\R})$ such that
 \[
\ds \p_{x_1}((1-e^{-2\la})\,\p_{x_1}v)+\p_{x_2}((1+e^{-2\la})\p_{x_2}v)=0\quad\mbox{ in }D^2\quad.
\]
Then for any $\delta\in (0,1)$ there exists $\theta\in (0,1)$ independent of $v$ such that
\be
\label{cc-5-b}
\|v\|_{L^{2,\infty}(B^2_{\theta})}\le\ \delta\ \|v\|_{L^{2,\infty}(B^2_1)}\quad.
\ee
\hfill $\Box$
\end{Lm}
This lemma is a straightforward consequence of the fact solutions to such an elliptic equation are smoother than $L^{2,\infty}$, from De Giorgi-Moser classical result we know that they are $W^{1,p}_{loc}$ for some $p>2$  with ad-hoc estimates
which imply (\ref{cc-5-b})
Finally we are going to make use of the following third lemma.
\begin{Lm}
\label{lm-cc-1}
Let $\la\in L^{2,1}(D^2,{\R})$ such that $\la\ge c_0>0$. For any $f\in L^{2,1}(D^2)$ there exists a unique $w\in W^{1,\infty}(D^2,{\R})$ such that
\[
\lf\{
\begin{array}{l}
\ds \p_{x_1}((1-e^{-2\la})\,\p_{x_1}w)+\p_{x_2}((1+e^{-2\la})\p_{x_2}w)=g\quad\mbox{ in }D^2\\[5mm]
\ds w=0\quad\quad\mbox{ on }\p D^2
\end{array}
\rg.
\]
moreover we have
\be
\label{cc-5}
\|\nabla w\|_{L^\infty(D^2)}\le C(\|\nabla\la\|_{2,1},c_0)\ \|g\|_{L^{2,1}(D^2)}\quad.
\ee
\hfill $\Box$
\end{Lm}
The proof of the previous lemma is classical if one replaces the Lorentz space $L^{2,1}$ by any $L^q$ space for $q>2$. Since the improved Sobolev embedding gives that identity map going from the space of functions with derivatives
in the Lorentz space $L^{2,1}$ into $C^0$ is continuous and since Calderon Zygmund theory extends  from the $L^q$ spaces to the Lorentz  $L^{2,1}$ space (see \cite{He}), the proof of lemma~\ref{lm-cc-1} follows easily.  

We go back to the proof of (\ref{cc-4}). We introduce $v\in W^{1,3/2}_0(B^2_{2r}(x^0))$ given by lemma~\ref{lm-cc-2} solving
\be
\label{cc-4-ab}
\begin{array}{l}
\ds \p_{x_1}\lf((1-e^{-2\la})\,\p_{x_1}v\rg)+ \p_{x_2}\lf((1+e^{-2\la})\,\p_{x_2}v\rg)\\[5mm]
\ds\quad=\,2\,\p_{x_1}\lf((1-e^{-2\la})\,\p_{x_1}(\chi_r^{x^0})\,u\rg)+ \,2\,\p_{x_2}\lf((1+e^{-2\la})\,\p_{x_2}(\chi_r^{x^0})\,u\rg)\\[5mm]
\ds\quad\quad\quad-\,(1-e^{-2\la})\,\p^2_{x^2_1}(\chi_r^{x^0})\,u-\,(1+e^{-2\la})\,\p^2_{x^2_2}(\chi_r^{x^0})\,u
\end{array}
\ee
Then $v$ satisfies
\be
\label{cc-4-abc}
\|v\|_{L^{2,\infty}(B^2_{2r}(x^0))}\le C(c_0)\ \|u\|_{L^{2\infty}(B^2_{2r}(x^0))}\quad.
\ee
Observe that the r.h.s of (\ref{cc-4-ab}) is supported on the annulus $B^2_{2r}(x^0))\setminus B^2_{r}(x^0))$ hence, making use of lemma~\ref{lm-cc-2} we deduce from the previous estimate that
\be
\label{cc-4-abc-1}
\forall\, \delta\in (0,1)\quad\exists\,  \theta\in (0,1)\quad\mbox{ s.t. }\quad\|v\|_{L^{2,\infty}(B^2_{\theta r}(x^0))}\le\ \delta\ \|u\|_{L^{2\infty}(B^2_{2r}(x^0))}
\ee
Now for any function $g\in L^{2,1}(B^2_{2r}(x^0))$ we consider $w_r^{x^0}$ given by lemma~\ref{lm-cc-1} for $\la$ on $B_{2r}^2(x^0)$ and satisfying then
\be
\label{cc-6}
\lf\{
\begin{array}{l}
\ds \p_{x_1}((1-e^{-2\la})\,\p_{x_1}w_r^{x^0})+\p_{x_2}((1+e^{-2\la})\p_{x_2}w_r^{x^0})=g\quad\mbox{ in }B_{2r}^2(x^0)\\[5mm]
\ds w=0\quad\quad\mbox{ on }\p B_{2r}^2(x^0) \quad.
\end{array}
\rg.
\ee
Due to scale dependancies of the norms we have
\be
\label{cc-7}
\|\nabla w_r^{x^0}\|_{L^\infty(B^2_{2r}(x^0))}\le C(\|\nabla\la\|_{L^{2,1}(B^2_{2r}(x^0))},c_0)\ \|f\|_{L^{2,1}(B^2_{2r}(x^0))}\quad.
\ee
Multiplying by $w_r^{x^0}-\ov{w_r^{x^0}}$ equation (\ref{cc-4-a}) to which we have subtracted equation  (\ref{cc-4-ab}) and integrating by parts, where $\ov{w_r^{x^0}}$ is the average of $w_r^{x^0}$ on $B^2_{2r}(x^0)$, gives
\be
\label{cc-8}
\int_{B^2_{2r}(x^0)} f\ (\chi_r^{x^0}\ u-v)\le C\ \|\nabla \la\|_{L^{2,1}(B^2_{2r}(x^0))}\ \|u\|_{L^{2,\infty}(B^2_{2r}(x^0))}\ \|\nabla w_r^{x^0}\|_{L^\infty(B^2_{2r}(x^0))}
\ee
Taking the sup over any $f$ supported in $B^2_{r}(x^0)$ whose $L^{2,1}$ norm 
is less than one gives then
\be
\label{cc-9}
\|u-v\|_{L^{2,\infty}(B^2_{r}(x^0))}\le C\ \|\nabla \la\|_{L^{2,1}(B^2_{2r}(x^0))}\ \|u\|_{L^{2,\infty}(B^2_{2r}(x^0))}\quad.
\ee
We choose now $r_0$ small enough in such a way that $ \sup_{x\in D^2_{1/2}}C\ \|\nabla \la\|_{L^{2,1}(B^2_{2r_0}(x^0))}<1/4$ and choosing $\theta$ in (\ref{cc-4-abc}) for $\delta=1/4$, we have found $\theta\in (0,1)$ independent
of $r$ such that
\[
\|u\|_{L^{2,\infty}(B^2_{\theta r}(x^0))}\le\ 2^{-1}\ \|u\|_{L^{2\infty}(B^2_{2r}(x^0))}
\]
A classical  iteration procedure of this inequality implies (\ref{cc-4}) and theorem~\ref{th-0.1} is proved.

\section{Isothermic conformally constrained surfaces.}
\reset

This section is devoted to the proof of theorem~\ref{th-III.1}.

\medskip
%\begin{Th}
%\label{th-III.1}
%Let $\vec{\Phi}$ be a weak immersion of the two torus $T^2$ in ${\mathcal E}_\Sigma(S^3)$. Assume that $\vec{\Phi}$ is a critical point of the area within the conformal class
%defined by $\vec{\Phi}$. Then $\vec{\Phi}$ is globally isothermic if and only if $\vec{\Phi}$ is either a minimal immersion or it is a flat CMC torus : a torus of the form $a \,S^1\times \sqrt{1-a^2}\, S^1$. 
%\hfill $\Box$
%\end{Th}
Let $\vec{\Phi}$ be a weak immersion of the torus $T^2$ equipped with a conformal structure $c$ into $S^3$. We assume that $\vec{\Phi}$ satisfies simultaneously the conformally constrained \'equation
\be
\label{Is-1}
H=\Re<Q,h^0>_{wp}
\ee
for some holomorphic quadratic differential $Q$ of $(T^2,c)$ and the isothermic condition
\be
\label{Is-2}
\Re<Q',h^0>_{wp}=0
\ee
for another non zero holomorphic quadratic differential $Q'$. 

Assume first $H\equiv0$, that is $\vec{\Phi}$ is minimal, it is clearly a critical point of the area under constrained conformal class. The Codazzi equation implies that the Weingarten quadratic form $h^0=\vec{n}\cdot\p^2_{z^2}\vec{\Phi}\, dz^2$ is a holomorphic. $h^0$ cannot be identically zero, otherwise one would contradict Liouville equation, since $H\equiv 0$. Taking $Q':=i\, h^0\ne 0$, we have (\ref{Is-2})
which implies that $\vec{\Phi}$ is isothermic.

 We shall exclude the case $H\equiv 0$ in this part from now on. Assuming (\ref{Is-1}), 
since $Q\ne 0$ and since the vector space of holomorphic quadratic differential on $T^2$ is a complex one dimensional space, we can then assume that there exist $\theta\in {\R}$ such that
\[
Q'=e^{\,i\,\theta}\, Q
\]
Since we are excluding the special case where $H\equiv 0$ we have $\theta\notin \pi\,{\Z}$. Modulo composition with a diffeomorphism we can assume that $\vec{\Phi}$ is conformal from
${\R}^2/\Lambda$ equipped with the flat metric where $\Lambda=\om_1{\Z}\oplus\om_2{\Z}$ with $\om_j\ne 0$ and $\om_2/\om_1\notin {\R}$ into $S^3$. In the canonical complex chart of ${\R}^2\simeq {\C}$
: $z=x_1+ix_2$ we have that $Q=f(z)\, dz^2$ where $f$ is holomorphic and periodic on ${\C}$ and therefore constant. Modulo rotations we can assume that $f(z)=1/4$ and then
\[
\Re<Q,h^0>_{wp}=e^{-2\la}\, H^0_{\Re}
\]
where $g=e^{2\la}\ (dx_1^2+dx_2^2)$ and the following system is satisfied 
\be
\label{Is-3}
\lf\{
\begin{array}{l}
\ds H=e^{-2\la}\, H^0_{\Re}\\[5mm]
\ds \cos\theta\, H^0_{\Re}-\sin\theta\, H^0_{\Im}=0
\end{array}
\rg.
\ee
and $\sin\theta> 0$. Let $t:=\cos\theta/\sin\theta$, the Codazzi identity (\ref{xsIII.20}) gives then
 \be
\label{Is-4}
(1+i\,t)\,\p_{\ov{z}}\lf(e^{2\la}\,H_{\Re}^0\rg)= e^{2\la}\ \p_z\lf(e^{-2\la}\,H^0_{\Re}\rg)\quad.
\ee
We have $2\ (1+i\,t)\,\p_{\ov{z}}=(\p_{x_1}-t\,\p_{x_2})+\,i\ (t\,\p_{x_1}+\p_{x_2})$. Hence   (\ref{Is-4}) becomes
\be
\label{Is-6}
\lf\{
\begin{array}{l}
\ds(\p_{x_1}-t\,\p_{x_2}){\mathbb I}^0_{11}=e^{-2\la}\,\p_{x_1}{\mathbb I}^0_{11}+e^{2\la}\,\p_{x_1}(e^{-4\la})\,{\mathbb I}^0_{11}\\[5mm]
\ds(t\,\p_{x_1}+\p_{x_2}){\mathbb I}^0_{11}=-e^{-2\la}\,\p_{x_2}{\mathbb I}^0_{11}-e^{2\la}\,\p_{x_2}(e^{-4\la})\,{\mathbb I}^0_{11}\quad.
\end{array}
\rg.
\ee
This gives
\be
\label{Is-7}
\lf\{
\begin{array}{l}
\ds (1-e^{-2\la})\, \p_{x_1}{\mathbb I}^0_{11}-\,t\ \p_{x_2}{\mathbb I}^0_{11}=2\, \p_{x_1}(e^{-2\la})\ {\mathbb I}^0_{11}\\[5mm]
\ds \ t\ \p_{x_1}{\mathbb I}^0_{11}+\, (1+e^{-2\la})\, \p_{x_2}{\mathbb I}^0_{11}=-2\, \p_{x_2}(e^{-2\la})\ {\mathbb I}^0_{11}\quad.
\end{array}
\rg.
\ee
From this identity we deduce
\be
\label{Is-8}
\lf\{
\begin{array}{l}
\ds (1+t^2-e^{-4\la})\ \p_{x_1}{\mathbb I}^0_{11}=2\, (1+e^{-2\la})\,\p_{x_1}(e^{-2\la})\ {\mathbb I}^0_{11}-2\,t\, \p_{x_2}(e^{-2\la})\ {\mathbb I}^0_{11}\\[5mm]
\ds  (1+t^2-e^{-4\la})\ \p_{x_2}{\mathbb I}^0_{11}=-2\, t\,\p_{x_1}(e^{-2\la})\ {\mathbb I}^0_{11}-2\, (1-e^{-2\la})\,\p_{x_2}(e^{-2\la})\ {\mathbb I}^0_{11}\quad.
\end{array}
\rg.
\ee
This implies the following conservation laws
\be
\label{Is-9}
\lf\{
\begin{array}{l}
\ds \p_{x_1}\lf[(1+t^2-e^{-4\la})\,{\mathbb I}^0_{11}\rg]=2\,\p_{x_1}(e^{-2\la})\ {\mathbb I}^0_{11}-2\,t\, \p_{x_2}(e^{-2\la})\ {\mathbb I}^0_{11}\\[5mm]
\ds \p_{x_2}\lf[(1+t^2-e^{-4\la})\,{\mathbb I}^0_{11}\rg]=-2\, t\,\p_{x_1}(e^{-2\la})\ {\mathbb I}^0_{11}-2\, \p_{x_2}(e^{-2\la})\ {\mathbb I}^0_{11}
\end{array}
\rg.
\ee
The fact that we are dealing with conservation laws can be seen as follows. Since, for a weak immersion, the conformal factor is continuous, we have
on the open set $\Om:=T^2\setminus (e^{-4\la})^{-1}\{1+t^2\}=T^2\setminus (e^{2\la})^{-1}\{\sin\theta\}$
\be
\label{Is-10}
\lf\{
\begin{array}{l}
\ds \p_{x_1}\lf[(1+t^2-e^{-4\la})\,{\mathbb I}^0_{11}\rg]=\lf[\p_{x_1}F\ -\,t\, \p_{x_2}F\rg]\ (1+t^2-e^{-4\la})\,{\mathbb I}^0_{11}\\[5mm]
\ds \p_{x_2}\lf[(1+t^2-e^{-4\la})\,{\mathbb I}^0_{11}\rg]=\lf[-\, t\,\p_{x_1}F- \p_{x_2}F\rg]\ (1+t^2-e^{-4\la})\,{\mathbb I}^0_{11}
\end{array}
\rg.
\ee
where
\be
\label{Is-11}
F(x):=\sin\theta\ \log\lf[\frac{e^{2\la}+\sin\theta}{e^{2\la}-\sin\theta}\rg]
\ee
Since $\nabla\la\in L^{2,1}(T^2)$ where $L^{2,1}$ is the Lorentz space whose dual is the weak Marcinkiewicz space $L^{2,\infty}$ (see the previous section) we have that $G:=(1+t^2-e^{-4\la})\,{\mathbb I}^0_{11}$
satisfies on $\Om$ the following PDE
\be
\label{Is-12}
\p_{\ov{z}}G= A(z)\ G(z)
\ee
where $A(z)\in L^{2,1}_{loc}(\Om)$. This equation implies $\Delta G=4\, \p_{z}(A\, G)$ which is critical for $G\in L^{2,\infty}$. The same approach as the one we employed for proving theorem~\ref{th-0.1} implies that $\nabla G\in L^{2,1}_{loc}$ and $G\in L^\infty_{loc}$.
Hence  $(1+t^2-e^{-4\la})\,{\mathbb I}^0_{11}\in L^\infty_{loc}(\Om)$ and using (\ref{Is-3}) together with the fact that $\la\in L^\infty$, the Liouville equation gives that $\Delta\la\in L^\infty_{loc}(\Om)$. Hence a standard bootstrap equation gives that $\vec{\Phi}$ is smooth on $\Om$. Thus both $G$ and $A$ are smooth on $\Om$.

 Carleman unique continuation argument applied to Beltrami equation (\ref{Is-12}) gives that $G$, and hence $H$, has only isolated zeros (we are excluding the fact $H\equiv 0$ since the beginning of this section.
Denote $\hat{\Om}$ the open set obtained by removing to $\Om$ the possibly existing isolated zeros of $H$, we then have
\be
\label{Is-13}
\lf\{
\begin{array}{l}
\ds \p_{x_1}\lf[\log\, G\rg]=\lf[\p_{x_1}F\ -\,t\, \p_{x_2}F\rg]\ \\[5mm]
\ds \p_{x_2}\lf[\log\, G\rg]=\lf[-\, t\,\p_{x_1}F- \p_{x_2}F\rg]\ 
\end{array}
\rg.
\ee
Or in other words we have
\be
\label{Is-13-a}
\p_{z}\log\,G=(1+i\,t)\,\p_{\ov{z}}F\quad.
\ee
This implies in particular
\be
\label{Is-14}
t\,\lf[\p^2_{x_1^2}F-\p^2_{x_2^2}F\rg]+2\,\p^2_{x_1x_2}F= 0\quad\quad\mbox{ on }\hat{\Om}
\ee
Since $\nabla F\in L^2_{loc}(\Om)$, we have $t\,\lf[\p^2_{x_1^2}F-\p^2_{x_2^2}F\rg]+2\,\p^2_{x_1x_2}F\in H^{-1}_{loc}(\Om)$ and is supported at isolated points.
A standard argument gives then
\be
\label{Is-15}
t\,\lf[\p^2_{x_1^2}F-\p^2_{x_2^2}F\rg]+2\,\p^2_{x_1x_2}F= 0\quad\quad\mbox{ on }{\Om}
\ee
Let $\sigma$ such that $\sinh 2\sigma=-t$, which is satisfied for
\[
\cosh^2\sigma=2^{-1}(1+\sin^{-1}\theta)\quad\mbox{ and }\quad\sinh^2\sigma=2^{-1}(-1+\sin^{-1}\theta)
\]
 and let $\tau\in {\R}$ such that
\[
\cos^2\tau=\frac{\cosh^2\sigma}{{\cosh^2\sigma+\sinh^2\sigma}}=\frac{1}{2}(1+\sin\theta)
\]
and 
\[
\sin^2\tau=\frac{\sinh^2\sigma}{{\cosh^2\sigma+\sinh^2\sigma}}=\frac{1}{2}(1-\sin\theta)
\]
In particular we choose $\cos\tau>0$ and sgn$(\sin\tau)=$sgn$(-t)$. Hence, since we are assuming $\theta\in (0,\pi)$, we take
\[
\tau=\frac{\theta}{2}-\frac{\pi}{4}
\]
With these notations equation (\ref{Is-15}) becomes
\be
\label{Is-16}
(\cos\tau\,\p_{x_1}+\sin\tau\,\p_{x_2})(-\sin\tau\,\p_{x_1}+\cos\tau\,\p_{x_2})\, F=0\quad\quad\mbox{ on }{\Om}
\ee
We proceed to a rotation of ${\R}^2$ by $\tau$ and denote $(y_1,y_2)$ the canonical coordinates after this rotation : 
\be
\label{Is-17}
\lf\{
\begin{array}{l}
\ds\frac{\p}{\p y_1}=\cos\tau\, \frac{\p}{\p x_1}+\sin\tau\, \frac{\p}{\p x_2}\\[5mm]
\ds\frac{\p}{\p y_2}=-\sin\tau\, \frac{\p}{\p x_1}+\cos\tau\, \frac{\p}{\p x_2}
\end{array}
\rg.
\ee
Then we have
in these rotated coordinates
\[
\p^2_{y_1y_2} F=0\quad\quad\mbox{ on }{\Om}
\]
Hence we deduce the existence of two functions $R$ and $S$ which are smooth on $\Om$ such that
\be
\label{Is-18}
F(y_1,y_2)=\sin\,\theta\ \log\lf[\frac{e^{2\la}+\sin\theta}{e^{2\la}-\sin\theta}\rg]=\sin\,\theta\ \lf[R(y_1)+S(y_2)\rg]\quad.
\ee
The system (\ref{Is-17}) can be rewritten as follows $\p_{x_1}+i\,\p_{x_2}=e^{i\tau}\, (\p_{y_1}+i\,\p_{y_2})$ or $\p_{\ov{z}}=e^{i\tau}\,\p_{\ov{w}}$ where we denote respectively $z=x_1+i\,x_2$ and $w=y_1+i\,y_2$. We have $\p_{{z}}=e^{-i\tau}\,\p_{{w}}$ hence $dz=e^{i\tau}\, dw$ which implies
\[
-\,i\,e^{i\,\theta}\,dw^2=dz^2
\]
Let $\hat{H}^0=\hat{H}^0_\Re+i\,\hat{H}^0_\Im$ be the expression of the Weingarten quadratic form in the $w$ coordinates. The intrinsic nature of the Weingarten quadratic form gives
\[
\hat{H}^0\,dw^2= H^0\,dz^2\quad\Longrightarrow\quad \hat{H}^0=-i\ e^{i\theta}\, H^0
\]
This gives
\be
\label{Is-18-a}
\lf\{
\begin{array}{l}
\ds\hat{H}^0_\Re=\sin\theta\, H^0_\Re+\cos\theta\, H^0_\Im=(\sin\theta)^{-1}\,{H^0_\Re}\\[5mm]
\ds\hat{H}^0_{\Im}=-\cos\theta\, H^0_\Re+\sin\theta\, H^0_\Im=0\quad.
\end{array}
\rg.
\ee
Hence, $y_1$ and $y_2$ are principal directions for $\mathbb{I}$ and the principal curvatures are given by
\be
\label{Is-18-b}
\lf\{
\begin{array}{l}
\ds\kappa_1=\hat{H}^0_\Re+H=\frac{e^{-2\la}}{\sin\theta}\,(e^{2\la}+\sin\theta)\, H^0_\Re\\[5mm]
\ds\kappa_2=-\hat{H}^0_\Re+H=\frac{e^{-2\la}}{\sin\theta}\,(-e^{2\la}+\sin\theta)\, H^0_\Re
\end{array}
\rg.
\ee
In these new coordinates the equation (\ref{Is-13-a}) becomes
\be
\label{Is-19}
\p_{w}\log\,G=(1+i\,t)\,e^{2\,i\,\tau}\,\p_{\ov{w}}F\quad.
\ee
We have $(1+i\,t)=(\sin\,\theta)^{-1} \, i\, e^{-i\theta}$ and since $e^{2\,i\,\tau}=-i\,e^{i\theta}$, we have
\be
\label{Is-20}
\p_{w}\log\,G=(\sin\,\theta)^{-1}\p_{\ov{w}}F=2^{-1}\,\lf[\dot{R}(y_1)+i\dot{S}(y_2)\rg]\quad.
\ee
Hence, for each connected component of $\Om$, we deduce the existence of a constant $C\in{\R}$ such that
\be
\label{Is-21}
\log\, G= R(y_1)-S(y_2)+C
\ee
The identity (\ref{Is-18}) gives
\be
\label{Is-22}
e^{2\la}=\sin\theta\ \frac{e^{R+S}+1}{e^{R+S}-1}
\ee
and then
\be
\label{Is-22-a}
1+t^2-e^{-4\la}=\frac{4}{\sin^2\theta}\,\frac{e^{R+S}}{(e^{R+S}+1)^2}
\ee
%In particular this gives
%\be
%\label{Is-23-a}
%H^0_\Re=e^C\ \sin^2\theta\,\frac{ e^{2\la}}{(e^{2\la}-\sin\theta)^2}\quad\mbox{ and }\quad H=e^C\ \frac{\sin^2\theta}{(e^{2\la}-\sin\theta)^2}
%\ee
%and
%\be
%\label{Is-23-b}
%\lf\{
%\begin{array}{l}
%\ds\kappa_1=e^C\,\sin\theta\ \frac{e^{2\la}+\sin\theta}{(e^{2\la}-\sin\theta)^2}\\[5mm]
%\ds\kappa_2=-\,e^C\, \frac{\sin\theta}{e^{2\la}-\sin\theta}
%\end{array}
%\rg.
%\ee
Liouville equation reads
\be
\label{Is-24}
-\Delta\la=e^{2\la}\ \lf[1+H^2-e^{-4\la}\,[({\mathbb I}^0_{11})^2+({\mathbb I}^0_{12})^2]\rg]
\ee
Using the fact that $H=e^{-4\la}\ {\mathbb I}_{11}^0$ and that ${\mathbb I}_{12}^0=t\,{\mathbb I}_{11}^0$ we have
\be
\label{Is-25}
-\Delta\la=e^{2\la}\ \lf[1-e^{-4\la}\ ({\mathbb I}^0_{11})^2\ [1+t^2-e^{-4\la}]\rg]
\ee
Recall that $G=(1+t^2-e^{-4\la})\,{\mathbb I}^0_{11}$, thus combining (\ref{Is-21}), (\ref{Is-22}) and ((\ref{Is-22-a}) we obtain
\be
\label{Is-23}
e^{-4\la}\ ({\mathbb I}^0_{11})^2\ [1+t^2-e^{-4\la}]=e^{-4\,\la}\ \frac{G^2}{1+t^2-e^{-4\la}}=\frac{e^{2C}}{4}\,{e^{R-S}}\, (e^{R}-e^{-S})^2
\ee
Using (\ref{Is-22}) a short computation gives
\be
\label{Is-29}
\nabla\la=\frac{2\,e^{R+S}}{e^{2(R+S)}-1}\ \lf(
\begin{array}{c}
\dot{R}\\[2mm]
\dot{S}
\end{array}
\rg)
\ee
and hence
\be
\label{Is-30}
\Delta\la=\frac{2\,e^{R+S}}{e^{2(R+S)}-1}\ \lf[\ddot{R}+\ddot{S}-[(\dot{R})^2+(\dot{S})^2]\,\frac{e^{2(R+S)}+1}{e^{2(R+S)}-1}\rg]
\ee
So the following ''double ODE'' is satisfied
\be
\label{Is-31}
\ddot{R}+\ddot{S}-[(\dot{R})^2+(\dot{S})^2]\,\frac{e^{2(R+S)}+1}{e^{2(R+S)}-1}=\frac{\sin\theta}{2}\  \frac{(e^{R+S}+1)^2}{e^{R+S}}\lf[1-\frac{e^{R+2C}}{4\,e^{S}}\, (e^{R}-e^{-S})^2\rg]\quad,
\ee
where $R$ depends on $y_1$ and $S$ depends on $y_2$. Hence both $R$ and $S$ are analytic functions. Assume both $R$ and $S$ are
non constant on each connected component of $\Om$ then each of them would satisfy infinitely many independent second order ODE which is a contradiction. So, on each connected component
of $\Om$, $\la$ and ${\mathbb I}$ depend either on $y_1$ or on $y_2$ exclusively.

Assume that $R$ is constant and that both $\la$ and $\mathbb I$ only depend on $y_2$ in a connected component of $\Om$. For any $y^0_2$ the curves $\Gamma_{y^0_2}$ given by
\[
\Gamma_{y^0_2}:=\lf\{\vec{\Phi}(y_1,y_2^0)\quad;\quad (y_1,y_2^0)\in \Om\rg\}
\]
are made of portions of planar circles. Indeed the unit tangent direction to $\Gamma_{y^0_2}$ are given by $\tau(y_1)=e^{-\la(y_2)}\,\p_{y_1}\vec{\Phi}(y_1,y_2^0)$. We have
\[
\p^2_{y_1^2}\vec{\Phi}= -\dot{\la}(y^0_2)\,\p_{y_2}\vec{\Phi}(y_1,y_2^0)+\kappa_1(y_2^0)\, e^{2\,\la(y_2^0)}\,\vec{n}(y_1,y_2^0)-e^{2\la(y^0_2)}\,\vec{\Phi}(y_1,y_2^0)\quad,
\]
and
\[
\p^3_{y_1^3}\vec{\Phi}= -\dot{\la}(y^0_2)\,\p^2_{y_1y_2}\vec{\Phi}(y_1,y_2^0)+\kappa_1(y_2^0)\, e^{2\,\la(y_2^0)}\,\p_{y_1}\vec{n}(y_1,y_2^0)-e^{2\la(y^0_2)}\,\p_{y_1}\vec{\Phi}(y_1,y_2^0)\quad,
\]
We have
\[
\p^2_{y_1y_2}\vec{\Phi}(y_1,y_2^0)=e^{-2\la(y_2^0)}\,\p^2_{y_1y_2}\vec{\Phi}(y_1,y_2^0)\cdot\p_{y_1}\vec{\Phi}\ \p_{y_1}\vec{\Phi}=\dot{\la}(y_2^0)\ \p_{y_1}\vec{\Phi}
\]
and
\[
\p_{y_1}\vec{n}(y_1,y_2^0)=e^{-2\la(y_2^0)}\,\p_{y_1}\vec{\Phi}\cdot\p_{y_1}\vec{n}(y_1,y_2^0)\ \p_{y_1}\vec{\Phi}=-\,\kappa_1(y_2^0)\ \p_{y_1}\vec{\Phi}\quad.
\]
Thus
\[
\p^3_{y_1^3}\vec{\Phi}= -\lf[\dot{\la}^2(y^0_2)+\kappa_1^2+e^{2\la}\rg]\,\p_{y_1}\vec{\Phi}(y_1,y_2^0)
\]
So the $k-$th osculating plane for $k\ge 2$ is constant equal to the 2 dimensional 2nd osculating plane. The fundamental theorem of curves in euclidian spaces asserts
that each connected component is included in a 2-dimensional plane. Moreover it has a constant curvature, so it is a portion of circle. The boundary of $\Om$
is made of points where $e^{2\la}=\sin\theta$ and $\la$ is continuous on ${\R}^2$, so if $\la(y_2^0)\ne 2^{-1}\,\log\sin\theta$, which is always the case for $y_2^0$ for
which there exists $y_1\in {\R}$ with $(y_1,y_2^0)\in\Om$, the closure of $\Gamma_{y^0_2}$ does not intersect the boundary of $\Om$ and $\Gamma_{y^0_2}$ is a closed
planar circle in ${S^3}$. Hence we have proved that each connected component of $\Om$ is a union of lines $y_2=cte$ or $y_1=cte$. Both situations cannot coexist and since again
$\la$ is $C^0$ on ${\R}^2$ and since $\la$ is constant on each component of the complement of $\Om$ and equal to $2^{-1}\,\log\sin\theta$, we have proved that $\la$ is either a function of $y_2$
globally or a function of $y_1$ globally. In any case, the immersed torus posses a foliation by planar circles giving principal directions.

We claim that these planar circles belong to parallel 2-planes. We have first
\be
\label{Is-32}
\p^2_{y_2y_1}\vec{\Phi}=\dot{\la}(y_2)\ e^{-2\la}\,\p_{y_1}\vec{\Phi}\quad.
\ee
Then we have
\[
\p^3_{y_2y_1^2}\vec{\Phi}=-\lf[\ddot{\la}+\dot{\la}^2+\kappa_1\,\kappa_2\, e^{2\la}+e^{2\la}\rg]\,\p_{y_2}\vec{\Phi}+e^{2\la}\,\lf[\dot{\kappa}_1+2\,\kappa_1\,\dot{\la}-\kappa_2\,\dot{\la}\rg]\ \vec{n}-\dot{\la}\ 
e^{2\la}\,\vec{\Phi}\quad.
\]
The Liouville equation reads
\[
-\ddot{\la}=\kappa_1\,\kappa_2\, e^{2\la}+e^{2\la}\quad,
\]
and the Codazzi equation reads
\[
\dot{\kappa}_1=-\dot{\la}\,(\kappa_1-\kappa_2)
\]
 Thus we have proved
\be
\label{Is-33}
 \p^3_{y_2y_1^2}\vec{\Phi}=\dot{\la}\,\p^2_{y_1^2}\vec{\Phi}\quad.
 \ee
The identities (\ref{Is-32}) and (\ref{Is-33}) imply that there exist generators of the osculating 2-planes, $\p_{y_1}\vec{\Phi}$ and $\p^2_{y_1^2}\vec{\Phi}$, whose $y_2$ derivatives still belong to the osculating 2-plane. We deduce that these 2-planes are independent
of $y_2$ and the immersed torus is then axially symmetric. 

Since in the $y$ coordinates we have $H=e^{-2\la} H^0_\Re=\sin\theta\,e^{-2\la}\ \hat{H}^0_\Re$ we have the following system
\be
\label{Is-34}
\lf\{
\begin{array}{l}
\ds -\ddot{\la}=\kappa_1\,\kappa_2\, e^{2\la}+e^{2\la}\\[5mm]
\ds \dot{\kappa}_1=-\dot{\la}\,(\kappa_1-\kappa_2)\\[5mm]
\ds \kappa_1+\kappa_2=\sin\theta\, e^{-2\la}\, (\kappa_1-\kappa_2)
\end{array}
\rg.
\ee
the last equation implies 
\be
\label{Is-34-a}
\kappa_2=\frac{\sin\theta-e^{2\la}}{\sin\theta+e^{2\la}}\,\kappa_1
\ee
Thus in particular
\be
\label{Is-35}
\kappa_2-\kappa_1=-\frac{2\,e^{2\la}}{\sin\theta+e^{2\la}}
\ee
Substituting (\ref{Is-35}) in the second equation of (\ref{Is-34}) gives
\[
\frac{\dot{\kappa}_1}{\kappa_1}=-\frac{2\, e^{2\la}\,\dot{\la}}{\sin\theta+e^{2\la}}
\]
from which we deduce the existence of $C\in {\R}$ such that
\be
\label{Is-36}
(e^{2\la}+\sin\theta)\, \kappa_1=e^C
\ee
Combining (\ref{Is-34-a}) and (\ref{Is-36}) gives
\be
\label{Is-37}
H=\frac{\kappa_1+\kappa_2}{2}=\frac{\sin\theta}{\sin\theta+e^{2\la}}\,\kappa_1=e^{-C}\ \sin\theta\, \kappa_1^2\quad.
\ee
Combining (\ref{Is-35}) and (\ref{Is-36}) gives
\be
\label{Is-38}
\kappa_2-\kappa_1=-2+2\,e^{-C}\,\sin\theta\ \kappa_1
\ee
 Combining (\ref{Is-37}) and (\ref{Is-38}) gives
\be
\label{Is-39}
\kappa_1= 1+e^{-C}\sin\theta\ (\kappa_1^2-\kappa_1)\quad.
\ee
Hence we deduce that $\kappa_1$ is constant on each connected component of $\Om$ and $\la$ is constant too. 
 So $\vec{\Phi}$ is flat on that component and $\la$ is constant with $e^{2\la}\ne\sin\theta$. The boundary of $\Om$
is made of points where $e^{2\la}=\sin\theta$, since $\la$ is continuous on ${\R}^2$, $\Om$ has no boundary and $\la$ and ${\mathbb I}$  are constant on ${\R}$. So $\vec{\Phi}$ is a flat CMC torus.
Moreover, because of the previous argument, it posses two compact foliations by planar circles and each circle represent a principal direction. So the torus is a rectangular one isometric
to a torus of the form $a \,S^1\times \sqrt{1-a^2}\, S^1$.

Assuming now $e^{2\la}=\sin\theta$ on the whole torus. Equations (\ref{Is-3}) and (\ref{Is-18-a}) give 
\[
\hat{H}^0_\Re=H\quad\quad\mbox{ and }\quad\quad\hat{H}^0_\Re=0\quad.
\]
Hence we deduce that the Gaussian curvature is identically equal to 1 which contradicts Liouville equation (\ref{Is-24}). 
 Hence we have proved theorem~\ref{th-III.1}.\hfill $\Box$
%$dw=e^{i\tau}\, dz$ which implies
%\[
%dw^2=-\,i\, e^{i\theta}\ dz^2
%\]
%Let $\hat{H}^0=\hat{H}^0_\Re+i\,\hat{H}^0_\Im$ be the expression of the Weingarten quadratic form in the $w$ coordinates. The intrinsic nature of the Weingarten quadratic form gives
%\[
%\hat{H}^0\,dw^2= H^0\,dz^2\quad\Longrightarrow\quad \hat{H}^0=-i\ e^{i\theta}\, H^0
%\]

%\be
%\label{Is-5}
%\lf\{
%\begin{array}{l}
%\ds(e^{2\la}-1)\,\p_{x_1}H^0_\Re-\, t\ e^{2\la}\,\p_{x_2}H^0_\Re=-2\ \lf[(1+e^{2\la})\,\p_{x_1}\la-\,t\ e^{2\la}\,\p_{x_2}\la\rg]\ H^0_\Re\\[5mm]
%\ds   t\ e^{2\la}\,\p_{x_1}H^0_\Re+(e^{2\la}+1)\,\p_{x_2}H^0_\Re=-2\ \lf[t\ e^{2\la}\,\p_{x_1}\la+(e^{2\la}-1)\, \p_{x_2}\la  \rg]\ H^0_\Re
%\end{array}
%\rg.
%\ee
%This implies
%\be
%\label{Is-6}
%\lf\{
%\begin{array}{l}
%\ds\lf[e^{4\la}\,(1+t^2)-1\rg]\, \p_{x_1}H^0_\Re=\lf[-2\, \lf[(1+e^{2\la})^2+t^2\,e^{4\la}\rg]\,\p_{x_1}\la+4\, t\, e^{2\la}\, \p_{x_2}\la\rg]\, H^0_\Re\\[5mm]
%\ds\lf[e^{4\la}\,(1+t^2)-1\rg]\, \p_{x_2}H^0_\Re=\lf[4\, t\, e^{2\la}\ \p_{x_1}\la-2\, \lf[(1-e^{2\la})^2+t^2\,e^{4\la}\rg]\, \p_{x_2}\la\rg]\, H^0_\Re
%\end{array}
%\rg.
%\ee

\section{Pairs of conformally constrained minimal surfaces which are conformally congruent.}
\reset
Let $\vec{\Phi}$ be a weak immersion of the torus $T^2$ which is a conformally constrained minimal surface in $S^3$ and assume there exists a conformal transformation $\Psi$ of $S^3$ 
which is not an isometry and such that $\Psi\circ\vec{\Phi}$ is again a conformally constrained minimal surface. 
Hence there exists $\vec{a}\in B^4\setminus \{0\}$ such that $\Psi_{\vec{a}}\circ\vec{\Phi}$ is a conformally constrained minimal surface where
\be
\label{IV.0}
\Psi_{\vec{a}}(\vec{y})=(1-|\vec{a}|^2)\frac{\vec{y}-\vec{a}}{|\vec{y}-\vec{a}|^2}-\vec{a}\quad.
\ee
Its conformal factor is given by 
\be
\label{IV.00}
\forall\ \vec{Y}\in T_{\vec{y}}S^3\quad\quad |d\Psi_{\vec{a}}\cdot \vec{Y}|=\frac{1-|\vec{a}|^2}{1+|\vec{a}|^2-2\,\vec{a}\cdot \vec{y}}\ |\vec{Y}|\quad.
\ee
Denote $\mu_{\vec{a}}(x)$ the function on $T^2$ given by
\[
e^{\mu_{\vec{a}}(x)}=\frac{1-|\vec{a}|^2}{1+|\vec{a}|^2-2\,\vec{a}\cdot \vec{\Phi}(x)}\
\]
From now on we shall omit the subscript $\vec{a}$ and simply write $\Psi$ and $\mu$.
We introduce the following parallelisation of $S^3$  
\[
\forall\,j\in {\Z}_3\quad\quad\vec{\eta}_j(y)=y_{j+1}\,\p_{y_{j-1}}-y_{j-1}\,\p_{y_{j+1}}+y_j\,\p_{y_4}-y_4\,\p_{y_j}
\]
One has
\be
\label{IV.1}
D_{\vec{\eta}_{j+1}}\vec{\eta}_j=\vec{\eta}_{j-1}\quad,\quad\quad D_{\vec{\eta}_{j-1}}\vec{\eta}_j=-\,\vec{\eta}_{j+1}\quad\mbox{ and }\quad D_{\vec{\eta}_{j}}\vec{\eta}_j=0
\ee
Where $D$ is the covariant derivative on $S^3$ associated to the Levi-Civita connection for the standard metric. One also verifies that
$(\vec{\eta}_1,\vec{\eta}_2,\vec{\eta}_3)$ realizes an orthonormal frame of $S^3$. Hence for each $y\in S^3$ the three unit vectors $e^{-\mu}\p_{\vec{\eta}_j}\Psi(y)$
realizes an orthonormal basis of $T_{\Psi(y)}S^3$. Let $a^j_{kl}$ be the real numbers such that
\[
D_{\vec{\eta}_l}\p_{\vec{\eta}_k}\Psi=\sum_{j=1}^3a_{lk}^j\  \p_{\vec{\eta}_j}\Psi\quad.
\]
Or in other words 
\[
a^j_{lk}=e^{-2\mu}\, D_{\vec{\eta}_l}\p_{\vec{\eta}_k}\Psi\cdot \p_{\vec{\eta}_j}\Psi=-\,e^{-2\mu}\, D_{\vec{\eta}_l}\p_{\vec{\eta}_j}\Psi\cdot \p_{\vec{\eta}_k}\Psi+2\,\delta_{jk}\ \p_{\vec{\eta}_l}\mu
\]
Let $\pi^{S^3}$ be the projection onto the tangent plane at $y\in S^3$ to $S^3$, we have
\be
\label{IV.2}
\begin{array}{l}
\ds D_{\vec{\eta}_l}\p_{\vec{\eta}_k}\Psi=\sum_{i,j=1}^4\pi^{S^3}\lf(\p_{y_j}(\p_{y_i}\Psi\,\eta_k^i)\, \eta^j_l\rg)\\[5mm]
\ds\quad=\pi^{S^3}\lf[\sum_{i,j=1}^4\lf(\p_{y_i}(\p_{y_j}\Psi\,\eta_l^j)\, \eta^i_k-\p_{y_j}\Psi\,\eta_k^i\,\p_{y_i}\eta^j_l+\p_{y_i}\Psi\,\eta_l^j\,\p_{y_j}\eta^i_k\rg)\rg]\\[5mm]
\ds\quad=D_{\vec{\eta}_k}\p_{\vec{\eta}_l}\Psi-\p_{D_{\vec{\eta}_k}\vec{\eta}_l}\Psi+\p_{D_{\vec{\eta}_l}\vec{\eta}_k}\Psi\\[5mm]
\ds\quad=D_{\vec{\eta}_k}\p_{\vec{\eta}_l}\Psi+2\, \p_{D_{\vec{\eta}_l}\vec{\eta}_k}\Psi
\end{array}
\ee
We are first computing $D_{\vec{\eta}_{j+1}}\p_{\vec{\eta}_j}\Psi$. We have
\be
\label{IV.3}
a^j_{j+1 j}=e^{-2\mu}\,D_{\vec{\eta}_{j+1}}\p_{\vec{\eta}_j}\Psi\cdot\p_{\vec{\eta}_j}\Psi=\p_{\vec{\eta}_{j+1}}\mu
\ee
and, using (\ref{IV.1}) and  (\ref{IV.2}) 
\be
\label{IV.4}
\begin{array}{l}
a^{j+1}_{j+1 j}=e^{-2\mu}\,D_{\vec{\eta}_{j+1}}\p_{\vec{\eta}_j}\Psi\cdot\p_{\vec{\eta}_{j+1}}\Psi=e^{-2\mu}\,D_{\vec{\eta}_{j}}\p_{\vec{\eta}_{j+1}}\Psi\cdot\p_{\vec{\eta}_{j+1}}\Psi=\p_{\vec{\eta}_{j}}\mu
\end{array}
\ee
Regarding the third coordinate $a^{j-1}_{j+1 j}$ we have in one hand using (\ref{IV.1}) and  (\ref{IV.2}) again
\be
\label{IV.5}
\begin{array}{l}
\ds D_{\vec{\eta}_{j+1}}\p_{\vec{\eta}_j}\Psi\cdot\p_{\vec{\eta}_{j-1}}\Psi=- D_{\vec{\eta}_{j+1}}\p_{\vec{\eta}_{j-1}}\Psi\cdot\p_{\vec{\eta}_{j}}\Psi\\[5mm]
\ds\quad=- D_{\vec{\eta}_{j-1}}\p_{\vec{\eta}_{j+1}}\Psi\cdot\p_{\vec{\eta}_{j}}\Psi-2\, \p_{D_{\vec{\eta}_{j+1}}{\vec{\eta}_{j-1}}}\Psi\cdot\p_{\vec{\eta}_{j}}\Psi\\[5mm]
\ds\quad=- D_{\vec{\eta}_{j-1}}\p_{\vec{\eta}_{j+1}}\Psi\cdot\p_{\vec{\eta}_{j}}\Psi+2\,e^{2\mu}\quad.
\end{array}
\ee
In the other hand we have
\be
\label{IV.6}
\begin{array}{l}
\ds D_{\vec{\eta}_{j+1}}\p_{\vec{\eta}_j}\Psi\cdot\p_{\vec{\eta}_{j-1}}\Psi=D_{\vec{\eta}_{j}}\p_{\vec{\eta}_{j+1}}\Psi\cdot\p_{\vec{\eta}_{j-1}}\Psi +2\,\p_{D_{\vec{\eta}_{j+1}}{\vec{\eta}_{j}}}\Psi\cdot\p_{\vec{\eta}_{j-1}}\Psi\\[5mm]
\quad=-D_{\vec{\eta}_{j}}\p_{\vec{\eta}_{j-1}}\Psi\cdot\p_{\vec{\eta}_{j+1}}\Psi +2\, e^{2\mu}\\[5mm]
\quad=-D_{\vec{\eta}_{j-1}}\p_{\vec{\eta}_{j}}\Psi\cdot\p_{\vec{\eta}_{j+1}}\Psi -\p_{D_{\vec{\eta}_{j}}\vec{\eta}_{j-1}}\Psi\cdot\p_{\vec{\eta}_{j+1}}\Psi +2\, e^{2\mu}\\[5mm]
\quad=-D_{\vec{\eta}_{j-1}}\p_{\vec{\eta}_{j}}\Psi\cdot\p_{\vec{\eta}_{j+1}}\Psi 
\end{array}
\ee
Summing (\ref{IV.5}) and (\ref{IV.6}) gives then
\be
\label{IV.7}
a^{j-1}_{j+1 j}= e^{-2\mu}\,D_{\vec{\eta}_{j+1}}\p_{\vec{\eta}_j}\Psi\cdot\p_{\vec{\eta}_{j-1}}\Psi=1\quad.
\ee
Combining now (\ref{IV.3}), (\ref{IV.4}) and (\ref{IV.7}) gives then
\be
\label{IV.8}
D_{\vec{\eta}_{j+1}}\p_{\vec{\eta}_j}\Psi=\p_{\vec{\eta}_{j+1}}\mu\ \p_{\vec{\eta}_j}\Psi+\p_{\vec{\eta}_{j}}\mu\ \p_{\vec{\eta}_{j+1}}\Psi+\p_{\vec{\eta}_{j-1}}\Psi
\ee
We are first computing $D_{\vec{\eta}_{j-1}}\p_{\vec{\eta}_j}\Psi$. We have
\be
\label{IV.3-a}
a^j_{j-1 j}=e^{-2\mu}\,D_{\vec{\eta}_{j-1}}\p_{\vec{\eta}_j}\Psi\cdot\p_{\vec{\eta}_j}\Psi=\p_{\vec{\eta}_{j-1}}\mu
\ee
and, using (\ref{IV.1}) and  (\ref{IV.2}) 
\be
\label{IV.4-a}
\begin{array}{l}
a^{j-1}_{j-1 j}=e^{-2\mu}\,D_{\vec{\eta}_{j-1}}\p_{\vec{\eta}_j}\Psi\cdot\p_{\vec{\eta}_{j-1}}\Psi=e^{-2\mu}\,D_{\vec{\eta}_{j}}\p_{\vec{\eta}_{j-1}}\Psi\cdot\p_{\vec{\eta}_{j-1}}\Psi=\p_{\vec{\eta}_{j}}\mu
\end{array}
\ee
Regarding the third coordinate $a^{j+1}_{j-1 j}$ we have seen
\be
\label{IV.9}
a^{j+1}_{j-1 j}= e^{-2\mu}\,D_{\vec{\eta}_{j-1}}\p_{\vec{\eta}_{j}}\Psi\cdot\p_{\vec{\eta}_{j+1}}\Psi =-1
\ee
Combining now (\ref{IV.3-a}), (\ref{IV.4-a}) and (\ref{IV.9}) gives then
\be
\label{IV.10}
D_{\vec{\eta}_{j-1}}\p_{\vec{\eta}_j}\Psi=\p_{\vec{\eta}_{j-1}}\mu\ \p_{\vec{\eta}_j}\Psi+\p_{\vec{\eta}_{j}}\mu\ \p_{\vec{\eta}_{j-1}}\Psi-\p_{\vec{\eta}_{j+1}}\Psi
\ee
A short computation gives also
\be
\label{IV.11}
D_{\vec{\eta}_{j}}\p_{\vec{\eta}_j}\Psi=\p_{\vec{\eta}_{j}}\mu\ \p_{\vec{\eta}_j}\Psi-\p_{\vec{\eta}_{j-1}}\mu\ \p_{\vec{\eta}_{j-1}}\Psi-\p_{\vec{\eta}_{j+1}}\mu\ \p_{\vec{\eta}_{j+1}}\Psi\quad.
\ee

We shall now establish the following lemma which is well known to experts in conformal geometry.
\begin{Lm}
\label{lm-IV-1}
Let $\vec{\Phi}$ be a weak immersion of the torus $T^2$ into $S^3$ and let $\Psi$ be a conformal transformation of the sphere $S^3$.
Denote  by $\vec{h}^0_{\vec{\Phi}}$ and by $\vec{h}^0_{\Psi\circ\vec{\Phi}}$ the Weingarten operators respectively for the immersions 
$\vec{\Phi}$ and $\Psi\circ\vec{\Phi}$ then the following identity holds
\be
\label{IV.12}
\vec{h}^0_{\Psi\circ\vec{\Phi}}= \Psi_\ast \vec{h}^0_{\vec{\Phi}}\quad.
\ee
\hfill $\Box$
\end{Lm}
\noindent{\bf Proof of lemma~\ref{lm-IV-1}.} We can assume that $\vec{\Phi}$ is conformal from a riemann surface into $S^3$. In some local complex coordinates
we denote by $e^{2\la}=|\p_{x_1}\vec{\Phi}|^2= |\p_{x_2}\vec{\Phi}|^2$. We have $\vec{h}^0_{\Phi}=e^{2\la}\,D_{z}(e^{-2\la}\p_z\vec{\Phi})\ dz^2$ where $D$ is the covariant derivative in $S^3$ issued from the Levi-Civita connection of the standard metric.  Let $e^{2\mu}$ be the conformal factor of the
conformal transformation $\Psi$. We have then $\vec{h}^0_{\Psi\circ\Phi}=e^{2\mu+2\la}\,D_{z}(e^{-2\mu-2\la}\p_z(\Psi\circ\vec{\Phi}))\ dz^2$. We compute
\be
\label{IV.13}
\begin{array}{l}
\ds D_{z}(e^{-2\mu-2\la}\p_z(\Psi\circ\vec{\Phi}))=-2\,\p_z\mu\,e^{-2\mu-2\la}\ \p_z(\Psi\circ\vec{\Phi})\\[5mm]
\ds\quad\quad+ e^{-2\mu}\ \sum_{j=1}^3 D_{z}\lf(e^{-2\la}\,\p_{\vec{\eta}_j}\Psi\ \vec{\eta}_j\cdot\p_z\vec{\Phi}\rg)
\end{array}
\ee
We have
\be
\label{IV.14}
\begin{array}{l}
\ds\sum_{j=1}^3 D_{z}\lf(e^{-2\la}\,\p_{\vec{\eta}_j}\Psi\ \vec{\eta}_j\cdot\p_z\vec{\Phi}\rg)=\Psi_\ast\lf(D_z\lf(e^{-2\la}\,\p_z\vec{\Phi}\rg)\rg)\\[5mm]
\ds\quad+e^{-2\la}\,\sum_{j=1}^3 D_{z}\lf(\p_{\vec{\eta}_j}\Psi\rg)\ \vec{\eta}_j\cdot\p_z\vec{\Phi}
+e^{-2\la}\,\sum_{j=1}^3 \p_{\vec{\eta}_j}\Psi\ D_z\vec{\eta}_j\cdot\p_z\vec{\Phi}
\end{array}
\ee
Observe first that, using (\ref{IV.1}), we have for any $j\in {\Z}_3$
\be
\label{IV.14a}
\begin{array}{l}
\ds D_z\vec{\eta}_j\cdot\p_z\vec{\Phi}=D_{\vec{\eta}_{j-1}}\vec{\eta}_j\cdot\p_z\vec{\Phi} \ \vec{\eta}_{j-1}\cdot\p_z\vec{\Phi}+D_{\vec{\eta}_{j+1}}\vec{\eta}_j\cdot\p_z\vec{\Phi} \ \vec{\eta}_{j+1}\cdot\p_z\vec{\Phi}\\[5mm]
\ds\quad=-\vec{\eta}_{j+1}\cdot\p_z\vec{\Phi}\ \vec{\eta}_{j-1}\cdot\p_z\vec{\Phi}+\vec{\eta}_{j-1}\cdot\p_z\vec{\Phi}\ \vec{\eta}_{j+1}\cdot\p_z\vec{\Phi}=0
\end{array}
\ee
Hence
\be
\label{IV.14b}
\begin{array}{l}
\ds\sum_{j=1}^3 D_{z}\lf(e^{-2\la}\,\p_{\vec{\eta}_j}\Psi\ \vec{\eta}_j\cdot\p_z\vec{\Phi}\rg)=\Psi_\ast\lf(D_z\lf(e^{-2\la}\,\p_z\vec{\Phi}\rg)\rg)\\[5mm]
\ds\quad+e^{-2\la}\,\sum_{j=1}^3 D_{z}\lf(\p_{\vec{\eta}_j}\Psi\rg)\ \vec{\eta}_j\cdot\p_z\vec{\Phi}
\end{array}
\ee
We have
\be
\label{IV.15}
D_{z}\lf(\p_{\vec{\eta}_j}\Psi\rg)=D_{\vec{\eta}_{j-1}}\lf(\p_{\vec{\eta}_j}\Psi\rg)\ \vec{\eta}_{j-1}\cdot\p_z\vec{\Phi}+D_{\vec{\eta}_{j}}\lf(\p_{\vec{\eta}_j}\Psi\rg)\ \vec{\eta}_{j}\cdot\p_z\vec{\Phi}+D_{\vec{\eta}_{j+1}}\lf(\p_{\vec{\eta}_j}\Psi\rg)\ \vec{\eta}_{j+1}\cdot\p_z\vec{\Phi}
\ee
We make now use of (\ref{IV.8}), (\ref{IV.10}) and (\ref{IV.11}) and we deduce
\be
\label{IV.16}
\begin{array}{l}
D_{z}\lf(\p_{\vec{\eta}_j}\Psi\rg)=\p_z\mu\ \p_{\vec{\eta}_j}\Psi\\[5mm]
\ds\quad+\p_{\vec{\eta}_{j-1}}\Psi\ \lf[\p_{\vec{\eta}_j}\mu\ \vec{\eta}_{j-1}\cdot\p_z\vec{\Phi}- \p_{\vec{\eta}_{j-1}}\mu\ \vec{\eta}_{j}\cdot\p_z\vec{\Phi}+\vec{\eta}_{j+1}\cdot\p_z\vec{\Phi}\rg]\\[5mm]
\ds\quad+\p_{\vec{\eta}_{j+1}}\Psi\ \lf[\p_{\vec{\eta}_j}\mu\ \vec{\eta}_{j+1}\cdot\p_z\vec{\Phi}- \p_{\vec{\eta}_{j+1}}\mu\ \vec{\eta}_{j}\cdot\p_z\vec{\Phi}-\vec{\eta}_{j-1}\cdot\p_z\vec{\Phi}\rg]
\end{array}
\ee
A relatively short computation gives
\be
\label{IV.16a}
\sum_{j=1}^3 D_{z}\lf(\p_{\vec{\eta}_j}\Psi\rg)\ \vec{\eta}_j\cdot\p_z\vec{\Phi}=2\,\p_z\mu\ \p_z(\Psi\circ\vec{\Phi})-\sum_{j,k=1}^3\p_{\vec{\eta}_j}\Psi\ \p_{\vec{\eta}_k}\Psi\ (\vec{\eta}_k\cdot\p_z\vec{\Phi})^2
\ee
Observe that
\be
\label{IV.16b}
\begin{array}{l}
\ds\sum_{k=1}^3(\vec{\eta}_k\cdot\p_z\vec{\Phi})^2=4^{-1}\sum_{k=1}^3\lf[(\p_{x_1}\vec{\Phi}\cdot\vec{\eta}_k)^2-(\p_{x_2}\vec{\Phi}\cdot\vec{\eta}_k)^2-2\,i\, \p_{x_1}\vec{\Phi}\cdot\vec{\eta}_k\,\p_{x_2}\vec{\Phi}\cdot\vec{\eta}_k\rg]\\[5mm]
\ds\quad\quad=4^{-1}\lf[|\p_{x_1}\vec{\Phi}|^2-|\p_{x_2}\vec{\Phi}|^2-2\,i\, \p_{x_1}\vec{\Phi}\cdot\p_{x_2}\vec{\Phi}\rg]=0
\end{array}
\ee
Combining (\ref{IV.14b}), (\ref{IV.16a}) and (\ref{IV.16b}) we obtain finally
\be
\label{IV.16c}
\ds\sum_{j=1}^3 D_{z}\lf(e^{-2\la}\,\p_{\vec{\eta}_j}\Psi\ \vec{\eta}_j\cdot\p_z\vec{\Phi}\rg)=\Psi_\ast\lf(D_z\lf(e^{-2\la}\,\p_z\vec{\Phi}\rg)\rg)+\,2\,\p_z\mu\ \p_z(\Psi\circ\vec{\Phi})
\ee
Combining (\ref{IV.13}) and (\ref{IV.16c}) gives (\ref{IV.12}) and lemma~\ref{lm-IV-1} is proved.\hfill $\Box$

\medskip

We establish now the following lemma.
\begin{Lm}
\label{lm-IV-2}
Let $\vec{\Phi}$ be a weak immersion of the torus $T^2$ into $S^3$ and let $\Psi$ be a conformal transformation of the sphere $S^3$.
Denote  by $\vec{H}_{\vec{\Phi}}$ and by $\vec{H}_{\Psi\circ\vec{\Phi}}$ the mean curvature vectors respectively for the immersions 
$\vec{\Phi}$ and $\Psi\circ\vec{\Phi}$ then the following identity holds
\be
\label{IV.17}
\vec{H}_{\Psi\circ\vec{\Phi}}=e^{-2\mu}\,\Psi_\ast\lf[\vec{H}_{\vec{\Phi}}-\nabla_{\vec{n}_\Phi}\mu\ \vec{n}_{\Phi}\rg]
\ee
where $e^{2\mu} g_{S^3}=\Psi^\ast g_{S^3}$ is the conformal factor of the conformal transformation $\Psi$, $\vec{n}_{\vec{\Phi}}$ is
the Gauss map of the immersion $\vec{\Phi}$ and $\nabla\mu$ is the gradient of  the function $\mu$ over $S^3$ for the standard metric on $S^3$.\hfill $\Box$
\end{Lm}
\noindent{\bf Proof of lemma~\ref{lm-IV-2}.} We can assume that $\vec{\Phi}$ is conformal from a riemann surface into $S^3$. In some local complex coordinates
we denote by $e^{2\la}=|\p_{x_1}\vec{\Phi}|^2= |\p_{x_2}\vec{\Phi}|^2$.
We have
\be
\label{IV.18}
\vec{H}_{\Psi\circ\vec{\Phi}}=2\ e^{-2\mu-2\la}\ D_{\ov{z}}\p_z(\Psi\circ\vec{\Phi})\quad.
\ee
We compute
\be
\label{IV.19}
\begin{array}{l}
\ds D_{\ov{z}}\p_z(\Psi\circ\vec{\Phi})=\psi_\ast D_{\ov{z}}\p_z\vec{\Phi}+\sum_{j=1}^3\p_{\vec{\eta}_j}\Psi\ \p_z\vec{\Phi}\cdot\p_{\ov{z}}\vec{\eta}_j\\[5mm]
\ds\quad+\sum_{j,k=1}^3D_{\vec{\eta}_k}\p_{\vec{\eta}_j}\Psi\ \p_{\ov{z}}\vec{\Phi}\cdot\vec{\eta}_k\, \p_{{z}}\vec{\Phi}\cdot\vec{\eta}_j
\end{array}
\ee
Using (\ref{IV.1}) we obtain in one hand
\be
\label{IV.20}
\sum_{j=1}^3\p_{\vec{\eta}_j}\Psi\ \p_z\vec{\Phi}\cdot\p_{\ov{z}}\vec{\eta}_j=\sum_{j=1}^3\p_{\vec{\eta}_j}\Psi\ \lf[\p_z\vec{\Phi}\cdot\vec{\eta}_{j-1}\ \p_{\ov{z}}\vec{\Phi}\cdot\vec{\eta}_{j+1}
-\p_z\vec{\Phi}\cdot\vec{\eta}_{j+1}\ \p_{\ov{z}}\vec{\Phi}\cdot\vec{\eta}_{j-1}\rg]
\ee
Using (\ref{IV.8}), (\ref{IV.10}) and (\ref{IV.11}) we obtain in the other hand
\be
\label{IV.21}
\begin{array}{l}
\ds\sum_{j,k=1}^3D_{\vec{\eta}_k}\p_{\vec{\eta}_j}\Psi\ \p_{\ov{z}}\vec{\Phi}\cdot\vec{\eta}_k\, \p_{{z}}\vec{\Phi}\cdot\vec{\eta}_j=\sum_{j,k=1}^3\p_{\vec{\eta}_j}\Psi\ \p_{z}\vec{\Phi}\cdot\vec{\eta}_j\ \p_{\vec{\eta}_k}\mu\ \p_{\ov{z}}\vec{\Phi}\cdot\vec{\eta}_k\\[5mm]
\ds\quad+\sum_{j,k=1}^3\p_{\vec{\eta}_j}\Psi\ \p_{\ov{z}}\vec{\Phi}\cdot\vec{\eta}_j\ \p_{\vec{\eta}_k}\mu\ \p_{{z}}\vec{\Phi}\cdot\vec{\eta}_k-\sum_{j,k=1}^3\p_{\vec{\eta}_j}\Psi\ \p_{\vec{\eta}_j}\mu\  \p_{{z}}\vec{\Phi}\cdot\vec{\eta}_k\ \p_{\ov{z}}\vec{\Phi}\cdot\vec{\eta}_k\\[5mm]
\ds\quad-\sum_{j=1}^3\p_{\vec{\eta}_j}\Psi\ \lf[\p_z\vec{\Phi}\cdot\vec{\eta}_{j-1}\ \p_{\ov{z}}\vec{\Phi}\cdot\vec{\eta}_{j+1}
-\p_z\vec{\Phi}\cdot\vec{\eta}_{j+1}\ \p_{\ov{z}}\vec{\Phi}\cdot\vec{\eta}_{j-1}\rg]
\end{array}
\ee
Combining (\ref{IV.19}), (\ref{IV.20}) and (\ref{IV.21}) gives then
\be
\label{IV.22}
\begin{array}{l}
\ds D_{\ov{z}}\p_z(\Psi\circ\vec{\Phi})=\psi_\ast D_{\ov{z}}\p_z\vec{\Phi}+\p_{z}(\Psi\circ\vec{\Phi})\ \p_{\ov{z}}\mu+\p_{\ov{z}}(\Psi\circ\vec{\Phi})\ \p_{{z}}\mu\\[5mm]
\ds\quad-\,|\p_z\vec{\Phi}|^2\ \sum_{j=1}^3\p_{\vec{\eta}_j}\mu\, \p_{\vec{\eta}_j}\Psi\quad.
\end{array}
\ee
Observe that we have
\be
\label{IV.23}
\p_{z}\vec{\Phi}\ \p_{\ov{z}}\mu+\p_{\ov{z}}\vec{\Phi}\ \p_{{z}}\mu=2^{-1}\lf[\p_{x_1}\mu\,\p_{x_1}\vec{\Phi}+\p_{x_2}\mu\,\p_{x_2}\vec{\Phi}\rg]
\ee
By definition $\nabla\mu:=\sum_{j=1}^3\p_{\vec{\eta}_j}\mu\,\vec{\eta}_j$. Denote $\pi_T$ the orthogonal projection onto the 2-plane tangent to $\vec{\Phi}(T^2)$. We have
at each point of the immersed surface
\be
\label{IV.24}
e^{2\la}\,\pi_T(\nabla\mu)=\sum_{i=1}^2\sum_{j=1}^3\p_{\vec{\eta}_j}\mu\,\vec{\eta}_j\cdot\p_{x_i}\vec{\Phi}\ \p_{x_i}\vec{\Phi}=\sum_{i=1}^2\p_{x_i}\mu\,\p_{x_i}\vec{\Phi}\quad.
\ee
Combining (\ref{IV.23}) and (\ref{IV.24}) gives
\be
\label{IV.25}
\p_{z}\vec{\Phi}\ \p_{\ov{z}}\mu+\p_{\ov{z}}\vec{\Phi}\ \p_{{z}}\mu=2^{-1}e^{2\la}\,\pi_T(\nabla\mu)\quad.
\ee
Thus
\be
\label{IV.26}
\begin{array}{l}
\ds\p_{z}(\Psi\circ\vec{\Phi})\ \p_{\ov{z}}\mu+\p_{\ov{z}}(\Psi\circ\vec{\Phi})\ \p_{{z}}\mu-\,|\p_z\vec{\Phi}|^2\ \sum_{j=1}^3\p_{\vec{\eta}_j}\mu\, \p_{\vec{\eta}_j}\Psi\\[5mm]
\ds\quad\quad=2^{-1}e^{2\la}\,\Psi_\ast\pi_{\vec{n}_{\vec{\Phi}}}(\nabla\mu)=2^{-1}e^{2\la}\,\Psi_\ast\lf[\nabla_{\vec{n}_{\vec{\Phi}}}\mu\ {\vec{n}_{\vec{\Phi}}}\rg]
\end{array}
\ee
Combining (\ref{IV.18}), (\ref{IV.22}) and (\ref{IV.26}) gives (\ref{IV.17}) and the lemma~\ref{lm-IV-2} is proved.\hfill $\Box$

We shall now prove the following lemma.

\begin{Lm}
\label{lm-IV-3}
Let $\vec{\Phi}$ be a weak immersion of the torus $T^2$. Assume $\vec{\Phi}$ is a conformally constrained minimal surface which is not isothermic, assume moreover it realizes it's conformal volume $A(\vec{\Phi})=V_c(\vec{\Phi})$
which is also assumed not to be a multiple of $4\pi$ and to be locally minimal in the conformal class defined by $\vec{\Phi}$. Assume moreover $\vec{\Phi}$ is a differentiable point of the conformal volume in $W^{1,\infty}\cap W^{2,2}(T^2)$ and that there exists
a conformal transformation $\Psi$ of $S^3$ which is not an isometry such that $A(\Psi\circ\vec{\Phi})=A(\vec{\Phi})$. Then  $\vec{\Phi}$ is a critical point of the area under the constraint $A(\Psi\circ\vec{\Phi})$ being constant.\hfill $\Box$
\end{Lm}
\noindent{\bf Proof of lemma~\ref{lm-IV-3}.}
Let $Q$ and $Q_\Psi$ be the two holomorphic quadratic differential of the torus $T^2$ equipped with the conformal class defined by $\vec{\Phi}^\ast g_{S^3}$, or equivalently $\vec{\Phi}^\ast \Psi^\ast g_{S^3}$ ,such that
\be
\label{IV.27}
\lf\{
\begin{array}{l}
H_{\vec{\Phi}}=\Re<Q,h^0_{\vec{\Phi}}>_{g_{\Phi}}\\[5mm]
H_{\Psi\circ\vec{\Phi}}=\Re<Q_\Psi,h^0_{\Psi\circ\vec{\Phi}}>_{g_{\Psi\circ\vec{\Phi}}}
\end{array}
\rg.
\ee
We now exclude the case $H_{\vec{\Phi}}=0$ or $H_{\Psi\circ\vec{\Phi}}=0$. Indeed if $\vec{\Phi}$ is minimal it is then isothermic and we are not considering this case and if $\Psi\circ\vec{\Phi}$
is minimal, it is well known that it is a strict maximum of the it's conformal volume and we could not have $A(\vec{\Phi})=A(\Psi\circ\vec{\Phi})$. 
Since we are excluding from now on the cases  $H_{\vec{\Phi}}=0$ or $H_{\Psi\circ\vec{\Phi}}=0$, then there exists $c\in {\C}^\ast$ such that
\[
Q=c\, Q_\Psi
\]
Assume first $c\in {\R}$. From lemma~\ref{lm-IV-1} we have 
\[
\Re<Q_\Psi,\vec{h}^0_{\Psi\circ\vec{\Phi}}>_{g_{\Psi\circ\vec{\Phi}}}=c\ e^{-4\mu}\, \Psi_\ast(\Re<Q,h^0_{\vec{\Phi}}>_{g_{\Phi}})
\]
Therefore for any $\vec{w}\in W^{1,\infty}\cap W^{2,2}(T^2)$ such that $\vec{w}\cdot \vec{\Phi}=0$ on $T^2$
\be
\label{IV.28}
dA_{\Psi\circ\vec{\Phi}}\cdot\Psi_\ast\vec{w}=-2\, c\,\int_{T^2}\ \vec{w}\cdot\Re<Q,\vec{h}^0_{\vec{\Phi}}>_{g_{\Phi}}\ dvol_g=\, c\ dA_{\vec{\Phi}}\cdot\vec{w}
\ee
So the differential of $\vec{\Phi}\rightarrow A(\Psi\circ\vec{\Phi})$ and $\vec{\Phi}\rightarrow A(\vec{\Phi})$ are proportional to each other. This implies the lemma in the case
when $c\in {\R}^\ast$.

Now we are going to rule out the case $c\in {\C}^\ast\setminus{\R}^\ast$. Let $\vec{w}\in W^{1,\infty}\cap W^{2,2}(T^2)$ such that $\vec{w}\cdot \vec{\Phi}=0$ on $T^2$ and $\vec{a}\in B^4_{1/2}$ and denote for $t$ small enough 
$\vec{\Phi}_t:=(\vec{\Phi}+t\vec{w})/|\vec{\Phi}+t\vec{w}|$. We denote $\vec{\Psi}_{\vec{a}}$ the conformal transformation
given by (\ref{IV.0}). Using (\ref{IV.00}) we have for $t$ small
\be
\label{IV.29}
\begin{array}{l}
\ds A(\vec{\Psi}_{\vec{a}}\circ\vec{\Phi}_t)-A(\vec{\Psi}_{\vec{a}}\circ\vec{\Phi})\\[5mm]
\ds\quad=(1-|\vec{a}|^2)^2\lf[\int_{T^2}\frac{dvol_{g_{\vec{\Phi}_t}}}{(1+|\vec{a}|^2-2\,\vec{a}\cdot\vec{\Phi}_t)^2}-\int_{T^2}\frac{dvol_{g_{\vec{\Phi}}}}{(1+|\vec{a}|^2-2\,\vec{a}\cdot\vec{\Phi})^2}\rg]
\end{array}
\ee
Hence we have
\be
\label{IV.30}
\begin{array}{l}
\ds\lf|A({\Psi}_{\vec{a}}\circ\vec{\Phi}_t)-A({\Psi}_{\vec{a}}\circ\vec{\Phi})-\lf|\frac{1-|\vec{a}|^2}{1+|\vec{a}|^2}\rg|^2(A(\vec{\Phi}_t)-A(\vec{\Phi}))\rg|\\[5mm]
\ds\quad\quad\le\ C_{\vec{\Phi}}\ |t|\, |\vec{a}|\, [\|\vec{w}\|_\infty+\|\nabla\vec{w}\|_\infty]
\end{array}
\ee
We have
\be
\label{IV.31}
\begin{array}{l}
\ds\lf|A({\Psi}_{\vec{a}}\circ\vec{\Phi}_t)-A({\Psi}_{\vec{a}}\circ\vec{\Phi})+2\, t\, \,\int_{T^2}\ \vec{w}\cdot\Re<Q,\vec{h}^0_{\vec{\Phi}}>_{g_{\vec{\Phi}}}\ dvol_{g_{\vec{\Phi}}}\rg|\\[5mm]
\quad\le\ C_{\vec{\Phi}}\ |t|\, |\vec{a}|\, [\|\vec{w}\|_\infty+\|\nabla\vec{w}\|_\infty]+o(t)
\end{array}
\ee
%Assume that $\Psi=\Psi_{\vec{a}_0}$ is the only $\Psi_{\vec{a}}\ne id$ such that $A(\vec{\Phi})=V_c(\vec{\Phi})=A(\Psi_{\vec{a}}\circ\vec{\Phi})$. 
Taking now the variation $\vec{\Phi}_t$ since
$V_c(\vec{\Phi})\notin 4\pi{\Z}$, for $t$ small enough there exists $\vec{a}_t$ such that $A(\Psi_{\vec{a}_t}\circ\vec{\Phi})=V_c(\vec{\Phi}_t)$ and $\vec{a}_t$ converges to the set $U\subset B^4$ of $\vec{a}\in B^4$ such that
$A(\Psi_{\vec{a}}\circ\vec{\Phi})=A(\vec{\Phi})=V_c(\vec{\Phi})$. Since $V_c(\vec{\Phi})$ is locally minimizing in the conformal class defined by $\vec{\Phi}$, for all $\vec{a}\in U$ 
$\Psi_{\vec{a}}\circ\vec{\Phi}$ is a conformally constrained minimal surface. Then from \cite{Ri4} there exist for each $\vec{a}$ an holomorphic quadratic form $Q_{\vec{a}}=c_{\vec{a}}\,Q$ such that
\be
\label{IV.32}
H_{\Psi_{\vec{a}}\circ\vec{\Phi}}=\Re<Q_{\vec{a}},h^0_{\Psi_{\vec{a}}\circ\vec{\Phi}}>_{g_{\Psi_{\vec{a}}\circ\vec{\Phi}}}
\ee
Recall that we are assuming $U\ne\{0\}$ and there exists $\vec{a}_0\in U\setminus\{0\}$ such that $\Psi=\Psi_{\vec{a}_0}$ and with our notations we have $Q_{\vec{a}_0}=Q_\Psi$.
It is also clear that $U$ is closed.

% Indeed, combining lemma~\ref{lm-IV-1}, lemma~\ref{lm-IV-2} and (\ref{IV.32}) we obtain for each $\vec{a}\in U$
%\be
%\label{IV.32a}
%H_{\vec{\Phi}}=e^{-2\mu_{\vec{a}}}\, \Re<c_{\vec{a}}\,Q,h^0_{\vec{\Phi}}>_{g_{\Phi}}+ 2^{-1}\nabla_{\vec{n}}e^{2\mu_{\vec{a}}}
%\ee 
%where we recall
%\[
%e^{\mu_{\vec{a}}(x)}=\frac{1-|\vec{a}|^2}{1+|\vec{a}|^2-2\,\vec{a}\cdot \vec{\Phi}(x)}\
%\]
%Hence
%\[
%\nabla_{\vec{n}}e^{2\mu_{\vec{a}}}=4\frac{(1-|\vec{a}|^2)^2\ \vec{a}\cdot\vec{n}(x)}{(1+|\vec{a}|^2-2\,\vec{a}\cdot \vec{\Phi}(x))^3}
%\]
%Observe that, as $\vec{a}$ varies  the functions of $x$ given by $ \Re<c_{\vec{a}}\,Q,h^0_{\vec{\Phi}}>_{g_{\Phi}}$ vary in a two dimensional space.

We are assuming that $V_c$ is differentiable at $\vec{\Phi}$ and hence there exists a linear form $L$ on $W^{1,\infty}\cap W^{2,2}$ such that
\[
|V_c(\vec{\Phi}_t)-V_c(\vec{\Phi})-\, t\, L(\vec{w})|=o(t)\quad.
\]
We claim now that
\be
\label{IV.33}
L(\vec{w})=\max_{\vec{a}\in U}\lf\{  -\, 2\, \int_{T^2}\ \vec{w}\cdot\Re<Q_{\vec{a}},\vec{h}^0_{\vec{\Phi}}>_{g_{\vec{\Phi}}}\ dvol_{g_{\vec{\Phi}}} \rg\}
\ee
\noindent{\bf Proof of the claim.}
Let $\vec{a}_t$ such that 
\[
V_c(\vec{\Phi}_t)=A(\Psi_{\vec{a}_t}\circ\vec{\Phi}_t)
\]
such a sequence exists for $t$ small enough since we are assuming $V_c(\vec{\Phi})\notin 4\pi{\Z}$. In order to simplify the presentation
we can assume that there exists $t_k\rightarrow 0$ such that $\vec{a}_{t_k}$ converges to $0$ since all the
$\Psi_{\vec{a}}$ for $\vec{a}\in U$ play the same role. We also omit to mention the subscript $k$. Using (\ref{IV.31}) we have
\be
\label{IV.34}
\begin{array}{l}
\ds\lf|V_c(\vec{\Phi}_t)-A({\Psi}_{\vec{a}_t}\circ\vec{\Phi})+2\, t\, \,\int_{T^2}\ \vec{w}\cdot\Re<Q,\vec{h}^0_{\vec{\Phi}}>_{g_{\vec{\Phi}}}\ dvol_{g_{\vec{\Phi}}}\rg|\\[5mm]
\quad\le\ C_{\vec{\Phi}}\ |t|\, |\vec{a}_t|\, [\|\vec{w}\|_\infty+\|\nabla\vec{w}\|_\infty]+o(t)
\end{array}
\ee
We deduce from this estimate
\be
\label{IV.35}
V_c(\vec{\Phi}_t)\le V_c(\vec{\Phi})-2\, t\, \,\int_{T^2}\ \vec{w}\cdot\Re<Q,\vec{h}^0_{\vec{\Phi}}>_{g_{\vec{\Phi}}}\ dvol_{g_{\vec{\Phi}}}+o(t)
\ee
This implies then
\be
\label{IV.36}
L(\vec{w})\le\max_{\vec{a}\in U}\lf\{  -\, 2\, \int_{T^2}\ \vec{w}\cdot\Re<Q_{\vec{a}},\vec{h}^0_{\vec{\Phi}}>_{g_{\vec{\Phi}}}\ dvol_{g_{\vec{\Phi}}} \rg\}
\ee
Let now $\vec{a}$ such that 
\be
\label{IV.37}
\max_{\vec{a}\in U}\lf\{  -\, 2\, \int_{T^2}\ \vec{w}\cdot\Re<Q_{\vec{a}},\vec{h}^0_{\vec{\Phi}}>_{g_{\vec{\Phi}}}\ dvol_{g_{\vec{\Phi}}} \rg\}=-\, 2\, \int_{T^2}\ \vec{w}\cdot\Re<Q_{\vec{a}},\vec{h}^0_{\vec{\Phi}}>_{g_{\vec{\Phi}}}\ dvol_{g_{\vec{\Phi}}}
\ee
We can assume without loss of generality that $\vec{a}=0$. We have
\[
V_c(\vec{\Phi}_t)\ge A(\vec{\Phi}_t)=A(\vec{\Phi})-\, 2\, \int_{T^2}\ \vec{w}\cdot\Re<Q,\vec{h}^0_{\vec{\Phi}}>_{g_{\vec{\Phi}}}\ dvol_{g_{\vec{\Phi}}}+o(t)
\]
This implies
\be
\label{IV.38}
L(\vec{w})\ge\max_{\vec{a}\in U}\lf\{  -\, 2\, \int_{T^2}\ \vec{w}\cdot\Re<Q_{\vec{a}},\vec{h}^0_{\vec{\Phi}}>_{g_{\vec{\Phi}}}\ dvol_{g_{\vec{\Phi}}} \rg\}
\ee
and the claim is proved.

 Since we are assuming that $Q_\Psi$ and $Q$ are not ${\R}-$parallel and since $\vec{\Phi}$ is assumed not to be isothermic, the space generated by the  linear forms $-\, 2\, \int_{T^2}\ \vec{w}\cdot\Re<Q_{\vec{a}},\vec{h}^0_{\vec{\Phi}}>_{g_{\vec{\Phi}}}\ dvol_{g_{\vec{\Phi}}}$ for $\vec{a}$ in $U$ is two dimensional. Then there exist two independent linear forms $L_1(\vec{w})$ and $L_2(\vec{w})$ and two functions $\al$ and $\beta$ on $U$ such that
  such that
 \[
 \begin{array}{l}
\forall \vec{a}\in U\quad\forall \vec{w}\in W^{1,\infty}\cap W^{2,2}(T^2,{\R}^4)\\[5mm]
\ds\quad-\, 2\, \int_{T^2}\ \vec{w}\cdot\Re<Q_{\vec{a}},\vec{h}^0_{\vec{\Phi}}>_{g_{\vec{\Phi}}}\ dvol_{g_{\vec{\Phi}}}=\al(\vec{a})\,L_1(\vec{w})+\beta(\vec{a})\,L_2(\vec{w})
\end{array}
\]
Since the space is bi-dimensional there exist $\vec{w}\ne\vec{w}'$ and $\vec{a}\ne\vec{a}'$ such that
\be
\label{IV.39}
\al(\vec{a}')\,L_1(\vec{w})+\beta(\vec{a}')\,L_2(\vec{w})<\al(\vec{a})\,L_1(\vec{w})+\beta(\vec{a})\,L_2(\vec{w})=L(\vec{w})
\ee
and 
\be
\label{IV.40}
\al(\vec{a})\,L_1(\vec{w}')+\beta(\vec{a})\,L_2(\vec{w}')<\al(\vec{a}')\,L_1(\vec{w}')+\beta(\vec{a}')\,L_2(\vec{w}')=L(\vec{w}')
\ee
For any $s,t\in {\R}$ one has, using the maximality property of $L$ we have for instance
\[
L(s\,\vec{w}+t\, \vec{w}')= s\, L(\vec{w})+t\, L(\vec{w}')\ge \al(\vec{a})\,L_1(s\vec{w}+t\vec{w}')+\beta(\vec{a})\,L_2(s\vec{w}+t\vec{w}')
\]
The three previous identities imply
\[
t\,\lf(\al(\vec{a}')\,L_1(\vec{w}')+\beta(\vec{a}')\,L_2(\vec{w}')\rg)\ge t\,\lf(\al(\vec{a}')\,L_1(\vec{w}')+\beta(\vec{a}')\,L_2(\vec{w}')\rg)
\]
This contradicts (\ref{IV.40}) for $t<0$. Hence the assumption that $Q_\Psi$ and $Q$ are not ${\R}-$parallel contradicts the fact that $V_c$ is differentiable at $\vec{\Phi}$. Lemma~\ref{lm-IV-3} is then proved.\hfill $\Box$

\section{Bounding the norm of the Lagrange Multiplier for Tori minimizing locally their conformal volume in their
conformal class.} 
\reset
Let $\vec{\Phi}$ be a weak immersion of the Torus i.e. an element in ${\mathcal E}_{T^2}$. We denote by
\[
{\mathcal T}_{T^2}({\vec{\Phi}})
\]
the subspace of ${\mathcal E}_{T^2}$ made of weak immersion which define the same Teichm\"uller class for some fixed generator
of the $\pi_1(T^2)$.

\medskip
In this subsection $\vec{\Phi}$ shall denote a weak immersion of $T^2$ minimizing locally the conformal volume in it's Teichm\"uller class.
\be
\label{0-I.1}
A(\vec{\Phi})=\inf_{\vec{\xi}\in {\mathcal U}}\sup_{\Psi\in G(S^3)} A(\Psi\circ\vec{\xi})
\ee
where $G(S^3)$ denotes the M\"obius group of conformal transformation of the sphere $S^3$ and ${\mathcal U}$ is an open neighborhood of $\vec{\Phi}$ in ${\mathcal T}_{T^2}({\vec{\Phi}})$ . 

\subsection{Families decreasing the area for small variations in the M\"obius group in a neighborhood of a non isothermic constrained minimal surface. }

The goal of the present subsection is to prove the following lemma
\begin{Lm}
\label{lm-0-I.2}
Let $\vec{\Phi}$ be a weak immersion of the torus. Assume $\vec{\Phi}$ is a critical point of the area under constrained conformal class which is not isothermic. Let $Q$ be an holomorphic quadratic differential such that
\[
\vec{H}=\Re(<Q,\vec{h^0}>_{wp})
\]
Then for ${\mathcal H}^2-$almost every point $x^0$ such that $2\,|Q|_{g_{\vec{\Phi}}}(x^0)>1$  there exists $\ep_{x^0,\vec{\Phi}}>0$,  $\nu_{x^0,\vec{\Phi}}>0$ and $\delta_{x^0,\vec{\Phi}}>0$  such that
\be
\label{0-I.18}
\begin{array}{l}
\forall \,\ep<\ep_{x^0,\vec{\Phi}}\quad\forall\,\vec{a}\in B_{\nu_{x^0,\vec{\Phi}}}^4(0)\quad
A(\Psi_{\vec{a}}\circ\vec{\Phi}_{x^0}^\ep(\delta_{x^0,\vec{\Phi}}\,\ep))< A(\Psi_{\vec{a}}\circ\vec{\Phi})-\, \ep^4\ C(x_0,\vec{\Phi})
\end{array}
\ee
for some positive constant $C(x_0,\vec{\Phi})>0$ depending only on $x^0$ and $\vec{\Phi}$ and not on $\ep$ and where $\vec{\Phi}_{x^0}^\ep(t)$ is the family defined in lemma~\ref{lm-xsIII.1a}.\hfill$\Box$
\end{Lm}
{\bf Proof of lemma~\ref{lm-0-I.2}.}
 Let $x_0$ be a point in $T^2$ which is a Lebesgue
point for $\nabla^2\vec{\Phi}$ and for $\nabla\vec{\Phi}$. Let ${\mathfrak U}_{x^0,\vec{\Phi}}$ be the open set given by lemma~\ref{lm-xsIII.1a} and denote by $\vec{\Phi}_{x^0}^\ep(t)$ the family of weak immersions of $T^2$ given by this lemma for $(t,\ep)$ in ${\mathfrak U}_{x^0,\vec{\Phi}}$. We have
\be
\label{0-I.2}
\frac{dA}{dt}(0)=0
\ee
and, due to lemma~\ref{lm-I}, $\vec{\Phi}$ satisfies the constraint minimal surfaces equation
\be
\label{0-I.2b}
\vec{H}=\Re<{Q},\vec{h}^0>_{wp}\quad.
\ee
for some holomorphic quadratic form $Q$. Hence in particular
\be
\label{0-I.3}
\begin{array}{l}
\ds A(\vec{\Phi}_{x^0}^\ep(t))=A(\vec{\Phi}) +\int_0^t(t-\tau)\,\frac{d^2 A(\vec{\Phi}_{x^0}^\ep(\tau))}{d\tau^2}\ d\tau\\[5mm]
\ds\quad\quad=A(\vec{\Phi}) -2\, \int_0^t(t-\tau)\,d\tau\,\frac{d}{d\tau}\lf[\int_\Sigma\vec{H}^\ep_{x^0}(\tau)\cdot\frac{d\vec{\Phi}^\ep_{x^0}}{d\tau}\ dvol_{g^\ep_{x^0}(\tau)}\rg]\\[5mm]
\ds\quad\quad=A(\vec{\Phi}) -2\, \int_0^t(t-\tau)\,d\tau\,\int_\Sigma\frac{d\vec{H}^\ep_{x^0}}{d\tau}\cdot\frac{d\vec{\Phi}^\ep_{x^0}}{d\tau}\ dvol_{g^\ep_{x^0}(\tau)}\\[5mm]
\ds\quad\quad\quad\quad\quad\quad-2\, \int_0^t(t-\tau)\,d\tau\,\int_\Sigma\vec{H}^\ep_{x^0}(\tau)\cdot\frac{d}{d\tau}\lf[\frac{d\vec{\Phi}^\ep_{x^0}}{d\tau}\ dvol_{g^\ep_{x^0}(\tau)}\rg]

\end{array}
\ee
Thus
\be
\label{0-I.4}
\begin{array}{l}
\ds A(\vec{\Phi}_{x^0}^\ep(t))=A(\vec{\Phi})-\, t^2\, \int_\Sigma\frac{d\vec{H}^\ep_{x^0}}{d\tau}(0)\cdot\frac{d\vec{\Phi}}{d\tau}(0)\ dvol_{g}\\[5mm]
\ds\quad\quad-2\, \int_0^t(t-\tau)\,d\tau\,\int_0^\tau\ ds\ \int_\Sigma\frac{d^2\vec{H}^\ep_{x^0}}{ds^2}\cdot\frac{d\vec{\Phi}}{ds}\ dvol_{g(s)}\\[5mm]
\ds\quad\quad-2\, \int_0^t(t-\tau)\,d\tau\,\int_0^\tau\ ds\ \int_\Sigma\frac{d\vec{H}^\ep_{x^0}}{ds}\cdot\frac{d}{ds}\lf[\frac{d\vec{\Phi}}{ds}\ dvol_{g(s)}\rg]\\[5mm]
\ds\quad\quad-2\, \int_0^t(t-\tau)\,d\tau\,\int_\Sigma[\vec{H}^\ep_{x^0}(\tau)-\vec{H}^\ep_{x^0}(0)]\cdot\frac{d}{d\tau}\lf[\frac{d\vec{\Phi}}{d\tau}\ dvol_{g(\tau)}\rg]\\[5mm]
\ds\quad\quad-2\, \int_0^t(t-\tau)\,d\tau\,\int_\Sigma\vec{H}^\ep_{x^0}(0)\cdot\frac{d}{d\tau}\lf[\frac{d\vec{\Phi}}{d\tau}\ dvol_{g(\tau)}\rg]
\end{array}
\ee
Since $\vec{H}^\ep_{x^0}(0)=\vec{H}$ satisfies (\ref{0-I.2b}) we have, omitting to write the subscript $x^0$ and the superscript $\ep$ from now on ,
\be
\label{0-I.5}
\begin{array}{l}
\ds-2\, \int_0^t(t-\tau)\,d\tau\,\int_\Sigma\vec{H}(0)\cdot\frac{d}{d\tau}\lf[\frac{d\vec{\Phi}}{d\tau}\ dvol_{g(\tau)}\rg]\\[5mm]
\ds=-2\, \int_0^t(t-\tau)\,d\tau\,\int_\Sigma\Re<{Q},\vec{h}^0(0)>_{g(0)}\cdot\frac{d}{d\tau}\lf[\frac{d\vec{\Phi}}{d\tau}\ dvol_{g(\tau)}\rg]\\[5mm]
\ds=-2\, \int_0^t(t-\tau)\,d\tau\,\int_\Sigma\Re<{Q},\vec{h}^0(\tau)>_{g(\tau)}\cdot\frac{d}{d\tau}\lf[\frac{d\vec{\Phi}}{d\tau}\ dvol_{g(\tau)}\rg]\\[5mm]
\ds-2\, \int_0^t(t-\tau)\,d\tau\,\int_\Sigma\Re\lf[<{Q},\vec{h}^0(0)>_{g(0)}-<Q,\vec{h}^0(\tau)>_{g(\tau)}\rg]\cdot\frac{d}{d\tau}\lf[\frac{d\vec{\Phi}}{d\tau}\ dvol_{g(\tau)}\rg]
\end{array}
\ee
Since ${\mathcal C}(\vec{\Phi}(t))\equiv {\mathcal C}(\vec{\Phi})$, using (\ref{xsIII.27}), we have then
\be
\label{0-I.6}
\begin{array}{l}
\ds-2\, \int_0^t(t-\tau)\,d\tau\,\int_\Sigma\vec{H}(0)\cdot\frac{d}{d\tau}\lf[\frac{d\vec{\Phi}}{d\tau}\ dvol_{g(\tau)}\rg]\\[5mm]
\ds=2\, \int_0^t(t-\tau)\,d\tau\,\int_\Sigma\frac{d}{d\tau}\lf[\Re<{Q},\vec{h}^0(\tau)>_{g(\tau)}\rg]\cdot\frac{d\vec{\Phi}}{d\tau}\ dvol_{g(\tau)}\\[5mm]
\ds-2\, \int_0^t(t-\tau)\,d\tau\,\int_\Sigma\Re\lf[<{Q},\vec{h}^0>_{g(0)}-<Q,\vec{h}^0(\tau)>_{g(\tau)}\rg]\cdot\frac{d}{d\tau}\lf[\frac{d\vec{\Phi}}{d\tau}\ dvol_{g(\tau)}\rg]
\end{array}
\ee
and then
\be
\label{0-I.7}
\begin{array}{l}
\ds-2\, \int_0^t(t-\tau)\,d\tau\,\int_\Sigma\vec{H}(0)\cdot\frac{d}{d\tau}\lf[\frac{d\vec{\Phi}}{d\tau}\ dvol_{g(\tau)}\rg]\\[5mm]
\ds=t^2\,\int_\Sigma\frac{d}{d\tau}\lf[\Re<{Q},\vec{h}^0(\tau)>_{g(\tau)}\rg](0)\cdot\frac{d\vec{\Phi}}{d\tau}(0)\ dvol_{g}\\[5mm]
\ds+2\, \int_0^t(t-\tau)\,d\tau\,\int_0^\tau\ ds\ \frac{d}{ds}\lf[\int_\Sigma\frac{d}{ds}\lf[\Re<{Q},\vec{h}^0(s)>_{g(s)}\rg]\cdot\frac{d\vec{\Phi}}{ds}\ dvol_{g(s)}\rg]\\[5mm]
\ds-2\, \int_0^t(t-\tau)\,d\tau\,\int_\Sigma\Re\lf[<{Q},\vec{h}^0>_{g(0)}-<Q,\vec{h}^0(\tau)>_{g(\tau)}\rg]\cdot\frac{d}{d\tau}\lf[\frac{d\vec{\Phi}}{d\tau}\ dvol_{g(\tau)}\rg]
\end{array}
\ee
Combining (\ref{0-I.4}) and (\ref{0-I.7}) gives then
\be
\label{0-I.8}
\begin{array}{l}
\ds A(\vec{\Phi}_{x^0}^\ep(t))=A(\vec{\Phi})-\, t^2\, \int_\Sigma\frac{d\vec{H}}{d\tau}(0)\cdot\frac{d\vec{\Phi}}{d\tau}(0)\ dvol_{g}\\[5mm]
\ds+t^2\,\int_\Sigma\frac{d}{d\tau}\lf[\Re<{Q},\vec{h}^0(\tau)>_{g(\tau)}\rg](0)\cdot\frac{d\vec{\Phi}}{d\tau}(0)\ dvol_{g}\\[5mm]
\ds-2\, \int_0^t(t-\tau)\,d\tau\,\int_0^\tau\ ds\ \int_\Sigma\frac{d^2\vec{H}}{ds^2}\cdot\frac{d\vec{\Phi}}{ds}\ dvol_{g(s)}\\[5mm]
\ds-2\, \int_0^t(t-\tau)\,d\tau\,\int_0^\tau\ ds\ \int_\Sigma\frac{d\vec{H}}{ds}\cdot\frac{d}{ds}\lf[\frac{d\vec{\Phi}^\ep_{x^0}}{ds}\ dvol_{g(s)}\rg]\\[5mm]
\ds-2\, \int_0^t(t-\tau)\,d\tau\,\int_\Sigma[\vec{H}(\tau)-\vec{H}(0)]\cdot\frac{d}{d\tau}\lf[\frac{d\vec{\Phi}}{d\tau}\ dvol_{g^\ep_{x^0}(\tau)}\rg]\\[5mm]
\ds+2\, \int_0^t(t-\tau)\,d\tau\,\int_0^\tau\ ds\ \int_\Sigma\frac{d^2}{ds^2}\lf[\Re<{Q},\vec{h}^0(s)>_{g(s)}\rg]\cdot\frac{d\vec{\Phi}}{ds}\ dvol_{g(s)}\\[5mm]
\ds+2\, \int_0^t(t-\tau)\,d\tau\,\int_0^\tau\ ds\ \int_\Sigma\frac{d}{ds}\lf[\Re<{Q},\vec{h}^0(s)>_{g(s)}\rg]\cdot\frac{d}{ds}\lf[\frac{d\vec{\Phi}}{ds}\ dvol_{g(s)}\rg]\\[5mm]
\ds-2\, \int_0^t(t-\tau)\,d\tau\,\int_\Sigma\Re\lf[<{Q},\vec{h}^0>_{g(0)}-<Q,\vec{h}^0(\tau)>_{g(\tau)}\rg]\cdot\frac{d}{d\tau}\lf[\frac{d\vec{\Phi}}{d\tau}\ dvol_{g(\tau)}\rg]
\end{array}
\ee
Using (\ref{xsIII.26}), (\ref{xsIII.29}) and (\ref{xsIII.31-zz}) we obtain 
\be
\label{0-I.9}
\begin{array}{l}
\ds -\, 2\, \int_\Sigma\lf[\frac{d\vec{H}}{d\tau}(0)-\frac{d}{d\tau}\lf[\Re<{Q},\vec{h}^0(\tau)>_{g(\tau)}\rg](0)\rg]\cdot\frac{d\vec{\Phi}}{d\tau}(0)\ dvol_{g}\\[5mm]
\ds\quad=\int_{\Sigma}\lf[|dv|_g^2-2\,
<\ov{q},dv\otimes dv>_g +8\,H^2\, v^2-(|\vec{\mathbb I}|^2_g+2)\, v^2\rg]\ dvol_g\quad.\\[5mm]
\ds\quad+\int_\Sigma\lf[2H\, dv+v\,dH+4\,v\,\Re(Q)\res dH+\Im\lf<Q,h^0\rg>\, \ast dv   \rg]\cdot \vec{w}_T\ dvol_g\\[5mm]
\ds\quad-\int_\Sigma\Im\lf<Q,h^0\rg>\ v\, \ast d\vec{w}_T^{\,\ast}\ dvol_g\quad.
\end{array}
\ee
where $v=\frac{d\vec{\Phi}^\ep_{x^0}}{d\tau}(0)\cdot\vec{n}$. Using the explicit formula of $\vec{\Phi}_{x^0}^\ep(t)$ we obtain
\be
\label{0-I.10}
v=\beta_{x^0}^\ep(x)\ \chi_{x^0}^\ep(x)\ \vec{n}(x^0)\cdot \vec{n}+\beta_{x^0}^\ep(x)\ \sum_{j=1}^2\al_j(0,\ep)\ \vec{a}_j\cdot\vec{n}
\ee
and
\[
|\vec{w}_T|=|\vec{w}\wedge\vec{n}|\le C\,\ep\, |\vec{n}-\vec{n}(x^0)|\ {\mathbf 1}_{x^0}^\ep+O_{\vec{\Phi},x^0}(\ep^3)
\]
Moreover
\[
|\nabla^g\vec{w}_T|\le C\, \lf(|\vec{n}-\vec{n}(x^0)|\ +\ep\, |\vec{\mathbb I}|_g\rg)\ {\mathbf 1}_{x^0}^\ep+O_{\vec{\Phi},x^0}(\ep^3)\,|\vec{\mathbb I}|_g
\]
Inserting the previous expressions in (\ref{0-I.9}) gives then
\be
\label{0-I.11}
\begin{array}{l}
\ds   -\, 2\, \int_\Sigma\lf[\frac{d\vec{H}}{d\tau}(0)-\frac{d}{d\tau}\lf[\Re<{Q},\vec{h}^0(\tau)>_{g(\tau)}\rg](0)>_{g(\tau)}\rg]\cdot\frac{d\vec{\Phi}^\ep_{x^0}}{d\tau}(0)\ dvol_{g}\\[5mm]
\ds\quad=\ep^2\ \int_{{\C}}|\nabla\chi|^2\ dx^2-\,4\,\ep^2\, e^{-2\la(x^0)}\, Q_1(x^0)\int_{\C}|\p_{x_1}\chi|^2-|\p_{x_2}\chi|^2\ dx^2\\[5mm]
\ds\quad\quad\quad\quad
-\,4\,\ep^2\, e^{-2\la(x^0)}\, Q_2(x^0)\int_{\C}2\,\p_{x_1}\chi\,\p_{x_2}\chi\ dx^2+o_{\vec{\Phi},x^0}(\ep^2)
\end{array}
\ee
We proceed to a rotation of the domain in such a way that $Q_2(x^0)=0$. Combining (\ref{0-I.8}), the estimates in lemma~\ref{lm-xsIII.1a} and (\ref{0-I.11}), one has finally
\be
\label{0-I.12}
\begin{array}{l}
\ds\lf|A(\vec{\Phi}_{x^0}^\ep(t))-A(\vec{\Phi})-{t^2\,\ep^2}\ F_\chi(x_0)\rg|\le C_{\vec{\Phi},x^0}\ \ep\ t^3+ t^2\ o_{\vec{\Phi},x^0}(\ep^2)\quad,
\end{array}
\ee
where $B(x^0)$ is the number given by
\be
\label{0-I.12a}
2\, F_\chi(x^0):=\int_{{\C}}|\nabla\chi|^2\ dx^2-\,4\, e^{-2\la(x^0)}\, Q_1(x^0)\int_{\C}|\p_{x_1}\chi|^2-|\p_{x_2}\chi|^2\ dx^2\quad.
\ee
For any $\vec{a}\in B^4$ we denote by $\Psi_{\vec{a}}$ the following element of the M\"obius group
\[
\Psi_{\vec{a}}(\vec{y})=(1-|\vec{a}|)\frac{\vec{y}-\vec{a}}{|\vec{y}-\vec{a}|^2}-\vec{a}\quad.
\]
Its conformal factor is given by 
\[
\forall\ \vec{Y}\in T_{\vec{y}}S^3\quad\quad |d\Psi_{\vec{a}}\cdot \vec{Y}|=\frac{1-|\vec{a}|}{1+|\vec{a}|^2-2\vec{a}\cdot \vec{y}}\ |\vec{Y}|\quad.
\]
Hence we have
\[
A(\Psi_{\vec{a}}\circ\vec{\Phi})=\int_{\Sigma}\frac{(1-|\vec{a}|)^2}{(1+|\vec{a}|^2-2\,\vec{a}\cdot \vec{\Phi})^2}\ dvol_g\quad.
\]
For $(t,\ep)$ in the neighborhood of the origin ${\mathfrak U}_{x^0,\vec{\Phi}}$ given by lemma~\ref{lm-xsIII.1a} we estimate for $|\vec{a}|<\nu<1$ the following difference
\be
\label{0-I.13}
\begin{array}{l}
\ds A(\Psi_{\vec{a}}\circ\vec{\Phi}_{x^0}^\ep(t))-A(\Psi_{\vec{a}}\circ\vec{\Phi})\\[5mm]
\ds\quad=(1-|\vec{a}|)^2\lf[\int_\Sigma\frac{dvol_{g_{x^0}^\ep(t)}}{(1+|\vec{a}|^2-2\,\vec{a}\cdot \vec{\Phi}_{x^0}^\ep(t))^2}-\frac{dvol_{g}}{(1+|\vec{a}|^2-2\,\vec{a}\cdot \vec{\Phi})^2}\rg]\\[5mm]
\ds\quad=(1-|\vec{a}|)^2\lf[\int_\Sigma\frac{dvol_{g_{x^0}^\ep(t)}-dvol_{g}}{(1+|\vec{a}|^2-2\,\vec{a}\cdot \vec{\Phi}_{x^0}^\ep(t))^2}\rg]\\[5mm]
\ds\quad\  +(1-|\vec{a}|)^2\lf[\int_\Sigma\lf(\frac{1}{(1+|\vec{a}|^2-2\,\vec{a}\cdot \vec{\Phi}_{x^0}^\ep(t))^2}-\frac{1}{(1+|\vec{a}|^2-2\,\vec{a}\cdot \vec{\Phi})^2}\rg)\ dvol_g\rg]
\end{array}
\ee
We have in one hand
\be
\label{0-I.14}
\begin{array}{l}
\ds \int_\Sigma\frac{dvol_{g_{x^0}^\ep(t)}-dvol_{g}}{(1+|\vec{a}|^2-2\,\vec{a}\cdot \vec{\Phi}_{x^0}^\ep(t)(x))^2}\\[5mm]
\ds\quad=\frac{1}{(1+|\vec{a}|^2-2\,\vec{a}\cdot \vec{\Phi}_{x^0}^\ep(t)(x^0))^2} \int_\Sigma dvol_{g_{x^0}^\ep(t)}-dvol_{g}\\[5mm]
\ds-\int_\Sigma\lf[\frac{dvol_{g_{x^0}^\ep(t)}-dvol_{g}}{(1+|\vec{a}|^2-2\,\vec{a}\cdot \vec{\Phi}_{x^0}^\ep(t)(x^0))^2}-\frac{dvol_{g_{x^0}^\ep(t)}-dvol_{g}}{(1+|\vec{a}|^2-2\,\vec{a}\cdot \vec{\Phi}_{x^0}^\ep(t)(x))^2}\rg]\ 
 \end{array}
\ee
Hence, using the estimate (\ref{xsIII.34a-z19p}), we obtain for $|\vec{a}|<1/4$
\be
\label{0-I.15}
\begin{array}{l}
\ds \lf|\int_\Sigma\frac{dvol_{g_{x^0}^\ep(t)}-dvol_{g}}{(1+|\vec{a}|^2-2\,\vec{a}\cdot \vec{\Phi}_{x^0}^\ep(t)(x))^2}-\frac{dvol_{g_{x^0}^\ep(t)}-dvol_{g} }{(1+|\vec{a}|^2-2\,\vec{a}\cdot \vec{\Phi}_{x^0}^\ep(t)(x^0))^2}\rg|\\[5mm]
\ds\le C_{\vec{\Phi},x^0}\ t\ \int_{\Sigma}\lf[  {\mathbf 1}^\ep_{x^0}(|\vec{n}-\vec{n}^0|+t)+ \ep^3+t\,\ep^2\rg]\ |\vec{a}|\ |\vec{\Phi}_{x^0}^\ep(t)(x)-\vec{\Phi}_{x^0}^\ep(t)(x^0)|\ dvol_g\\[5mm]
\ds\le 4\ |\vec{a}|\ C_{\vec{\Phi},x^0}\,\int_{\Sigma} [t\, \ep^3+t^2\,\ep^2]\ dvol_g+\, |\vec{a}|\ \ep\ t\ C_{\vec{\Phi}}\ \int_{B_\ep(x^0)} (|\vec{n}-\vec{n}^0|+t)\ dvol_g\\[5mm]
\ds\le \, C_{\vec{\Phi},x^0}\  |\vec{a}|\ (t\ \ep^3+t^2\ \ep^2)
\end{array}
\ee
In the other hand we have, using the explicit expression of $\vec{\Phi}_{x^0}^\ep(t)$ and the estimate (\ref{xsIII.34a-z1}),
\be
\label{0-I.16}
\begin{array}{l}
\ds\lf|\int_\Sigma\lf(\frac{1}{(1+|\vec{a}|^2-2\,\vec{a}\cdot \vec{\Phi}_{x^0}^\ep(t))^2}-\frac{1}{(1+|\vec{a}|^2-2\,\vec{a}\cdot \vec{\Phi})^2}\rg)\ dvol_g\rg|\\[5mm]
\ds\quad\le |\vec{a}|\ \lf[\int_{B_\ep(x^0)}\,\ep\ t\,\|\chi\|_\infty dvol_g+ C_{\vec{\Phi},x^0}\int_\Sigma t\ \ep^3\ dvol_g\rg]\le \, C_{\vec{\Phi},x^0}\, |\vec{a}|\ t\, \ep^3
\end{array}
\ee
Combining (\ref{0-I.12}), (\ref{0-I.13}), (\ref{0-I.14}), (\ref{0-I.15}) and (\ref{0-I.16}) we obtain for $|\vec{a}|<1/4$ and $(t,\ep)\in{\mathfrak U}_{x^0,\vec{\Phi}}$
\be
\label{0-I.17}
\begin{array}{l}
\ds\lf| A(\Psi_{\vec{a}}\circ\vec{\Phi}_{x^0}^\ep(t))-A(\Psi_{\vec{a}}\circ\vec{\Phi})-\frac{t^2\,\ep^2\, F_\chi(x^0)}{(1+|\vec{a}|^2-2\,\vec{a}\cdot \vec{\Phi}_{x^0}^\ep(t)(x^0))^2}\rg|\\[7mm]
\ds\quad\le\ C_{\vec{\Phi},x^0}\ \ep\ t^3+ t^2\ o_{\vec{\Phi},x^0}(\ep^2)+C_{\vec{\Phi},x^0}\  |\vec{a}|\ (t\ \ep^3+t^2\ \ep^2)
\end{array}
\ee
Let now $\varphi(s)\in C^\infty_0({\R})$ be an arbitrary non zero function compactly supported on ${\R}$. For any $\tau\in{\R}_+^\ast$ we denote
\[
\chi_\tau(x_1,x_2):=\,\varphi(\tau\, x_1)\ \varphi(x_2)
\] 
We have
\[
F_{\chi_\tau}(x^0)= \int_{\R}\dot{\varphi}^2(s)\, ds\, \int_{\R}{\varphi}^2(s)\, ds\, \lf[ \frac{\tau}{2}\, (1- 4\, e^{-2\la(x^0)}\, Q_1(x^0))+\frac{1}{2\,\tau}\, (1+ 4\, e^{-2\la(x^0)}\, Q_1(x^0))\rg]
\]
Assuming 
\[
|4\, e^{-2\la(x^0)}\, Q_1(x^0)|>1
\]
which corresponds to assume $4\,|Q|^2_g>1$ since $Q_2=0$ and $|dz^2|^2_g=4\, e^{-4\la}$, we can choose $\tau\in{\R}^\ast_+$ such that $F_{\chi_\tau}(x^0)<0$. Hence by choosing $t=\delta_{\vec{\Phi},x^0}\ \ep$ for $\ep<\ep^0_{\vec{\Phi},x^0}$  and $\ep^0_{\vec{\Phi},x^0}$ is chosen small enough and $|\vec{a}|<\nu_{{\vec{\Phi},x^0}}$ for some
$\nu_{{\vec{\Phi},x^0}}>0$, also chosen small enough but independent of $\ep$, we obtain from (\ref{0-I.17}) the lemma~\ref{lm-0-I.2}.\hfill $\Box$

\section{Proof of theorem~\ref{th-0-2}.}
\reset

Let $\vec{\Phi}$ be a weak immersion of the torus $T^2$ in $S^3$ and assume that $\vec{\Phi}$ is a local minimizer of the conformal volume in it's conformal class. From section II we deduce that
$\vec{\Phi}$ is a conformally constrained minimal surface. If $\vec{\Phi}$ is additionally isothermic then we apply section IV and the theorem is proved in that case. If now $\vec{\Phi}$ is not isothermic
and posses a conformally congruent immersion which is not isometric to $\vec{\Phi}$ then lemma~\ref{lm-IV-3} implies the theorem in that case. If now $\vec{\Phi}$ satisfies
\[
\forall \ \vec{a}\in B^4(0)\setminus\{ 0\}\quad \quad A(\vec{\Phi})> A(\Psi_{\vec{a}}\circ\vec{\Phi})\quad,
\]
then, since $A(\vec{\Phi})\notin 4\pi{\N}$, for all $\nu>0$ there exists $\delta>0$ such that
\be
\label{cc-bb}
\forall \ \vec{a}\in B^4(0)\setminus B^4_\nu(0)\quad \quad A(\vec{\Phi})> A(\Psi_{\vec{a}}\circ\vec{\Phi})+\delta\quad.
\ee
Assume there would exists a non zero measure set of points $x^0\in T^2$ such that $2\ |Q|_{g_{\vec{\Phi}}}(x^0)>1$ then, using lemma~\ref{lm-0-I.2}, there exists $\ep_0>0$, $\nu_0>0$ and a family of deformations $\vec{\Phi}^\ep$ in the same conformal class defined by $\vec{\Phi}$ such that such that $\vec{\Phi}^\ep$ converges to $\vec{\Phi}$ as $\ep\rightarrow 0$ in the space of weak immersions and
\be
\label{0-I.18-cb}
\begin{array}{l}
\forall \,\ep<\ep_{0}\quad\forall\,\vec{a}\in B_{\nu_{0}}^4(0)\quad
A(\Psi_{\vec{a}}\circ\vec{\Phi}^\ep)< A(\Psi_{\vec{a}}\circ\vec{\Phi})-\, \ep^4\ C_0
\end{array}
\ee
for some positive constant $C_0>0$ which is independent of $\ep$. Let $\nu:=\nu_0$ and let $\delta_0>0$ such that (\ref{cc-bb}) holds. Then there exists $\ep_1>0$ such that
\be
\label{cc-bb-1}
\forall \ \vec{a}\in B^4(0)\setminus B^4_{\nu_0}(0)\quad \forall\ep<\ep_1\quad\quad A(\vec{\Phi})> A(\Psi_{\vec{a}}\circ\vec{\Phi}^\ep)+\delta/2\quad.
\ee
Hence for any $0<\ep<\min\{\ep_0,\ep_1\}$ we would have
\[
V_c(\vec{\Phi}^\ep)< V_c(\vec{\Phi})\quad.
\]
This contradicts the assumption that $\vec{\Phi}$ locally minimizes the conformal volume in it's conformal class and theorem~\ref{th-0-2} is proved in all cases.

\appendix
\section{Appendix}
\reset

\subsection{Codazzi identity.}

Recall the definition of the Weingarten form of the immersion $\vec{\Phi}$
\[
\vec{h}^0:=2\, \pi_{\vec{n}}(\p^2_{z^2}\vec{\Phi})\ dz\otimes dz=e^{2\la}\ \vec{H}^0\ dz\otimes dz
\]
where
\[
\vec{H}^0=2\,\p_z\lf(e^{-2\la}\p_z\vec{\Phi}\rg)=2^{-1} e^{-2\la}\ \pi_{\vec{n}}(\p^2_{x^2_1}\vec{\Phi}-\p^2_{x^2_2}\vec{\Phi})-
i\ e^{-2\la}\ \pi_{\vec{n}}(\p^2_{x_1x_2}\vec{\Phi})
=[H^0_{\Re}+i\,H^0_{\Im}]\ \vec{n}\quad,
\]
where $\p_z=2^{-1}(\p_{x_1}-i\,\p_{x_2})$. This gives
\be
\label{xsIII.17}
\lf\{
\begin{array}{l}
\ds H^0_{\Re}=-\frac{e^{-2\la}}{2}\,\lf[(\p_{x_1}\vec{n},\p_{x_1}\vec{\Phi})-(\p_{x_2}\vec{n},\p_{x_2}\vec{\Phi})\rg]=e^{-2\la}{\mathbb I}^0_{11}=-e^{-2\la}{\mathbb I}^0_{22}\\[5mm]
\ds H^0_{\Im}=e^{-2\la}\,(\p_{x_1}\vec{n},\p_{x_2}\vec{\Phi})=e^{-2\la}\,(\p_{x_2}\vec{n},\p_{x_1}\vec{\Phi})=-e^{-2\la}{\mathbb I}^0_{12}\quad.
\end{array}
\rg.
\ee
We have
\[
\begin{array}{l}
\p_{\ov{z}}\vec{H}^0=2\,\p_{\ov{z}}\p_z\lf(e^{-2\la}\,\p_z\vec{\Phi}\rg)=2\,\p_{{z}}\p_{\ov{z}}\lf(e^{-2\la}\,\p_z\vec{\Phi}\rg)\\[5mm]
\quad=-4\,\p_z\lf(e^{-2\la}\,\p_{\ov{z}}\la\,\p_z\vec{\Phi}\rg)+2\,\p_z\lf(e^{-2\la}\,\p^2_{z\ov{z}}\vec{\Phi}\rg)=-4\,\p_z\lf(e^{-2\la}\,\p_{\ov{z}}\la\,\p_z\vec{\Phi}\rg)+2^{-1}\,\p_z\lf(e^{-2\la}\,\Delta\vec{\Phi}\rg)
\end{array}
\]
Hence we have obtained the following identity so far
\be
\label{xsIII.18}
\p_{\ov{z}}\vec{H}^0=-4\,\p_z\lf(e^{-2\la}\,\p_{\ov{z}}\la\,\p_z\vec{\Phi}\rg)+\p_z\vec{H}-\p_z\vec{\Phi}\quad.
\ee
Taking the scalar product with $\vec{n}$ gives then
\be
\label{xsIII.19}
\p_{\ov{z}}H^0=-2\,\p_{\ov{z}}\la\,H^0+\p_zH
\ee
from which we deduce
\be
\label{xsIII.20}
\p_{\ov{z}}\lf(e^{2\la}\,H^0\rg)= e^{2\la}\ \p_zH\quad.
\ee
This can also be rewritten locally as follows
\be
\label{xsIII.20.0}
\lf\{
\begin{array}{l}
\ds\p_{x_1}{\mathbb I}^0_{11}+\p_{x_2}{\mathbb I}^0_{12}=e^{2\la}\ \p_{x_1}H\\[5mm]
\ds\p_{x_2}{\mathbb I}^0_{11}-\p_{x_1}{\mathbb I}^0_{12}=-e^{2\la}\ \p_{x_2}H
\end{array}
\rg.
\ee
Hence we have proved the following
\begin{Lma}
\label{Codazzi}
Let $h^0:=2\,\vec{n}\cdot\p^2_{z^2}\vec{\Phi}\,dz^2$ and denote $g_{{\C}}:=e^{2\la}\ d\ov{z}\otimes dz$ we have
\be
\label{xsIII.21}
\ov{\p}h^0=g_{{\C}}\otimes\p H\quad.
\ee
\hfill $\Box$
\end{Lma}
\subsection{Construction of infinitesimal perturbations within a conformal class.}

 Recall that ${\mathcal T}_{T^2}$ denotes the
Teichm\"uller Space of $T^2$ and, having fixed generators of the $\pi_1(T^2)$, for any metric $g$ on $T^2$ we denote
by $[g]$ the Teichm\"uller class associated to $g$. 
Let $\vec{\Phi}$ be a weak immersion of $T^2$ and assume it is \underbar{not isothermic}. Using \cite{Ri3}, we decuce the existence of 
 two  maps $\vec{a}_j$  in the space $W^{2,2}\cap W^{1,\infty}(T^2,{\R}^4)$ such that
\[
(d{\mathcal C}_{\vec{\Phi}}\cdot\vec{a}_j)_{j=1,2}\quad\quad\mbox{ forms a basis to the tangent space at $[\vec{\Phi}^\ast g_{{\R}^3}]$ of }{\mathcal T}_{T^2}
\]
and $\vec{a}_j(x)\cdot\vec{\Phi}(x)\equiv 0$ on $T^2$. %Because of the explicit formula of $d{\mathcal C}$ given in \cite{Ri3} we ca assume that $\vec{a}_j=a_j\,\vec{n}$.
\begin{Lma}
\label{lm-xsIII.1a}
Let $x_0$ be a point in $T^2$ which is a Lebesgue
point for $\nabla^2\vec{\Phi}$ and for $\nabla\vec{\Phi}$. Let $\chi$ be a smooth function on ${\C}$ supported in $B_1(0)$. In some fixed conformal chart
in a neighborhood of $x^0$ we denote $\chi_{x^0}^\ep(x):=\ep\,\chi(\ep^{-1}(x-x^0))$. There exists a neighborhood of $0$ in ${\R}^2$, ${\mathfrak U}_{x^0,\vec{\Phi}}$, and  there exists two $C^1$ maps $\sigma_j(t,\ep)$ such that
\be
\label{xsIII.32a}
\begin{array}{l}
\ds\forall\  (t,\ep)\in{\mathfrak U}_{x^0,\vec{\Phi}} \quad\quad{\mathcal C}\lf(\vec{\Phi}_{x^0}^\ep(t)\rg)\equiv{\mathcal C}(\vec{\Phi})\\[5mm]
\mbox{where }\quad\vec{\Phi}_{x^0}^\ep(t):=\beta^\ep(x,t)\,\lf[\vec{\Phi}(x)+t\ \chi_{x^0}^\ep(x)\,\vec{n}(x_0)+t\ \sum_{j=1}^2\al_j(t,\ep)\, \vec{a}_j(x)\rg]\quad,
\end{array}
\ee
\be
\label{xsIII.33a}
|\al_j(0,\ep)|=O_{\vec{\Phi},x^0}(\ep^3)\quad,
\ee	
\be
\label{xsIII.34a}
\beta^\ep(x,t) \mbox{ is chosen such that }|\vec{\Phi}_{x^0}^\ep(t)|\equiv 1
\ee
Moreover $\al_j(t,\ep)$ is $C^2$ with respect to $t$ for $\ep>0$. Finally the following estimates hold for all  $ (t,\ep)\in{\mathfrak U}_{x^0,\vec{\Phi}}$ and $k=0,1,2$
\be
\label{xsIII.34a-z1}
\sum_{j=1}^2\left|\frac{d^k}{dt^k}\lf(t\,\al_j(t,\ep)\rg)\right|\le\,  C_{\vec{\Phi},x^0}\,\ep^2\,\lf[\inf\{\ep\,|t|^{1-k}, 1\}+|t|^{2-k}\rg]
\ee
and
\be
\label{xsIII.34a-z19pp}
\ds\lf|\frac{d}{dt}\nabla\vec{\Phi}_{x^0}^\ep(t)\rg|\le C_{\vec{\Phi}}\ \lf[{\mathbf 1}^\ep_{x^0}+ C_{\vec{\Phi},x^0}\,(\ep^3+|t|\,\ep^2)   \rg]\quad .
\ee
and
\be
\label{xsIII.34a-z19p}
\lf|\frac{d g_{\vec{\Phi}_{x^0}^\ep(t)}}{dt}\rg|\le C_{\vec{\Phi}}\,\lf[  {\mathbf 1}^\ep_{x^0}(|\vec{n}-\vec{n}^0|+t)+ O_{x^0,\vec{\Phi}}(\ep^3+|t|\,\ep^2) \rg]
\ee
and
\be
\label{xsIII.34a-z22p}
\begin{array}{l}
\ds\lf|\frac{d}{dt}\nabla^2\vec{\Phi}_{x^0}^\ep(t)\rg|+\lf|\frac{d}{dt}{\mathbb I}_{\vec{\Phi}_{x^0}^\ep(t)}\rg|\le C_{\vec{\Phi}}\ {\mathbf 1}^\ep_{x^0}\ \lf[ \ep^{-1}+ \ep\, |\nabla^2\vec{\Phi}|\rg]\\[5mm]
\ds\quad\quad\quad+C_{\vec{\Phi},x^0}\,(\ep^3+|t|\,\ep^2)\, \lf[1+\sum_{j=1}^2|\nabla^2 \vec{a}_j|\rg] \quad .
\end{array}
\ee
and finally
\be
\label{xsIII.34a-z30p}
\ds\lf|\frac{d^2}{dt^2}\nabla^2\vec{\Phi}_{x^0}^\ep(t)\rg|\le\ C_{\vec{\Phi},x^0}\   {\mathbf 1}^\ep_{x^0}+  C_{\vec{\Phi},x^0}\,[1+|\nabla^2\vec{\Phi}|]\ \,(\ep^2+|t|\,\ep)  \quad ,
\ee
where $C_{\vec{\Phi}}>0$ is independent of $(t,\ep)\in{\mathfrak U}_{x^0,\vec{\Phi}}$,  where ${\mathbf 1}^\ep_{x^0}$ is the characteristic function of the geodesic ball
of center $x^0$ and radius $\ep>0$ for $g_0$, the constant scalar curvature metric of volume 1, and $C_{\vec{\Phi},x^0}$ is a constant depending only on $\chi$ , on
\[
\sup_{\ep>0} \ep^{-2}\int_{B_\ep(x^0)}[|{\mathbb I}_{\vec{\Phi}}|_g^2\ +1]\ dvol_g
\]
and on $\|\nabla a_j\|_{L^\infty(\Sigma)}+\|\log|\nabla\vec{\Phi}|\|_{L^\infty(\Sigma)}$. \hfill $\Box$
\end{Lma}
\noindent{\bf Proof of Lemma~\ref{lm-xsIII.1a}.}
We Assume 
\be
\label{xsIII.36a}
\ep<\ep_0\quad\quad\mbox{ and }\quad\quad |t|\le 4^{-1}\, \|e^{\la}\|_{L^\infty(B_{\ep_0}(x^0))}
\ee 
for some $\ep_0>0$, where the ball $B_\ep(x^0)$ is the ball of radius $\epsilon$ for the flat metric of volume one $g_0$ conformally 
equivalent to $g_{\vec{\Phi}}$. We use a chart in a neighborhood of $x^0$ in which $g_0=dx_1^2+dx_2^2$ and $g_{\vec{\Phi}}=e^{2\la}\,g_0$. $\ep_0$ has been then chosen sufficiently small in order for this chart $x=(x_1,x_2)$ to contain the support of $\chi_{x^0}^\ep(x)$.

\medskip

\noindent Denote $$\vec{n}^0:=\vec{n}(x^0)\quad.$$

\medskip

\noindent Under these assumptions, for $\ep_0$ small enough, the perturbation of $\vec{\Phi}$
\[
\vec{\Phi}+ t\, \chi^\ep_{x^0}\ \vec{n}(x_0)
\]
is still an immersion. Consider now $\al_j\in {\R}$ such that 
\be
\label{xsIII.37a}
\sum_{j=1}^2|\al_j|\,\|\vec{a}_j\|_{L^\infty(\Sigma)}<4^{-1} \|e^{-\la}\|_\infty\quad.
\ee
Under the assumptions (\ref{xsIII.36a}) and (\ref{xsIII.37a}), the map
\[
\vec{\Phi}^\ep_{x^0}(t,\al_1,\al_2):=\beta(t,\ep,\al_1,\al_2)(x)\,\lf[\vec{\Phi}+ t\, \chi^\ep_{x^0}\ \vec{n}^0+t \sum_{j=1}^2\al_j\ \vec{a}_j\rg]
\]
is again an immersion where $\beta$ is chosen in such a way that $|\vec{\Phi}^\ep_{x^0}|\equiv 1$ on $T^2$ :

\be
\label{xsIII.37b}
\begin{array}{l}
\ds\beta(t,\ep,\al_1,\al_2)(x)=\lf[1+2\,t\, \vec{\Phi}\cdot\vec{n}^0\,\chi^\ep_{x^0}+t^2(\chi^\ep_{x^0})^2+\lf[\sum_{j=1}^2t\, \al_j\,\vec{a}_j\rg]^2\rg.\\[5mm]
\ds\quad\quad\quad\quad\quad\quad\quad+\lf.2\, t\,\chi^\ep_{x^0}\lf[\sum_{j=1}^2t\,\al_j\,\vec{a}_j\cdot\vec{n}^0\rg]\rg]^{-1/2}
\end{array}
\ee

%\be
%\label{xsIII.37b}
%\begin{array}{l}
%\ds\beta(t,\ep,\al_1,\al_2)(x)=\lf[1+2\,t\, \vec{\Phi}\cdot\vec{n}^0\,\chi^\ep_{x^0}+t^2\ \lf[\sum_{j=1}^2\al_j\,a_j+\chi^\ep_{x^0}\rg]^2\rg.\\[5mm]
%\ds\quad\quad\quad\quad\quad\quad\quad+\lf.t^2\, [1-\vec{n}\cdot\vec{n}^0]\ \lf[\sum_{j=1}^2\al_j\,a_j+\chi^\ep_{x^0}\rg]\rg]^{-1/2}
%\end{array}
%\ee

 For $(t,\ep,\al_1,\al_2)$ in the open neighborhood ${\mathcal U}_{\vec{\Phi}}$ of $(0,0,0,0)$  given by (\ref{xsIII.36a}) and (\ref{xsIII.37a})
we define $\Gamma(t,\ep,\al_1,\al_2)$ by
\[
\forall\ (t,\ep,\al_1,\al_2)\in {\mathcal U}_{\vec{\Phi}}\quad,\quad t\ne 0\quad \Gamma(t,\ep,\al_1,\al_2):={\mathcal C}(\vec{\Phi}^\ep_{x^0}(t,\al_1,\al_2))
\]
For any $\ep<\ep_0$ we denote by $(\al_1^0(\ep),\al_2^0(\ep))$ the unique element of ${\R}^2$ such that 
\[
d{\mathcal C}_{\vec{\Phi}}\cdot \lf[\chi^\ep_{x^0}\ \vec{n}^0+\sum_{j=1}^2\al_j^0(\ep)\ \vec{a}_j\rg]=0
\]
We are using in this assertion the fact that $d{\mathcal C}_{\vec{\Phi}}\cdot\vec{a}_j$ realizes a basis of the tangent space to the Teichm\"uller space ${\mathcal T}_{T^2}$
at ${\mathcal C}(\vec{\Phi})$. 
Let $(Q^1,Q^2)$ be a pair of 2 independent holomorphic quadratic forms of $(T^2, [\vec{\Phi}^\ast g_{{\R}^3}])$, chosen to be orthonormal 
with respect to the Weil-Peterson metric. One has (see \cite{Ri3}), denoting by $g_0$ the flat metric of volume one on 
$(T^2, [\vec{\Phi}^\ast g_{{\R}^3}])$
\[
\begin{array}{l}
\ds d{\mathcal C}_{\vec{\Phi}}\cdot \lf[\chi^\ep_{x^0}\ \vec{n}\rg]=\sum_{j=1}^2Q^j\ \int_{\Sigma}\chi^\ep_{x^0}\ \Re<Q^j,h^0>_{wp, g_0}\ \vec{n}\cdot\vec{n}(x_0)dvol_{g_0}\\[5mm]
\ds\quad\quad=\sum_{j=1}^2Q^j\ \int_{B_\ep(x^0)}\chi^\ep_{x^0}\ <\ov{q}^j,{\mathbb I}^0_{\vec{\Phi}}>_{g_{\vec{\Phi}}}\ \vec{n}\cdot\vec{n}(x_0)\ dvol_{g_{\vec{\Phi}}}
\end{array}
\]
where $\ov{q}^j$ is the following trace free symmetric real 2-differential given in any conformal coordinates by
\[
\ov{q}^j:=2\,Q^j_1\, [dx_1^2-dx_2^2]-2\, Q^j_2\, [dx_1\, dx_2+dx_2\, dx_1]\quad\mbox{ where }\quad Q^j= (Q^j_1+i\,Q^j_2)\ dz^2\quad.
\]
and then in any conformal chart one has (see previous sections)
\[
<\ov{q}^j,{\mathbb I}^0_{\vec{\Phi}}>_{g_{\vec{\Phi}}}= 4\, e^{-4\la}\,[Q^j_1\ {\mathbb I}^0_{11}-Q^j_2\ {\mathbb I}^0_{12}]
\]
where $g_{\vec{\Phi}}=e^{2\la}\ [dx_1^2+dx_2^2]$ in these coordinates. Therefore we have
\be
\label{xsIII.38t}
\lf|d{\mathcal C}_{\vec{\Phi}}\cdot \lf[\chi^\ep_{x^0}\ \vec{n}(x_0)\rg]\rg|\le C\  \|e^{-2\la}\|_\infty\,\ep^2\ \lf[\int_{B_\ep(x^0)}|{\mathbb I}^0_{\vec{\Phi}}|^2\rg]^{1/2}\le O_{\vec{\Phi},x^0}( \ep^3)\quad,
\ee
Since $x_0$ is chosen to be a Lebesgue point for the second fundamental form we have
\be
\label{xsIII.38s}
\lim_{\ep\rightarrow 0}\ep^{-2}\int_{B^{g_0}_\ep(x^0)}|{\mathbb I}^0_{\vec{\Phi}}|^2 dvol_{g_{\vec{\Phi}}}=\pi\ |{\mathbb I}^0_{\vec{\Phi}}|^2(x_0)\ e^{2\la(x_0)}
\ee
So we deduce from (\ref{xsIII.38t}) and (\ref{xsIII.38s})
\be
\label{xsIII.39-v}
|\al_j^0|=O_{\vec{\Phi},x^0}(\ep^3)\quad.
\ee
We claim that the mapping $\Gamma$ is $C^1$ in ${\mathcal U}_{\vec{\Phi}}$. The $C^1$ property of $\Gamma$ at $(t,0,\al_1,\al_2)$ requires a justification : We have
\[
\p_\ep\Gamma(t,\ep,\al_1,\al_2)=d{\mathcal C}_{\vec{\Phi}_{x^0}^\ep} \lf[\p_\ep\log\beta\ \vec{\Phi}^{\ep}_{x^0}+\beta\ \ep^{-1}\nabla\chi(\ep^{-1}(x-x^0))\cdot(x-x^0)\rg]\quad.
\]
Observe that
\[
\|\ep^{-1}\nabla\chi(\ep^{-1}(x-x^0))\cdot(x-x^0)\|_{L^\infty(B_\ep(x^0))}\le\|\nabla\chi\|_\infty
\]
This implies in particular using the explicit expression of $\beta$ given by (\ref{xsIII.37b})
\[
|\p_\ep\log\beta|\le |t|\ \|\nabla\chi\|_\infty [1+|t|]\ \lf[1+\sum_{j=1}^2|\al_j|\ \|a_j\|_\infty\rg]
\]
and using the explicit expression of $d{\mathcal C}_{\vec{\Phi}_{x^0}^\ep}$ we obtain that
\[
\begin{array}{l}
\ds|\p_\ep\Gamma(t,\ep,\al_1,\al_2)|\le\, C\  \|e^{2\la}\|_\infty\,\ep\ \|\nabla\chi\|_\infty [1+|t|]^2\ \lf[1+\sum_{j=1}^2|\al_j|\ \|a_j\|_\infty\rg]\lf[\int_{B_\ep(x^0)}|{\mathbb I}^0_{\vec{\Phi}}|^2\rg]^{1/2}\end{array}
\]
thus $\p_\ep\Gamma(t,\ep,\al_1,\al_2)$ extends continuously by $0$ at $(t,0,\al_1,\al_2)$ and the claim is proved.

Observe that we have
\be
\label{xsIII.40}
\p_{\al_j}\Gamma(0,0,0,0)=d{\mathcal C}_{\vec{\Phi}}\cdot\vec{a}_j
\ee
By assumption on $(\vec{a}_j)_{j=1,2}$ the pair $(\p_{\al_j}\Gamma(0,0,0,0))_{j=1,2}$ generates the tangent space to the Teichm\"uller space ${\mathcal T}_\Sigma$
at ${\mathcal C}(\vec{\Phi})$. We can then apply the local inversion theorem and deduce the existence of a neighborhood ${\mathfrak U}$
of $(0,0)$ and the existence of two $C^1$ functions $\al_1(t,\ep)$ and $\al_2(t,\ep)$ defined on ${\mathfrak U}$ such that there exists a neighborhood ${\mathcal U}$ of $(0,0,0,0)$ included in ${\mathcal U}_{\Phi}$ for which
\[
{\mathcal U}\cap {\mathcal C}^{-1}({\mathcal C}(\vec{\Phi}))=\{(t,\ep,\al_1(t,\ep),\al_2(t,\ep))\}\quad.
\]
It is clear that 
\[
\al_j(0,\ep)=\al_j^0(\ep)\quad.
\]
Denote $\beta^\ep(x,t):=\beta(t,\ep,\al_1(t,\ep),\al_2(t,\ep))$ and 
\[
\vec{\Phi}^\ep_{x^0}(t):=\beta^\ep(x,t)\,\lf[\vec{\Phi}+ t\, \chi^\ep_{x^0}\ \vec{n}^0+t \sum_{j=1}^2\al_j(t,\ep)\ \vec{a}_j\rg]
\]
We have for any $(t,\ep)$ in ${\mathfrak U}$
\be
\label{xsIII.42s}
0\equiv d{\mathcal C}_{\vec{\Phi}_{x^0}^\ep(t)}\cdot\frac{\p \vec{\Phi}_{x^0}^\ep(t)}{\p t}=\beta^\ep_{x^0}(x,t)\ d{\mathcal C}_{\vec{\Phi}_{x^0}^\ep(t)}\cdot
\lf[\chi_{x^0}^\ep\ \vec{n}^0
+\sum_{j=1}^2\ \frac{d(t\,\al_j(t,\ep))}{dt}\ \vec{a}_j\rg]
\ee
where we have used the fact that
\[
d{\mathcal C}_{\vec{\Phi}_{x^0}^\ep}\vec{\Phi}_{x^0}^\ep=\sum_{j=1}^2Q^j\ \int_{\Sigma}\vec{\Phi}_{x^0}^\ep\cdot <\ov{q}^j,\vec{\mathbb I}^0_{\vec{\Phi}}>_{g_{\vec{\Phi}}}\ dvol_{g_{\vec{\Phi}_{x^0}^\ep}}=0
\]
Hence
\be
\label{xsIII.34a-z11}
\sum_{j=1}^2\frac{d(t\,\al_j(t,\ep))}{dt}\    d{\mathcal C}_{\vec{\Phi}_{x^0}^\ep(t)}\cdot\vec{a}_j=d{\mathcal C}_{\vec{\Phi}_{x^0}^\ep(t)}\cdot
\chi_{x^0}^\ep\,\vec{n}^0
\ee
For any $\vec{w}\in L^2$ the map
\[
d{\mathcal C}_{\vec{\Phi}_{x^0}^\ep(t)}\cdot\vec{w}=\sum_{j=1}^2Q^j\ \int_{\Sigma}\vec{w}\cdot <\ov{q}^j,\vec{\mathbb I}^0_{\vec{\Phi}_{x^0}^\ep}>_{g_{\vec{\Phi}_{x^0}^\ep}}\ dvol_{g_{\vec{\Phi}_{x^0}^\ep}}
\]
Hence we have the existence of two maps $\vec{F}^j(x,p,q)$ for $x\in\Sigma$, $p\in {\R}^2\otimes{\R}^4$ and $q\in {\R}^2\otimes{\R}^2\otimes{\R}^4$ which is algebraic in $p$ and linear in $q$ and such that for any $0<\delta<1$and for $p=p_1\otimes p_2$ satisfying $\delta |p|^2\le p_1\wedge p_2\le \delta^{-1} |p|^2$
\[
|F^j(x,p,q)|\le C_\delta\ |p|^{-2}\, |q| 
\]  
where $C_\delta$ only depends on $\delta$ and such that
\[
d{\mathcal C}_{\vec{\Phi}_{x^0}^\ep(t)}\cdot\vec{w}=\sum_{j=1}^2Q^j\ \int_{\Sigma}\vec{w}\cdot \vec{F}^j(x,\nabla\vec{\Phi}_{x^0}^\ep(t),\nabla^2\vec{\Phi}_{x^0}^\ep(t))\ dvol_{g_0}
\]
Combining this fact and the explicit expression of $\vec{\Phi}_{x^0}^\ep(t)$ we deduce that, for $\ep$ positive,
\[
t\, \longrightarrow d{\mathcal C}_{\vec{\Phi}_{x^0}^\ep(t)}\cdot\vec{w}
\]
is a $C^1$ function. Thus the following functions are $C^1$ with respect to $t$ as long as $\ep>0$ :
\[
d{\mathcal C}_{\vec{\Phi}_{x^0}^\ep(t)}\cdot\vec{a}_j\quad\mbox{ and }\quad d{\mathcal C}_{\vec{\Phi}_{x^0}^\ep(t)}\cdot
\chi_{x^0}^\ep\,\vec{n}^0
\]
Since $d{\mathcal C}_{\vec{\Phi}_{x^0}^\ep(0)}\cdot\vec{a}_j=d{\mathcal C}_{\vec{\Phi}}\cdot\vec{a}_j$ forms a basis of the
tangent space to ${\mathcal T}_{T^2}$ at ${\mathcal C}(\vec{\Phi})$, for $t$ small enough $(d{\mathcal C}_{\vec{\Phi}_{x^0}^\ep(t)}\cdot\vec{a}_j)$ is a $C^1$ map into the space of basis of ${\mathcal T}_{T^2}$ at ${\mathcal C}(\vec{\Phi})$. The functions 
$\frac{d(t\,\al_j(t,\ep))}{dt}$ are then the coordinates of the $C^1$ map $d{\mathcal C}_{\vec{\Phi}_{x^0}^\ep(t)}\cdot
\chi_{x^0}^\ep\,\vec{n}^0$ in the $C^1$ frame $(d{\mathcal C}_{\vec{\Phi}_{x^0}^\ep(t)}\cdot\vec{a}_j)_{j=1,2}$. Thus 
$\frac{d(t\,\al_j(t,\ep))}{dt}$ are $C^1$ functions of $t$ in a neighborhood of $0$ for $\ep>0$.

\medskip

We shall establish now the estimates (\ref{xsIII.34a-z1}), (\ref{xsIII.34a-z19p}), (\ref{xsIII.34a-z22p}) and (\ref{xsIII.34a-z30p}). We take  $t>0$ to simplify the notations. Using the definition (\ref{xsIII.37b})
 denoting by ${\mathbf 1}^\ep_{x^0}$ the characteristic function of the geodesic ball $B_\ep(x^0)$ for the metric $g_0$, we have
\be
\label{xsIII.34a-z3}
\begin{array}{l}
\ds\lf|\frac{d\beta^\ep_{x^0}}{dt}(x,t)\rg|\le\, C\, \ep\,{\mathbf 1}^\ep_{x^0}\lf[ |\vec{n}(x)-\vec{n}^0|+\,t\,\ep+\sum_{j=1}^2|t\, \al_j|+t\, \sum_{j=1}^2\lf| \frac{d[t\,\al_j]}{dt}\rg|\rg]\\[5mm]
\ds\quad\quad\quad\quad\quad+C\,\sum_{j=1}^2|t\, \al_j| \lf|\frac{d[t\,\al_j]}{dt}\rg|
\end{array}
\ee
and
\be
\label{xsIII.34a-z4}
\begin{array}{l}
\ds\lf|\frac{d^2\beta^\ep_{x^0}}{dt^2}(x,t)\rg|\le\, C\, \ep\,{\mathbf 1}^\ep_{x^0}\lf[\,\ep+ \sum_{j=1}^2\lf| \frac{d[t\,\al_j]}{dt}\rg|+t\,\sum_{j=1}^2\lf| 
\frac{d^2[t\,\al_j]}{dt^2}\rg|\rg]\\[5mm]
\ds\quad\quad\quad\quad\quad+C\,\sum_{j=1}^2|t\, \al_j|\, \lf|\frac{d^2[t\,\al_j]}{dt^2}\rg|+C\,\sum_{j=1}^2\lf| \frac{d[t\,\al_j]}{dt}\rg|^2+C\ \lf|\frac{d\beta^\ep}{dt}(x,t)\rg|^2
\end{array}
\ee
We have moreover
\be
\label{xsIII.34a-z5}
\begin{array}{l}
\ds|\nabla\beta^\ep_{x^0}|\le \, C_{\vec{\Phi}}\, \,{\mathbf 1}^\ep_{x^0}\lf[ \ep\,t\, e^{\la}|\vec{n}(x)-\vec{n}^0| + t\, |\vec{n}(x)-\vec{n}^0| +\ep\, t^2+t\,
\sum_{j=1}^2|t\, \al_j|\rg]\\[5mm]
\ds\quad\quad\quad\quad\quad+C_{\vec{\Phi}}\,\sum_{j=1}^2|t\, \al_j|^2 
\end{array}
\ee
and
\be
\label{xsIII.34a-z6}
\begin{array}{l}
\ds|\nabla^2\beta^\ep_{x^0}|\le \, C_{\vec{\Phi}}\,{\mathbf 1}^\ep_{x^0}\lf[   t\ \ep^{-1}\  |\vec{n}(x)-\vec{n}^0|+ t\,\ep\ |\nabla^2\vec{\Phi}|+ t^2 +t\,
\ep^{-1}\sum_{j=1}^2|t\, \al_j|\rg]\\[5mm]
\ds\quad\quad\quad\quad\quad+C_{\vec{\Phi}}\,\sum_{j=1}^2|t\, \al_j|^2\,[1+|\nabla^2a_j|]+C\,|\nabla\beta^\ep_{x^0}|^2
\end{array}
\ee
We have also
\be
\label{xsIII.34a-z7}
\begin{array}{l}
\ds\lf|\nabla\frac{d\beta^\ep_{x^0}}{dt}\rg|\le \, C_{\vec{\Phi}}\,\,{\mathbf 1}^\ep_{x^0}\lf[ [\,|\vec{n}(x)-\vec{n}^0| +\ep\, t+
\sum_{j=1}^2|t\, \al_j|+t\,\sum_{j=1}^2\lf|\frac{d}{dt}[t\, \al_j]\rg|\rg]\\[5mm]
\ds\quad\quad\quad\quad\quad+C_{\vec{\Phi}}\,\sum_{j=1}^2 |t\, \al_j|\, \lf|\frac{d[t\,\al_j]}{dt}\rg|+|\nabla\beta^\ep_{x^0}|\,\lf|\frac{d\beta^\ep_{x^0}}{dt}(x,t)\rg|
\end{array}
\ee
and
\be
\label{xsIII.34a-z8}
\begin{array}{l}
\ds\lf|\nabla\frac{d^2\beta^\ep_{x^0}}{dt^2}\rg|\le \, C_{\vec{\Phi}}\,\,{\mathbf 1}^\ep_{x^0}\lf [ \ep+ \sum_{j=1}^2\lf| \frac{d[t\,\al_j]}{dt}\rg| +   t\,\sum_{j=1}^2\lf| 
\frac{d^2[t\,\al_j]}{dt^2}\rg|  \rg]    \\[5mm]
\ds\quad\quad\quad\quad\quad+C_{\vec{\Phi}}\,\sum_{j=1}^2|t\, \al_j|\,\lf| \frac{d^2[t\,\al_j]}{dt^2}\rg|+C\,\sum_{j=1}^2\lf| \frac{d[t\,\al_j]}{dt}\rg|^2\\[5mm]
\ds\quad\quad\quad\quad\quad+|\nabla\beta^\ep_{x^0}|\,\lf|\frac{d\beta^\ep_{x^0}}{dt}(x,t)\rg|^2+|\nabla\beta^\ep_{x^0}|\,\lf|\frac{d^2\beta^\ep_{x^0}}{dt^2}(x,t)\rg|+\lf|\nabla\frac{d\beta^\ep_{x^0}}{dt}\rg|\,\lf|\frac{d\beta^\ep_{x^0}}{dt}(x,t)\rg|
\end{array}
\ee
Finally, regarding this time the second derivatives of $\beta_{x^0}^\ep$ we have the following pointwise estimates
\be
\label{xsIII.34a-z9}
\begin{array}{l}
\ds\lf|\nabla^2\frac{d\beta^\ep_{x^0}}{dt}\rg|\le \, C_{\vec{\Phi}}\,\,{\mathbf 1}^\ep_{x^0}\lf[ \ep^{-1}\,|\vec{n}(x)-\vec{n}^0| +\ep\, |\nabla^2\vec{\Phi}|+ t+\ep^{-1}
\sum_{j=1}^2|t\, \al_j|+t\,\sum_{j=1}^2\lf|\frac{d}{dt}[t\, \al_j]\rg|\rg]\\[5mm]
\ds\quad\quad\quad\quad\quad+C_{\vec{\Phi}}\,\sum_{j=1}^2|t\, \al_j|\, \lf|\frac{d[t\,\al_j]}{dt}\rg|\,[1+|\nabla^2a_j|]\\[5mm]
\ds\quad\quad\quad\quad\quad+|\nabla\beta^\ep_{x^0}|^2\,\lf|\frac{d\beta^\ep_{x^0}}{dt}(x,t)\rg|+|\nabla^2\beta^\ep_{x^0}|\,\lf|\frac{d\beta^\ep_{x^0}}{dt}(x,t)\rg|+|\nabla\beta^\ep_{x^0}|\,\lf|\nabla\frac{d\beta^\ep_{x^0}}{dt}(x,t)\rg|\\[5mm]
\end{array}
\ee
and
\be
\label{xsIII.34a-z10}
\begin{array}{l}
\ds\lf|\nabla^2\frac{d^2\beta^\ep_{x^0}}{dt^2}\rg|\le \, C_{\vec{\Phi}}\,\,{\mathbf 1}^\ep_{x^0}\lf[ 1+\ep^{-1}\sum_{j=1}^2\lf|\frac{d}{dt}[t\, \al_j]\rg|+\ep^{-1}t\,\sum_{j=1}^2\lf|\frac{d^2}{dt^2}[t\, \al_j]\rg|\rg]\\[5mm]
 \ds\quad\quad+C_{\vec{\Phi}}\,\sum_{j=1}^2[|t\, \al_j|\, \lf|\frac{d^2[t\,\al_j]}{dt^2}\rg|+C_{\vec{\Phi}}\,\sum_{j=1}^2[1+|\nabla^2a_j|]\lf|\frac{d[t\,\al_j]}{dt}\rg|^2\\[5mm]
\ds\quad\quad+|\nabla^2\beta^\ep_{x^0}|\,\lf|\frac{d^2\beta^\ep_{x^0}}{dt^2}(x,t)\rg|+\lf|\nabla\frac{d\beta^\ep_{x^0}}{dt}(x,t)\rg|\,\lf|\nabla\frac{d\beta^\ep_{x^0}}{dt}(x,t)\rg|+|\nabla\beta^\ep_{x^0}|\,\lf|\nabla\frac{d^2\beta^\ep_{x^0}}{dt^2}(x,t)\rg|\\[5mm]
\ds\quad\quad+\lf|\nabla^2\frac{d\beta^\ep_{x^0}}{dt}(x,t)\rg|\,\lf|\frac{d\beta^\ep_{x^0}}{dt}(x,t)\rg|+|\nabla\beta^\ep_{x^0}|^2\,\lf|\frac{d\beta^\ep_{x^0}}{dt}(x,t)\rg|^2+|\nabla\beta^\ep_{x^0}|^2\,\lf|\frac{d^2\beta^\ep_{x^0}}{dt^2}(x,t)\rg|\\[5mm]
\ds\quad\quad+|\nabla^2\beta^\ep_{x^0}|\,\lf|\frac{d\beta^\ep_{x^0}}{dt}(x,t)\rg|^2\\[5mm]
\end{array}
\ee

We shall now estimate successively $(t\,\al_j(t,\ep))$, $d/dt(t\,\al_j(t,\ep))$ and $d^2/dt^2(t\,\al_j(t,\ep))$. Since $d{\mathcal C}_{\vec{\Phi}_{x^0}^\ep(0)}\cdot\vec{a}_j=d{\mathcal C}_{\vec{\Phi}}\cdot\vec{a}_j$ forms a basis of the
tangent space to ${\mathcal T}_{T^2}$ at ${\mathcal C}(\vec{\Phi})$, using  (\ref{xsIII.38t}) and the pointwize inequality
\[
\lf|{\mathbb I}_{\vec{\Phi}_{x^0}^\ep(t)}\rg|\le \lf|{\mathbb I}_{\vec{\Phi}}\rg| + C_{\Phi}\, \lf[t\,\ep^{-1}+ t\sum_{k=1}^2|\al_k(t)|\ [|\nabla^2a_k|+1]\rg]
\]
  we have
\be
\label{xsIII.34a-z12}
\begin{array}{l}
\ds\lf|\frac{d}{dt}(t\,\al_j(t,\ep))\rg|\le C\,\lf|d{\mathcal C}_{\vec{\Phi}^\ep_{x^0}(t)}\cdot \vec{n}^0\, \chi_{x^0}^\ep\rg|=O_{x^0,\vec{\Phi}}(\ep^3+t\,\ep^2)
%\ds\quad\le C\, \lf|[d{\mathcal C}_{\vec{\Phi}^\ep_{x^0}(t)}-d{\mathcal C}_{\vec{\Phi}}]\cdot \vec{n}^0\, \chi_{x^0}^\ep\rg|+C\,\lf|d{\mathcal C}_{\vec{\Phi}}\cdot \vec{n}^0\, \chi_{x^0}^\ep\rg|
\end{array}
\ee
%Using (\ref{xsIII.38t}) we have then
%\be
%\label{xsIII.34a-z13}
%\begin{array}{l}
%\ds\lf|\frac{d}{dt}(t\,\al_j(t,\ep))\rg|\le C\,\ep\,\int_{B_\ep(x^0)}|\nabla^2\vec{\Phi}_{x^0}^\ep(t)- \nabla^2\vec{\Phi}|\ e^{-2\la}\ dx \\[5mm]
%\ds\quad\quad+C\,\ep\,\int_{B_\ep(x^0)}| \nabla^2\vec{\Phi}|\ |\nabla\vec{\Phi}_{x^0}^\ep(t)- \nabla\vec{\Phi}| e^{-3\la}\ dx      +O_{x^0,\vec{\Phi}}(\ep^3)
%\end{array}
%\ee
%We have that
%\be
%\label{xsIII.34a-z14}
%\begin{array}{l}
%\ds|\nabla^2\vec{\Phi}_{x^0}^\ep(t)- \nabla^2\vec{\Phi}|\le C\, [|\nabla^2\beta_{x^0}^\ep|+ e^\la\,|\nabla \beta_{x^0}^\ep|+|\beta_{x^0}^\ep-1|\,|\nabla^2\vec{\Phi}|+t\,\ep^{-1}\,{\mathbf 1}^\ep_{x^0}]\\[5mm]
%\ds\quad\quad+C\, |t\,\al_j(t,\ep)|
%\end{array}
%\ee
%and
%\be
%\label{xsIII.34a-z15}
%|\nabla\vec{\Phi}_{x^0}^\ep(t)-\nabla\vec{\Phi}|\le C\, [|\nabla\beta_{x^0}^\ep|+ e^\la\,| \beta_{x^0}^\ep-1|+\,t\ {\mathbf 1}^\ep_{x^0}+ |t\,\al_j(t,\ep)|]
%\ee
%Combining (\ref{xsIII.34a-z5}), (\ref{xsIII.34a-z6}), (\ref{xsIII.34a-z13}), (\ref{xsIII.34a-z14}) and (\ref{xsIII.34a-z15}), we obtain after some computations
%\be
%\label{xsIII.34a-z16}
%\lf|\frac{d}{dt}(t\,\al_j(t,\ep))\rg|\le C\, t\, \ep^2+O_{x^0,\vec{\Phi}}(\ep^2)\ (t\,\al_j(t,\ep))+O_{x^0,\vec{\Phi}}(\ep^3)\quad.
%\ee 
%From which we deduce
%\be
%\label{xsIII.34a-z17}
%\lf\{
%\begin{array}{l}
%\ds|t\,\al_j(t,\ep)|\le O_{x^0,\vec{\Phi}}(\ep^3)\ |t|+C\, t^2\ \ep^2\\[5mm]
%\ds \lf|\frac{d}{dt}(t\,\al_j(t,\ep))\rg|\le C\, |t|\, \ep^2+O_{x^0,\vec{\Phi}}(\ep^3)\quad .
%\end{array}
%\rg.
%\ee 
Inserting this estimate in (\ref{xsIII.34a-z3}) gives
\be
\label{xsIII.34a-z30}
\lf|\frac{d\beta^\ep_{x^0}}{dt}\rg|\le\, C_{\vec{\Phi}}\ \ep\ {\mathbf 1}^\ep_{x^0}\ \lf[|\vec{n}(x)-\vec{n}^0|+t\,\ep\rg]+\, t\ O_{x^0,\vec{\Phi}}((\ep^3+t\,\ep^2)^2)\quad.
\ee
in (\ref{xsIII.34a-z5}) in also gives
\be
\label{xsIII.34a-z31-a}
\lf|\nabla\beta^\ep_{x^0}\rg|\le\, C_{\vec{\Phi}}\  {\mathbf 1}^\ep_{x^0}\ \lf[t\,|\vec{n}(x)-\vec{n}^0|+t^2\,\ep\rg]+\, t^2\ O_{x^0,\vec{\Phi}}((\ep^3+t\,\ep^2)^2)\quad,
\ee
and 
\be
\label{xsIII.34a-z31}
\begin{array}{l}
\ds\lf|\nabla^2\beta^\ep_{x^0}\rg|\le\, C_{\vec{\Phi}}\  {\mathbf 1}^\ep_{x^0}\ \lf[t\,\ep^{-1}\,|\vec{n}(x)-\vec{n}^0|+t\,\ep\,|\nabla^2\vec{\Phi}|+t^2\,\rg]\\[5mm]
\ds\quad\quad\quad+\, t^2\ O_{x^0,\vec{\Phi}}((\ep^3+t\,\ep^2)^2)\,[1+|\nabla^2a|]\quad.
\end{array}
\ee
Inserting  (\ref{xsIII.34a-z3}) and the previous two estimates in (\ref{xsIII.34a-z7}) give
\be
\label{xsIII.34a-z32}
\lf|\nabla\frac{d\beta^\ep_{x^0}}{dt}\rg|\le\, C_{\vec{\Phi}}\  {\mathbf 1}^\ep_{x^0}\ \lf[\,|\vec{n}(x)-\vec{n}^0|+t\,\ep\rg]+\, t\ O_{x^0,\vec{\Phi}}((\ep^3+t\,\ep^2)^2)\quad.
\ee
Finally  inserting  (\ref{xsIII.34a-z3})and the previous two estimates in (\ref{xsIII.34a-z9}) give
\be
\label{xsIII.34a-z33}
\begin{array}{l}
\ds\lf|\nabla^2\frac{d\beta^\ep_{x^0}}{dt}\rg|\le\, C_{\vec{\Phi}}\  {\mathbf 1}^\ep_{x^0}\ \lf[\,\ep^{-1}\,|\vec{n}(x)-\vec{n}^0|+ \ep \, |\nabla^2\vec{\Phi}|+t\,\rg]\\[5mm]
\ds\quad+\, t\ O_{x^0,\vec{\Phi}}((\ep^3+t\,\ep^2)^2)[1+|\nabla^2a|]\quad.
\end{array}
\ee
We shall now estimate $d(g_{\vec{\Phi}_{x^0}^\ep(t)})/dt$ as well as $d({\mathbb I}^0_{\vec{\Phi}_{x^0}^\ep(t)})/dt$. We first have
\be
\label{xsIII.34a-z18}
\begin{array}{l}
\ds\lf|\frac{d}{dt}\nabla\vec{\Phi}_{x^0}^\ep(t)\rg|\le \lf|\frac{d}{dt}\nabla\beta_{x^0}^\ep(t)\rg|+C_{\vec{\Phi}}\,\lf|\frac{d}{dt}\beta_{x^0}^\ep(t)\rg|\,\lf[1+ |t\,\al_j(t,\ep)|\rg]\\[5mm]
\ds\quad\quad\quad\quad+C_{\vec{\Phi}}\ \lf[  {\mathbf 1}^\ep_{x^0}+ \lf|\frac{d}{dt}(t\,\al_j(t,\ep))\rg|+|\nabla\beta^\ep_{x^0}|\ \lf[  \ep\ {\mathbf 1}^\ep_{x^0}+ \lf|\frac{d}{dt}(t\,\al_j(t,\ep))\rg| \rg]\rg]\quad .
\end{array}
\ee
From which we deduce using the previous estimates
\be
\label{xsIII.34a-z19}
\ds\lf|\frac{d}{dt}\nabla\vec{\Phi}_{x^0}^\ep(t)\rg|\le C_{\vec{\Phi}}\ \lf[ {\mathbf 1}^\ep_{x^0}+O_{x^0,\vec{\Phi}}(\ep^3+t\,\ep^2)   \rg]\quad .
\ee
Regarding the metric $g_{\vec{\Phi}_{x^0}^\ep(t)}$ we have
\[
g_{\vec{\Phi}_{x^0}^\ep(t)}=\nabla\vec{\Phi}_{x^0}^\ep(t)\dot{\otimes}\nabla\vec{\Phi}_{x^0}^\ep(t)
\]
So we deduce
\[
\lf|\frac{dg_{ij}}{dt}-\beta^2\frac{d}{dt}\lf(\p_{x_i}(\beta^{-1}\, \vec{\Phi}_{x^0}^\ep(t)) \cdot\p_{x_j}(\beta^{-1}\, \vec{\Phi}_{x^0}^\ep(t))\rg)\rg|\le C_{\vec{\Phi}}\lf[ \lf|\nabla\frac{d\beta}{dt}\rg|+  \lf|\frac{d\beta}{dt}\rg|
+|\nabla\beta|\ \lf|\frac{d}{dt}\nabla\vec{\Phi}_{x^0}^\ep(t)\rg|\rg]
\]
Observe that we have
\be
\label{xsIII.34a-z19-a}
\begin{array}{l}
\ds\lf(\p_{x_i}(\beta^{-1}\, \vec{\Phi}_{x^0}^\ep(t)) \cdot\p_{x_j}(\beta^{-1}\, \vec{\Phi}_{x^0}^\ep(t))\rg)=\delta_{ij}\, e^{2\la}\\[5mm]
\ds\quad+t\ [\vec{n}^0-\vec{n}]\ \lf[\p_{x_i}\chi(\ep^{-1}(x-x^0))\, \p_{x_j}\vec{\Phi} +\p_{x_j}\chi(\ep^{-1}(x-x^0))\, \p_{x_i}\vec{\Phi} \rg]\\[5mm]
\ds\quad+t\,\sum_{k=1}^2\al_k(t)\ [\p_{x_j}\vec{a}_k\cdot\p_{x_i}\vec{\Phi}+\p_{x_i}\vec{a}_k\cdot\p_{x_j}\vec{\Phi}]\\[5mm]
\ds\quad+t^2\,\sum_{k=1}^2\lf[\p_{x_i}\chi(\ep^{-1}(x-x^0))\ \vec{n}^0\cdot \p_{x_j}\vec{a}_k+\p_{x_j}\chi(\ep^{-1}(x-x^0))\ \vec{n}^0\cdot \p_{x_i}\vec{a}_k\rg]\\[5mm]
\ds\quad+t^2\p_{x_i}\chi(\ep^{-1}(x-x^0))\p_{x_j}\chi(\ep^{-1}(x-x^0))+t^2\sum_{k=1}^2\al_k(t)\ \p_{x_i}\vec{a}_k\cdot\sum_{k=1}^2\al_k(t)\ \p_{x_j}\vec{a}_k
\end{array}
\ee
We deduce from the previous identities
\be
\label{xsIII.34a-z20}
\lf|\frac{d g_{\vec{\Phi}_{x^0}^\ep(t)}}{dt}\rg|\le C_{\vec{\Phi}}\,\lf[  {\mathbf 1}^\ep_{x^0}(|\vec{n}-\vec{n}^0|+t)+ O_{x^0,\vec{\Phi}}(\ep^3+t\,\ep^2) \rg]
\ee
Now we have
\be
\label{xsIII.34a-z21}
\begin{array}{l}
\ds\lf|\frac{d}{dt}\nabla^2\vec{\Phi}_{x^0}^\ep(t)\rg|\le \lf|\frac{d}{dt}\nabla^2\beta_{x^0}^\ep(t)\rg|+C_{\vec{\Phi}}\,|\nabla\beta_{x^0}^\ep|\,\lf[{\mathbf 1}^\ep_{x^0}+  \lf|\frac{d}{dt}(t\,\al_j(t,\ep))\rg|\rg]\\[5mm]
\ds\quad+C_{\vec{\Phi}}\,\lf|\frac{d}{dt}\nabla\beta_{x^0}^\ep(t)\rg|\,+C_{\vec{\Phi}}\ \lf[ \ep^{-1}\, {\mathbf 1}^\ep_{x^0}+ \lf|\frac{d}{dt}(t\,\al_j(t,\ep))\rg|\ |\nabla^2 a_j|\rg]\\[5mm]
\ds\quad+ C_{\vec{\Phi}}\ |\nabla^2\beta^\ep_{x^0}|\ \lf[  \ep\ {\mathbf 1}^\ep_{x^0}+ \lf|\frac{d}{dt}(t\,\al_j(t,\ep))\rg| \rg]
\end{array}
\ee
From which we deduce using the previous estimates
\be
\label{xsIII.34a-z22}
\ds\lf|\frac{d}{dt}\nabla^2\vec{\Phi}_{x^0}^\ep(t)\rg|\le C_{\vec{\Phi}}\ {\mathbf 1}^\ep_{x^0}\ \lf[ \ep^{-1}+ \ep\, |\nabla^2\vec{\Phi}|\rg]+O_{x^0,\vec{\Phi}}(\ep^3+t\,\ep^2)\, \lf[1+\sum_{j=1}^2|\nabla^2 \vec{a}_j|\rg]  \quad .
\ee
This gives then
\be
\label{xsIII.34a-z23}
\ds\lf|\frac{d}{dt}{\mathbb I}_{\vec{\Phi}_{x^0}^\ep(t)}\rg|\le C_{\vec{\Phi}}\ {\mathbf 1}^\ep_{x^0}\ \lf[ \ep^{-1}+ \ep\, |\nabla^2\vec{\Phi}|\rg]+O_{x^0,\vec{\Phi}}(\ep^3+t\,\ep^2)\,\lf[1+ \sum_{j=1}^2|\nabla^2 \vec{a}_j|\rg]  \quad .
\ee
Taking now the $t$ derivative of the identity (\ref{xsIII.34a-z11})
\be
\label{xsIII.34a-z24}
\begin{array}{l}
\ds\sum_{j=1}^2\frac{d^2(t\,\al_j(t,\ep))}{dt^2}\    d{\mathcal C}_{\vec{\Phi}_{x^0}^\ep(t)}\cdot\vec{a}_j=-\sum_{j=1}^2\frac{d(t\,\al_j(t,\ep))}{dt}\    \frac{d\,\lf(d{\mathcal C}_{\vec{\Phi}_{x^0}^\ep(t)}\cdot\vec{a}_j\rg)}{dt}\\[5mm]
\ds\quad\quad\quad\quad\quad\quad\quad\quad\quad\quad\quad+\frac{d \lf(d{\mathcal C}_{\vec{\Phi}_{x^0}^\ep(t)}\cdot
\chi_{x^0}^\ep\,\vec{n}^0\rg)}{dt}
\end{array}
\ee
Since again $d{\mathcal C}_{\vec{\Phi}_{x^0}^\ep(0)}\cdot\vec{a}_j=d{\mathcal C}_{\vec{\Phi}}\cdot\vec{a}_j$ forms a basis of the
tangent space to ${\mathcal T}_{T^2}$ at ${\mathcal C}(\vec{\Phi})$ we deduce that for $j=1,2$
\be
\label{xsIII.34a-z25}
\begin{array}{l}
\ds\lf|\frac{d^2(t\,\al_j(t,\ep))}{dt^2}\rg|\le\sum_{j=1}^2\lf|\frac{d(t\,\al_j(t,\ep))}{dt}\rg|\    \lf|\frac{d\,\lf(d{\mathcal C}_{\vec{\Phi}_{x^0}^\ep(t)}\cdot\vec{a}_j\rg)}{dt}\rg|\\[5mm]
\ds\quad\quad\quad\quad\quad\quad\quad\quad\quad\quad\quad+\lf|\frac{d \lf(d{\mathcal C}_{\vec{\Phi}_{x^0}^\ep(t)}\cdot
\chi_{x^0}^\ep\,\vec{n}^0\rg)}{dt}\rg|
\end{array}
\ee
Using (\ref{xsIII.34a-z20}) and (\ref{xsIII.34a-z23}) we obtain in one hand
\be
\label{xsIII.34a-z26}
\lf|\frac{d\,\lf(d{\mathcal C}_{\vec{\Phi}_{x^0}^\ep(t)}\cdot\vec{a}_j\rg)}{dt}\rg|\le O_{x^0,\vec{\Phi}}(\ep)\ 
\ee
and in the other hand
\be
\label{xsIII.34a-z27}
\lf|\frac{d\,\lf(d{\mathcal C}_{\vec{\Phi}_{x^0}^\ep(t)}\cdot   \chi_{x^0}^\ep\,\vec{n}^0   \rg)}{dt}\rg|\le O_{x^0,\vec{\Phi}}(\ep^2)
\ee
Hence combining (\ref{xsIII.34a-z25}), (\ref{xsIII.34a-z26}) and (\ref{xsIII.34a-z27}) we obtain
\be
\label{xsIII.34a-z28}
\lf|\frac{d^2(t\,\al_j(t,\ep))}{dt^2}\rg|\le O_{x^0,\vec{\Phi}}(\ep^2)
\ee
We bound now $d^2(\nabla^2\vec{\Phi}_{x^0}^\ep)/dt^2$. We have
\be
\label{xsIII.34a-z29}
\begin{array}{l}
\ds\lf|\frac{d^2}{dt^2}\nabla^2\vec{\Phi}_{x^0}^\ep(t)\rg|\le \lf|\frac{d^2}{dt^2}\nabla^2\beta_{x^0}^\ep(t)\rg|+\lf|\frac{d}{dt}\nabla^2\beta_{x^0}^\ep(t)\rg|\,\lf[|\nabla \vec{\Phi}|+\, C\, t\ {\mathbf 1}^\ep_{x^0} +C\, |t\,\al_j(t,\ep)|\rg]\\[5mm]
\ds\quad\quad\quad\quad+|\nabla^2\beta^\ep_{x^0}|\ \lf[   \lf|\frac{d^2}{dt^2}(t\,\al_j(t,\ep))\rg| \rg]+ \lf|\frac{d^2}{dt^2}\nabla\beta_{x^0}^\ep(t)\rg|\ \lf[|\nabla \vec{\Phi}|+\, C\, t\ {\mathbf 1}^\ep_{x^0} +C\, |t\,\al_j(t,\ep)|\rg]\\[5mm]
\ds\quad\quad\quad\quad+ \lf|\frac{d}{dt}\nabla\beta_{x^0}^\ep(t)\rg|\ \lf[ C\, {\mathbf 1}^\ep_{x^0} +C\, \lf|\frac{d(t\,\al_j(t,\ep))}{dt}\rg|\ \rg] +C\, |\nabla\beta_{x^0}^\ep(t)|\ \lf|\frac{d^2(t\,\al_j(t,\ep))}{dt^2}\rg|\\[5mm]
\ds\quad\quad\quad\quad+\lf|\frac{d^2}{d t^2}\beta_{x^0}^\ep(t)\rg|\ \lf[   |\nabla \vec{\Phi}|+\, C\, t\ \ep^{-1}\,{\mathbf 1}^\ep_{x^0} +C\, |t\,\al_j(t,\ep)\rg]\\[5mm]
\ds\quad\quad\quad\quad+\lf|\frac{d}{dt}\beta_{x^0}^\ep(t)\rg|\ \lf[C\, \ep^{-1}\,{\mathbf 1}^\ep_{x^0}+C\, \lf|\frac{d(t\,\al_j(t,\ep))}{dt}\rg|\rg]+C\ \lf|\frac{d^2}{dt^2}(t\,\al_j(t,\ep))\rg|
\end{array}
\ee
Inserting the previous estimates in this inequality gives, after a lengthy computation,
\be
\label{xsIII.34a-z30-b}
\ds\lf|\frac{d^2}{dt^2}\nabla^2\vec{\Phi}_{x^0}^\ep(t)\rg|\le\ C_{\vec{\Phi},x^0}\   {\mathbf 1}^\ep_{x^0}+  [1+|\nabla^2\vec{\Phi}|]\ \, O_{x^0,\vec{\Phi}}(\ep^2+t\,\ep)  \quad .
\ee
which concludes the proof of lemma~\ref{lm-xsIII.1a}.
\hfill $\Box$

\end{document}